\documentclass[11pt,reqno,a4paper]{amsart}
\usepackage[latin1]{inputenc}
\usepackage{amsmath}
\usepackage{amsfonts}
\usepackage{amssymb}
\usepackage{graphicx}
\usepackage{bbm}
\usepackage{float}
\usepackage[left=2.0cm, right=2.0cm, top=2.3cm, bottom=2.5cm]{geometry}

\usepackage{color}
\definecolor{uibred}{RGB}{170, 0, 0}
\definecolor{uibblue}{RGB}{0, 84, 115}
\definecolor{uibgreen}{RGB}{119, 175, 0}
\definecolor{uiborange}{RGB}{217, 89, 0}
\definecolor{MyDarkGreen}{rgb}{0.00,0.41,0.24} % This is the color used for comments
\definecolor{hatter}{rgb}{0.85,0.85,0.85}
\definecolor{kek}{rgb}{0.00,0.00,1.00}
\definecolor{kekes}{rgb}{0.33,0.52,0.62}
\definecolor{lila}{rgb}{0.77,0.00,0.77}
\definecolor{szurke}{rgb}{0.74,0.74,0.76}

\usepackage{mathrsfs}

\usepackage{graphicx}%
\usepackage{color}
\usepackage{cite}
\usepackage{hyperref}
\hypersetup{
	colorlinks = true,%
	citecolor = blue,%[rgb]{0.65,0.0,0.0},
	filecolor=red,%
	linkcolor = [rgb]{0.65,0.0,0.0},%
	anchorcolor = red,
	pagecolor = red,
	urlcolor= [rgb]{0.65,0.0,0.0}
	linktocpage=true,
	pdfpagelabels=true,
	bookmarksnumbered=true,
}

\numberwithin{equation}{section}

\newcommand{\R}{\mathbb{R}}

\newtheorem{proposition}{Proposition}[section]
\newtheorem{theorem}{Theorem}[section]
\newtheorem{lemma}{Lemma}[section]
\newtheorem{corollary}{Corollary}[section]

\newtheorem{remark}{Remark}[section]
\newtheorem{example}{Example}[section]

\title[Equality in Borell-Brascamp-Lieb inequalities]{Equality in Borell-Brascamp-Lieb inequalities on curved spaces
}

\author{Zolt\'an M. Balogh and Alexandru Krist\'aly}	

\thanks{Z. M. Balogh was
	supported by the Swiss National Science Foundation, Grant Nr. {200020\_165507}.  A. Krist\'aly  was supported by the STAR-UBB Advanced Fellowship-Intern (Project CNFIS-FDI-2016-0056).}

\begin{document}
%\linenumbers

		\begin{abstract} 
					By using optimal mass transportation and a quantitative H\"older inequality, we provide estimates for the   Borell-Brascamp-Lieb deficit on complete Riemannian  manifolds. Accordingly, equality cases  in Borell-Brascamp-Lieb inequalities (including Brunn-Minkowski and Pr\'ekopa-Leindler inequalities) are characterized in terms of the optimal transport map between suitable marginal probability measures. These results provide several  qualitative applications both in the flat and non-flat frameworks. In particular, by using Caffarelli's regularity result for the Monge-Amp\`ere equation, we {give a new proof} of Dubuc's characterization of the equality in Borell-Brascamp-Lieb inequalities in the Euclidean setting. When the  $n$-dimensional 
				Riemannian manifold has Ricci curvature  ${\rm Ric}(M)\geq (n-1)k$ for some $k\in \mathbb R$, it turns out that 
				equality in the Borell-Brascamp-Lieb inequality is expected only when a particular region of the manifold between the marginal supports has constant sectional curvature $k$. A precise characterization is provided for the equality  in the			
				Lott-Sturm-Villani-type distorted Brunn-Minkowski inequality on Riemannian manifolds. 
				Related results for (not necessarily reversible) Finsler manifolds are also presented.\\

\begin{center}
	\normalsize\textit{Dedicated to our friend, Professor Csaba Varga.}
\end{center}

	\end{abstract}
	
	%	\vspace*{-1.8cm}

%		By using optimal mass transportation and a quantitative H\"older inequality, we provide estimates for the   Borell-Brascamp-Lieb deficit on complete Riemannian  manifolds. As a consequence, equality cases  in Borell-Brascamp-Lieb inequalities (including Brunn-Minkowski and Pr\'ekopa-Leindler inequalities) are characterized in terms of the optimal transport map. In particular, when the  %$n$-dimensional 
%	Riemannian manifold has Ricci curvature bounded  below,  %${\rm Ric}(M)\geq (n-1)k$ for some $k\in \mathbb R$, 
%	equality in the Lott-Sturm-Villani-type distorted Brunn-Minkowsi inequality is expected only when the manifold is a space form. 
%	Some weak stability results are also established for log-Brunn-Minkowski inequalities. 
%	At the end, some related results for (not necessarily reversible) Finsler manifolds are also presented.

	\maketitle

	
	
%	\tableofcontents

	\section{Introduction}

	\subsection{Background and motivation}
	The Borell-Brascamp-Lieb inequality in the Euclidean space $\R^n$ states that for every
	fixed $s \in (0,1)$, $p \geq -\frac{1}{n}$ and integrable functions  $f, g, h: \R^n \to
	[0, \infty)$ which satisfy 
	\begin{eqnarray}\label{elso-BBL-feltetel}
	h((1-s)x + sy)
	\geq \mathcal M_s^p \left( f(x), g(y) \right) \quad \mbox{ for all } x, y \in \R^n,
	\end{eqnarray}
	one has
	\begin{eqnarray}\label{BBL-euklidesz}
		\int_{\R^n} h \geq \mathcal M_s^{\frac{p}{1 + np}} \left( \int_{\R^n} f,
		\int_{\R^n} g \right).
	\end{eqnarray}
	Here, for every $s\in (0,1)$, $p\in \mathbb R\cup
	\{\pm\infty\}$ and $a,b\geq 0$,  the $p$-mean is defined by
	$$\mathcal M_s^p(a,b)=\left\{\begin{array}{lll}
	\left( (1-s)a^p + s b^p \right)^{1/p} &\mbox{if} &  ab\neq 0, \\
	% u\geq 0 &\mbox{in} &   \Omega;\\
	0 &\mbox{if} &  ab=0,
	\end{array}\right.$$
	{with the  conventions $\mathcal M_s^{-\infty}(a,b)=\min\{a,b\}$; 
		$\mathcal M_s^{0}(a,b)=a^{1-s}b^s;$ and $\mathcal M_s^{+\infty}(a,b)=\max\{a,b\}$ if $ab\neq 0$ and $\mathcal M_s^{+\infty}(a,b)=0$ if $ab= 0$.} 
	
	 In terms of entropy, Borell-Brascamp-Lieb inequality implies that if a  Radon measure  $\mu$,   
has as density $\rho$  a $p$-concave function (i.e., $f=g=h=\rho$ satisfies (\ref{elso-BBL-feltetel})), then $\mu$ is a $q$-concave measure with the parameter $q = \frac{p}{1 + np}\in [-\infty,\frac{1}{n}]$. A sort of converse of the latter statement is given by Borell \cite{Borell}, characterizing the $q$-concave measures by means of $p$-concave functions. Further  contributions to this subject can be found  in Bobkov and Ledoux  \cite{Bobkov-Ledoux}, Brascamp and Lieb \cite{Brascamp-Lieb}. In particular, $0$-concave measures are characterized by log-concave density functions, while $\frac{1}{n}$-concave measures on convex sets are equal to the Lebesgue $\mathcal L^n$-measure up to a multiplicative constant.

	Another important consequence of the Borell-Brascamp-Lieb inequality (for $p=+\infty$ and indicator functions) is the usual Brunn-Minkowski inequality,  implying e.g. the isoperimetric inequality, which relates the $\mathcal L^n$-measure  of two measurable sets $A$ and $B$ in $\mathbb R^n$  with the (outer) $\mathcal L^n$-measure of their Minkowski sum $(1-s)A+sB=\{(1-s)x+sy:x\in A,y\in B\}$ as 
	\begin{equation}\label{BM-1}
	\mathcal L^n((1-s)A+sB)^\frac{1}{n}\geq (1-s) \mathcal L^n(A)^\frac{1}{n}+s\mathcal L^n(B)^\frac{1}{n}.
	\end{equation}
 An equivalent form of (\ref{BM-1}), coming also from the Borell-Brascamp-Lieb inequality (for $p=0$ and indicator functions), is the dimension-free-Brunn-Minkowski inequality, -- or the geometric form of the Pr\'ekopa-Leindler inequality, -- which states that 
\begin{equation}\label{log-BM}
\mathcal L^n((1-s)A+sB)\geq  \mathcal L^n(A)^{1-s}\mathcal L^n(B)^s.
\end{equation}

Characterizations of cases of {\it equality} and the problem of  {\it stability}  in the aforementioned inequalities (\ref{BBL-euklidesz})-(\ref{log-BM}) are still subjects for further investigation.  After the pioneering works by Brunn and Minkowski, it is well known for more than a century that 
	 equality in (\ref{BM-1}) holds if and only if the sets $A$ and $B$
	 are homothetic  convex  bodies  from  which  sets  of  measure  zero  have  been  removed; similarly,  equality in (\ref{log-BM}) holds if and only if the sets $A$ and $B$ are translated  convex  bodies up to a null measure set.  
	 The equality case in the generic Borell-Brascamp-Lieb inequality (\ref{BBL-euklidesz}) has been studied in the mid of seventies by Dubuc \cite{Dubuc} on $\mathbb R^n$ by using deep convexity and measure theoretical results together with a careful inductive argument w.r.t. the dimension of the space $\mathbb R^n$. Later on, Dancs and Uhrin \cite{DU, DU-2} obtained some qualitative Borell-Brascamp-Lieb inequalities on $\mathbb R$, providing also some higher-dimensional versions. A few years ago, Ball and B\"or\"oczky \cite{BB-1, BB-2} obtained stability results for the one-dimensional functional Pr\'ekopa-Leindler inequality with some extensions also to higher-dimensions. Very recently,  various stability results are established in $\mathbb R^n$ for the generic Borell-Brascamp-Lieb inequality by Ghilli and Salani \cite{GS}, Rossi \cite{Rossi-PhD}, Rossi and Salani \cite{Rossi-Salani, Rossi-Salani-AA}, for the Pr\'ekopa-Leindler inequality by Bucur and Fragal\`a \cite{BF},
	  and for the Brunn-Minkowski inequality  by Christ \cite{Christ}, Figalli and Jerison \cite{Figalli-Jerison, Figalli-Jerison-2, Figalli-Jerison-3} and Figalli, Maggi and Pratelli \cite{FMP-Inventiones, FMP-AIHP}. The common strategy in the aforementioned papers, up to the latter two papers, is the use of various arguments from convex analysis combined usually with some inductive step w.r.t. the dimension, by fully exploring the Euclidean character of the space. 
	  In \cite{FMP-Inventiones, FMP-AIHP}, quantitative Brunn-Minkowski inequalities are established by using optimal mass transportation arguments in $\mathbb R^n.$ Further results concerning equality and stability in the Brunn-Minkowski inequality in $\mathbb R^n$ can be found in Milman and Rotem \cite{MR} and Colesanti,  Livshyts and  Marsiglietti \cite{CLM}.

%	 It is well known  more than a hundred years that 
%	 equality in (\ref{BM-1}) holds if and only if the sets $A$ and $B$ are homothetic, while equality in (\ref{log-BM}) holds if and only if the sets $A$ and $B$ are translates.  The equality case in the generic Borell-Brascamp-Lieb inequality (\ref{BBL-euklidesz}) has been studied in the mid of seventies by Dubuc \cite{Dubuc} on $\mathbb R^n;$ his approach is based on convex analysis  and on  induction w.r.t. the dimension of the space $\mathbb R^n$. Later on, Dancs and Uhrin \cite{DU, DU-2} obtained some qualitative Borell-Brascamp-Lieb inequalities on $\mathbb R$, providing also some higher-dimensional versions by induction. A few years ago, Ball and B\"or\"oczky \cite{BB-1, BB-2} obtained stability results for the one-dimensional functional Pr\'ekopa-Leindler inequality with some extensions also to higher-dimensional cases. Very recently, further stability results are established for the generic Borell-Brascamp-Lieb inequality by Rossi and Salani \cite{Rossi-Salani} and for the Brunn-Minkowski inequality by Figalli and Jerison \cite{Figalli-Jerison, Figalli-Jerison-2}. The common strategy in the aforementioned papers is the use of induction w.r.t. the dimension,  exploring the Euclidean character of the space structure. 

	As far as we know, no equality/stability results are available for Borell-Brascamp-Lieb inequalities on  {\it curved spaces}. It is  clear that the arguments from the aforementioned papers (see \cite{BB-1}, \cite{BB-2}, \cite{BF}, \cite{Christ}, \cite{DU}, \cite{Dubuc}, \cite{Figalli-Jerison-2}, \cite{Figalli-Jerison},  and references therein) cannot be applied in such a nonlinear setting. 
	The starting point of our investigation is the celebrated work by  Cordero-Erausquin, McCann and Schmuckenschl{\"a}ger
	\cite{CMS} who established a Riemannian version of the Borell-Brascamp-Lieb inequality via optimal mass transportation culminating in a distorted Jacobian determinant inequality. 
	The Finslerian counterparts of the results from \cite{CMS} are provided by Ohta \cite{Ohta}.  
%	In turn, the latter works have their genesis in 
	We point out that the first optimal mass transportation approaches to geometric inequalities have  been  provided by Gromov in \cite{MS} (via the Knothe map) and 
%	to the 	
%	Pr\'ekopa-Leindler, Brunn-Minkowski and Borell-Brascamp-Lieb inequalities by 
McCann \cite{McCann_Adv_Math}, \cite[Appendix D]{McCann-PhD} (via the Brenier map).  It is worth
mentioning that Knothe \cite{Knothe} himself used his map to prove the  generalized
Brunn-Minkowski inequality \eqref{log-BM}.

%which encapsulates in particular the Pr\'ekopa-Leindler inequality. 

		The main purpose of our paper is {to characterize the equality in  Borell-Brascamp-Lieb inequalities on complete $n$-di\-men\-sional Riemannian/Finsler manifolds} for the whole spectrum of the parameter $p\geq -\frac{1}{n}$ by exploring a quantitative H\"older inequality and the theory of optimal mass transportation. 
%		Although our approach is more appropriate for characterizing equality cases, we furnish also
%		 some rigidity and weak stability results for geometric inequalities (e.g. for log-Brunn-Minkowski inequality).   
		
		In the sequel, we roughly present some of our achievements.

\subsection{Brief description of main results and consequences}	Let $(M,w)$
	be a  complete $n$-dimensional Riemannian manifold $(n\geq 2)$ with the induced distance function $d:M\times M\to [0,\infty)$;  without mentioning explicitly,   we assume throughout the whole paper that $(M,w)$ is connected.	For a fixed
	$s \in (0,1)$ and $(x,y) \in M \times M$ let
	%\begin{eqnarray}\label{Riemannian-intermediate}
	$$Z_s(x,y) = \{ z \in M : d(x,z) = s d(x,y),\
	d(z,y) = (1-s) d(x,y)\}$$
	be the set of $s$-intermediate points between $x$ and $y$, replacing the convex combination in (\ref{elso-BBL-feltetel}). 
	%\end{eqnarray}
	%and
	Since
	$(M,d)$ is complete, $Z_s(x,y)\neq \emptyset$ for
	every $x,y\in M.$  Accordingly,  the Minkowski interpolation set $$
	Z_s(A,B) = \bigcup\limits_{(x,y) \in A \times B} Z_s(x,y)
	$$ replaces the Minkowski sum of the nonempty sets $A,B \subset M$.

	Let  $s\in (0,1)$ and $p\geq -\frac{1}{n}$.  
	If $f,g,h:M\to [0,\infty)$ are three nonzero, compactly supported integrable functions, the natural Riemannian reformulation of (\ref{elso-BBL-feltetel})  reads as 
	\begin{eqnarray}\label{ConditionRescaledBBLWithWeights-vege} h(z)
	\geq \mathcal M^{p}_s
	\left(\frac{f(x)}{v_{1-s}(y,x)},\frac{g(y)}{v_s(x,y)} \right) \
	\  {\rm for\ all}\ (x,y)\in M\times M, z\in Z_s(x,y),
	\end{eqnarray}
	where $v_s$ is the volume distortion coefficient (see (\ref{vol-distortion}) for its precise definition).  
	Under the assumption (\ref{ConditionRescaledBBLWithWeights-vege}), the main result of Cordero-Erausquin, McCann and Schmuckenschl{\"a}ger
	\cite{CMS}  says that
	\begin{eqnarray}\label{BBL-eredeti-Riemann}
		\int_M h \geq \mathcal  M_s^{\frac{p}{1 + np}} \left( \int_M f,
		\int_M g \right),
	\end{eqnarray}
	where the integrals are considered w.r.t. the  Riemannian measure  $\textsf{m}$ on $(M,w)$.  This will be referred  throughout the paper
	as the Borell-Brascamp-Lieb inequality with exponent $p$.

	%$\delta_{\mathbb R^n,s}^p(f,g,h)\geq 0.$

	For simplicity of notation, let $\|\cdot\|_1$ be the $L^1$-norm of any integrable function on $M.$
	For the above functions $f,g$ and $h$, let us consider the {\it Borell-Brascamp-Lieb deficit} given by
	$$\delta_{M,s}^p(f,g,h)=\frac{\| h\|_1}{\mathcal  M_s^\frac{p}{1+pn}\left(\|f\|_1,\| g\|_1\right)}-1. $$

	We first provide an estimate for the  Borell-Brascamp-Lieb deficit on a  general Riemannian manifold that will be achieved by using optimal mass transportation and a quantitative H\"older inequality:

	\begin{theorem}\label{Theorem-Riemannian}  {\bf (Estimate of the Borell-Brascamp-Lieb deficit)} Let {\rm $(M,w)$} be a complete $n$-di\-men\-sional Riemannian manifold, $s\in (0,1),$ $p\geq -\frac{1}{n}$ and  $f,g,h:M\to [0,\infty)$ be three nonzero, compactly supported integrable functions  satisfying $(\ref{ConditionRescaledBBLWithWeights-vege}).$ Then 
		$$\displaystyle\delta_{M,s}^p(f,g,h)\geq \displaystyle\int_M \tilde f(x)G_s^{p,n}\left(\frac{f(x)}{v_{1-s}(\psi(x),x)},\frac{g(\psi(x))}{v_s(x,\psi(x))},\frac{1}{\|f\|_1},\frac{1}{\|g\|_1}\right){\rm d}\sf{m},
		$$
		where $\psi:M\to M$ is the unique optimal transport
		map from the measure $\mu=\tilde f {\rm d}\sf{m}$ to $\nu=\tilde g{\rm d}
		\sf{m}$ with densities $\tilde f=f/\|f\|_1$, $\tilde g=g/\|g\|_1$, 
		and $G_s^{p,n}\geq 0 $ is the gap-function given in Lemma $\ref{lemma-p-mean}.$
	\end{theorem}
	
%	Some remarks are in order. 
	
	The uniqueness of the optimal transport
map $\psi:M\to M$ from the probability measure $\mu=\tilde f {\rm d}\textsf{m}$ to $\nu=\tilde g{\rm d}
\textsf{m}$ is well known by McCann \cite{McCann} having the form $\psi(x)=\exp_x(-\nabla \varphi(x))$ for a.e. $x\in \operatorname{supp} f$ for some $d^2/2$-concave function $\varphi:M\to \mathbb R,$ where $\nabla$ denotes the Riemannian gradient.  Let $\psi_s:M\to M$ be the $s$-interpolant optimal transport map $\psi_s(x)=\exp_x(-s\nabla \varphi(x))$ for a.e. $x\in \operatorname{supp} f$, and Jac$(\psi_s)(x)$ its Jacobian in a.e. $x\in \operatorname{supp} f$. 

%, and $\mu_s=(\psi_s)_{\#}\mu$
%be the push-forward measure whose density is $\rho_s=\frac{\dd \mu_s}{ {\rm d}V_w}.$

By Theorem \ref{Theorem-Riemannian}  the equality in the Borell-Brascamp-Lieb inequality can be characterized by studying the properties of the gap-function $G_s^{p,n},$  leading us to the following result:

		\begin{theorem}\label{Theorem-Riemannian-egyenloseg} {\bf (Equality in Borell-Brascamp-Lieb inequality; $p>-\frac{1}{n}$)} 
%			Let $s\in (0,1)$ and $p>-\frac{1}{n}$. Under the same assumptions as in Theorem $\ref{Theorem-Riemannian},$ 
Let {\rm $(M,w)$} be a complete $n$-di\-men\-sional Riemannian manifold, $s\in (0,1),$ $p>-\frac{1}{n}$ and  $f,g,h:M\to [0,\infty)$ be three nonzero, compactly supported integrable functions  satisfying $(\ref{ConditionRescaledBBLWithWeights-vege}).$ Then 
			the following two assertions are equivalent: 
		\begin{itemize}
			\item[(a)] $\delta_{M,s}^p(f,g,h)=0,$ i.e.,  equality holds in the Borell-Brascamp-Lieb inequality$;$
%			\item[(b)]  $\rho_s=\mathcal \mathcal M_s^\frac{p}{pn+1}(\|f\|_1,\|g\|_1)h;$
			\item[(b)] the following statements simultaneously hold$:$ 
			\begin{itemize}
				\item[(i)] $\operatorname{supp}h=\psi_s(\operatorname{supp}f)$ up to a null measure set$;$
				\item[(ii)] ${\rm Jac}(\psi_s)(x)={v_{1-s}(\psi(x),x)}\left[\mathcal  M_s^\frac{p}{pn+1}\left(1,\frac{\|g\|_1}{\|f\|_1}\right)\right]^\frac{pn}{pn+1}$ for a.e. $x\in \operatorname{supp}f;$
				\item[(iii)]
				for
				a.e. $x\in \operatorname{supp} f$, one has $$\frac{h(\psi_s(x))}{\left[\mathcal  M_s^\frac{p}{pn+1}(\|f\|_1,\|g\|_1)\right]^\frac{1}{pn+1}}=\frac{ f(x)}{
					v_{1-s}(\psi(x),x)\|f\|_1^\frac{1}{pn+1}}=\frac{ g(\psi(x))}{ v_s(x,\psi(x))\|g\|_1^\frac{1}{pn+1}}.$$
			\end{itemize}
		\end{itemize}
	\end{theorem}
\medskip 

\noindent  The equality in the Borell-Brascamp-Lieb inequality for $p=-\frac{1}{n}$ is genuinely different than the case  $p>-\frac{1}{n}$ and it will be treated separately in Theorem \ref{Theorem-Riemannian-2}. 

\medskip 

%\begin{remark}\rm \label{remark-elso}
%	 We point out that  usually the inclusion $\psi_s(\operatorname{supp}f)\subset \operatorname{supp}h$ is strict. According to Theorem \ref{Theorem-Riemannian-egyenloseg}(b), the equality in the Borell-Brascamp-Lieb inequality implies the equality $\psi_s(\operatorname{supp}f)= \operatorname{supp}h$, which corresponds in $\mathbb R^n$ to the   {Alesker-Dar-Milman}  parametrization of the Minkowski sum of two sets; for further details, see Remark \ref{remark-adm}.\\

 	Theorem \ref{Theorem-Riemannian-egyenloseg} provides both well known and genuinely new rigidity results;  we briefly  present some of them in the sequel (for details, see \S \ref{section-Euclidean-0} and \S \ref{section-Riemannian-0}):\\ 

%\textcolor{red}{
$\bullet$ {\it Equality in the Borell-Brascamp-Lieb inequality in $\mathbb R^n$: a new approach to Dubuc's characterization.} As a first consequence of Theorem \ref{Theorem-Riemannian-egyenloseg}  we prove that equality in the Borell-Brascamp-Lieb inequality in $\mathbb R^n$ holds if and only if the functions $f,$ $g$ and $h$ are obtained as compositions of fixed $(t,p)$-concave function $\Phi$ with appropriate homotheties, where the support of $\Phi$ is convex up to a null set;  for the precise statement, see Theorem \ref{Theorem-Euklidean-egyenloseg}. This result provides a new qualitative formulation of Dubuc's characterization, see \cite[Th\'eor\`eme 12]{Dubuc}. Our strategy relies on applying  
 	Theorem \ref{Theorem-Riemannian-egyenloseg} in order to reduce the problem to the equality case in the Brunn-Minkowski inequality for the marginal supports $\operatorname{supp}f$ and $\operatorname{supp}g$, implying the convexity of these sets. Using the convexity of the support of the target measure,   a suitable application of the celebrated regularity result of Caffarelli \cite{Caffarelli} provides smoothness of the optimal mass transport map which turns to be an affine function in $\mathbb R^n$.  
% well known fact that 
%equality in the  Euclidean  Brunn-Minkowski inequality (\ref{BM-1}) holds for the convex sets $A$ and $B$ in $\mathbb R^n$ if and only if they are homothetic; see   Corollary \ref{Cor-Euclidean-2}. 
We notice that Caffarelli's regularity has been already employed in order to establish sharp stability results in $\mathbb R^n$ for the Brunn-Minkowski inequality (\ref{BM-1}), see  Figalli, Maggi and Pratelli \cite{FMP-Inventiones, FMP-AIHP}. 
 %are established for the Brunn-Minkowski inequality    involving convex sets 
%}

\medskip	
%\textcolor{red}{	
$\bullet$ {\it Equality in Borell-Brascamp-Lieb inequality implies constant curvature.} We  state that the equality in the Borell-Brascamp-Lieb inequality on an $n$-dimensional Riemannian manifold with Ricci curvature  Ric$(M)\geq k(n-1)$ for some $k\in \mathbb R$ can be expected to hold only when a particular region of the manifold between the marginal supports has {constant sectional curvature} $k$; see Theorem  \ref{Theorem-rigiditas} for details. The proof is based on Theorem \ref{Theorem-Riemannian-egyenloseg} and a careful comparison argument \`a la Bishop-Crittenden of the volume distortion coefficients with suitable quantities involving Jacobi fields on space forms. 
%}
	
%\end{remark}

%When $p=+\infty$ in Theorem \ref{Theorem-Riemannian-egyenloseg}, the expressions $\frac{p}{pn+1}$ and $\frac{1}{pn+1}$ are understood as $\frac{1}{n}$ and $0$, respectively. 

%Some remarks are in order concerning Theorems \ref{Theorem-Riemannian} \& \ref{Theorem-Riemannian-egyenloseg}, respectively. 

%Theorem \ref{Theorem-Riemannian} is proved by using the distorted Jacobian inequality from \cite{CMS} together with a quantitative H\"older inequality which contains the gap-function $G_s^{p,n}\geq 0$.  

 \medskip 
$\bullet$ {\it Equality in distorted Brunn-Minkowski inequality \`a la Lott-Sturm-Villani.} 
%As a particular case of the previous item, we establish rigidity result 
For some $s\in (0,1)$, $k\in \mathbb R$ and $n\geq 2$, let
$$\tau_s^{k,n}(\theta)=\left\{
\begin{array}{lll}
s^\frac{1}{n}\left(\sinh\left(\sqrt{-k}s\theta\right)\big/\sinh\left(\sqrt{-k}\theta\right)\right)^{1-\frac{1}{n}}&
{\rm if} & k\theta^2<0;\\
s & {\rm if} & k\theta^2=0;\\
s^\frac{1}{n}\left(\sin\left(\sqrt{k}s\theta\right)\big/\sin\left(\sqrt{k}\theta\right)\right)^{1-\frac{1}{n}}& {\rm if} &
0<k\theta^2<\pi^2;\\
+\infty & {\rm if} & k\theta^2\geq \pi^2,
\end{array}
\right.$$
be the distortion coefficient  introduced independently by Lott and Villani \cite{LV} and Sturm \cite{Sturm-2} in order to define their famous curvature-dimension condition  $\textsf{CD}(k,n)$ on metric measure spaces. Let $(M,w)$ be a complete $n$-dimensional Riemannian manifold with Ricci curvature bounded below, i.e., Ric$(M)\geq k(n-1)$ for some $k\in \mathbb R$ (which is equivalent to the validity of $\textsf{CD}(k,n)$) and 
% (which is equivalent to the validity of $\textsf{CD}(k,n)$, see e.g. Sturm \cite[Theorem 1.7]{Sturm-2}). 
let us denote by $\textsf{m}$ the Riemannian measure on $M$.
The distorted Brunn-Minkowski inequality reads as
\begin{equation}\label{BM-eredeti}
\textsf{m}(Z_s(A,B))^\frac{1}{n}\geq \tau_{1-s}^{k,n}(\Theta_{A,B})\textsf{m}(A)^\frac{1}{n}+\tau_{s}^{k,n}(\Theta_{A,B})\textsf{m}(B)^\frac{1}{n},
\end{equation}
 where $A,B\subset M$ are measurable sets with $\textsf{m}(A)\neq 0\neq \textsf{m}(B)$ and  
 \begin{eqnarray}\label{theta-nak}
 	\Theta_{A,B}= \left\{
 	\begin{array}{lll}
 		\inf_{(x,y)\in A\times B}d(x,y)
 		\  &\mbox{if} &  k\geq 0; \\
 		% u\geq 0 &\mbox{in} &   \Omega;\\
 		\sup_{(x,y)\in A\times B}d(x,y)
 		\  &\mbox{if} &  k< 0,
 	\end{array}\right.
 \end{eqnarray}
  see Sturm \cite[Proposition 2.1]{Sturm-2} and Villani \cite[Theorem 18.5]{Villani-1}.    
  Hereafter, the measure of $Z_s(A,B)$ in (\ref{BM-eredeti}) is always  understood w.r.t. the outer measure $\textsf{m}^*$ of $\textsf{m}.$ In fact, inequality (\ref{BM-eredeti}) is a direct consequence of the Borell-Brascamp-Lieb inequality (\ref{BBL-eredeti-Riemann}) even on metric measure spaces verifying the $\textsf{CD}(k,n)$ condition, see Bacher \cite[Proposition 3.2]{Bacher}; we also recall its proof in subsection \ref{subsection-4.2}.
% 
% 
%
% . The number $\tau_s^{k,n}$ 
% encodes information on the curvature of space forms serving as model structures in the definition of the Lott-Sturm-Villani's curvature-dimension condition  $\textsf{CD}(k,n)$ on metric measure spaces. 
 Theorem \ref{Theorem-Riemannian-egyenloseg} provides the following scenario concerning the equality in  (\ref{BM-eredeti}) (for details, see  Theorem \ref{Theorem-Riemannian-rigiditas-gorbulet-CD}): 
% in the sequel, $A\triangle B=(A\setminus B)\cup (B\setminus A)$ denotes the symmetric difference of $A$ and $B$.
%\textcolor{red}{
 %\noindent   we have the following conclusions on space forms. 
  \begin{enumerate}
  	\item[$\bullet$] {\it Positively curved case}: if $k>0$ (e.g., the round sphere $\mathbb S^n$), then equality in (\ref {BM-eredeti}) is characterized by the overlapping of the sets $Z_s(A,B)$, $A$ and $B$  up to a null measure set. Moreover, if $A\times B$ does not contain cut locus pairs,  equality holds in (\ref {BM-eredeti}) if and only if
  	there exists an open, geodesic convex set in $M$ which differs from $A$ and $B$ by a null set.
  	\item[$\bullet$] {\it Negatively curved case}: if $k<0$ and $(M,g)$ has nonpositive, nonzero sectional curvature (e.g., the hyperbolic space $\mathbb H^n$), 
  	equality in (\ref {BM-eredeti})  cannot hold for any positive measure sets $A$ and $B$.
%  	\item[(c)]  is characterized by the homothetic position  of the convex sets $A$ and $B$ in $\mathbb R^n$. 
  \end{enumerate}
The proof of the latter statements are based on a porosity argument and a geometric form of the Steinhaus density theorem (concerning the 'difference' of two sets).

%}      
  % In Section \ref{section-4} further rigidity results are presented. 
  
%
  %\textbf{
%  	Corollary \ref{log-BM-stability} is a simple consequence of Theorem \ref{Theorem-Riemannian}, providing a quantitative version of the log-Brunn-Minkowski inequality in $\mathbb R^n$; 
%  KELL???
%  	\noindent Hereafter, ${\rm Conv}(A)$ denotes the convex envelope\footnote{Conv$(A)$ is not the convex hull of $A$ in general; indeed, the convex hull $C(A)$ of $A$ is the smallest geodesic convex set containing $A$. In particular, if $A$ is a pair of antipodal points on a sphere, $C(A)$ cannot be uniquely determined while  ${\rm Conv}(A)$ is the whole sphere.} of $A\subset M$, i.e., ${\rm Conv}(A)=\cup_{i=0}^\infty A_i$, where $A_0=A$ and $A_i$ is the union of all geodesic segments between points of $A_{i-1}$, $i\geq 1$. }
%  

\medskip

\subsection{Organization of the paper}
In \S\ref{section-2} we first state a quantitative H\"older inequality (see Lemma \ref{lemma-p-mean}) which is crucial in the proof of  Theorem \ref{Theorem-Riemannian}. We then prove simultaneously Theorems \ref{Theorem-Riemannian}   and \ref{Theorem-Riemannian-egyenloseg}.  
In \S\ref{section-Euclidean-0} we prove Theorem \ref{Theorem-Euklidean-egyenloseg} which provides a qualitative version of Dubuc's characterization in $\mathbb R^n$ concerning the equality case in the Borell-Brascamp-Lieb inequality. In   \S\ref{section-Riemannian-0} we deal  with Riemannian manifolds by proving that the equality in Borell-Brescamp-Lieb inequality implies constant curvature, see  Theorem  \ref{Theorem-rigiditas},    and we discuss the equality cases in the distorted Brunn-Minkowski inequality (\ref{BM-eredeti}), see Theorem \ref{Theorem-Riemannian-rigiditas-gorbulet-CD}. 
%
%
% in \S\ref{subsection-Euclidean}, as a simple consequence of our main results we prove 
%Corollary \ref{log-BM-stability} and its slightly more general form. 
 In \S\ref{section-5-0} certain Borell-Brascamp-Lieb inequalities are presented on not necessarily reversible Finsler manifolds, highlighting some subtle differences between Riemannian/Euclidean and Finslerian frameworks, respectively, (see e.g. Corollary \ref{corollary-Finsler-24}). Finally, in \S\ref{appendix}, we provide the proof of  Lemma \ref{lemma-p-mean}.

	%In order to show the flavor of our results, we shall present in the sequel an abstract quantitative  Borell-Brascamp-Lieb inequality on a metric measure space having the differential structure required in the specific settings (i.e., Riemannian manifolds, Finsler manifolds and Heisenberg groups); details are presented later.  To be more precise, let $(M,d,\textsf{m})$ be a geodesic metric measure space with topological dimension $N\in \mathbb N$,  and $Z_s(x,y)$ be the $s$-intermediate points between $x$ and $y$ w.r.t. the  metric $d$ on $M$ defined by
	% $$ Z_s(x,y)=  \{ z \in M : d(x,z) = s d(x,y),\
	% d(z,y) = (1-s) d(x,y)\}.$$	Let $\mu_0$ and $\mu_1$ be two compactly supported probability measures on $M$ that are  absolutely continuous  w.r.t. the reference measure $\textsf{m}$. %  with probability %densities $\rho_0$ and  $\rho_1,$ respectively.
	% 	Assume that there exists a unique optimal transport map 
	% 	$\psi:M\to M$ transporting $\mu_0$ to $\mu_1$ associated to the cost function $\frac{d^2}{2}$.  If  $\psi_s$ denotes the interpolant optimal transport map associated to $\psi$, defined as
	% 	$$\psi_s(x) = Z_s(x,\psi(x)) \mbox{ for } \mu_0\mbox{-a.e. } x \in M,$$
	% we assume	the push-forward measure $\mu_s=(\psi_s)_{\#}\mu_0$ is also absolutely continuous w.r.t. $\textsf{m}$. 

	\section{A quantitative H\"older inequality and the Borell-Brascamp-Lieb deficit: proof of main results}\label{section-2}

	It is known (see e.g. Gardner \cite[Lemma
	10.1]{Gardner})  the following version of the H\"older inequality
	\begin{eqnarray}\label{MspIneq}
	\mathcal   M_s^{p}(a,b)\mathcal  M_s^{q}(c,d) \geq \mathcal M_s^{\eta}(ac, bd),
	\end{eqnarray} 
	holds for every $a,b,c,d \geq 0, s \in (0,1)$ and
	$p, q \in \R$ such that $p+q\geq0$ with ${\eta}=\frac{pq}{p+q}$ when
	$p$ and $q$ are not both zero, and $\eta=0$ if $p=q=0$.

 We first provide a technical  improvement of (\ref{MspIneq}) needed to  prove Theorem \ref{Theorem-Riemannian} whose proof is presented in the Appendix.

 \begin{lemma}{\bf (Quantitative H\"older inequality)}\label{lemma-p-mean}
		Let $n\in \mathbb N\setminus \{0\}$, $s\in (0,1)$ and $a,b,c,d>
		0$ be arbitrarily fixed numbers. For  $p >-\frac{1}{n}$ denote by  
		$\tilde p=\frac{p}{pn+1},$. If $p=+\infty$ then $\tilde p = \frac{1}{n}$  and if $p = -\frac{1}{n}$ we set  $\tilde p=-\infty$ .
		\begin{itemize}
			\item[{(i)}]  If $p\in (-\frac{1}{n},\infty)\setminus \{0\},$  then
			$$	\mathcal M_s^p(a,b)\mathcal M_s^{-\tilde p}(c,d)\geq \mathcal M_s^{-\frac{1}{n}}(ac,bd)\left[1+G_s^{p,n}(a,b,c,d) \right],$$
			where for $p >0$
			{\small \begin{eqnarray*}
					G_s^{p,n}(a,b,c,d)&=&(1-s)\frac{n}{\max(pn,1)}\left|\left[\mathcal M_s^{-p}\left(1,\frac{a}{b}\right)\right]^\frac{p\tilde p n}{\max(pn,1)}-\left[\mathcal M_s^{-\frac{1}{n}}\left(1,\frac{bd}{ac}\right)\right]^\frac{\tilde p}{\max(pn,1)}\right|^{\frac{\max(pn,1)}{\tilde pn}}\\
					&&		+s\frac{n}{\max(pn,1)}\left|\left[\mathcal M_s^{-p}\left(\frac{b}{a},1\right)\right]^\frac{p\tilde p n}{\max(pn,1)}-\left[\mathcal M_s^{-\frac{1}{n}}\left(\frac{ac}{bd},1\right)\right]^\frac{\tilde p}{\max(pn,1)}\right|^{\frac{\max(pn,1)}{\tilde pn}}.
			\end{eqnarray*}}
			and for $p<0$ 
			$$G_s^{p,n}(a,b,c,d)=G_s^{-\tilde p,n}(c,d,a,b).$$ Moreover,  $G_s^{p,n}(a,b,c,d)=0$ if and only if $\frac{a}{b}=\left(\frac{d}{c}\right)^\frac{1}{pn+1}.$ 
			\item[{(ii)}] If $p=0$  then
			$$	\mathcal M_s^0(a,b)\mathcal M_s^0(c,d)=\mathcal M_s^0(ac,bd)\geq \mathcal M_s^{-\frac{1}{n}}(ac,bd)\left[1+G_s^{0,n}(a,b,c,d) \right],$$
			where $$G_s^{0,n}(a,b,c,d)=n{\min(s,1-s)}
			\left[\mathcal M_s^\frac{1}{n}(bd,ac)\right]^{-\frac{1}{n}}\left|(bd)^\frac{\min(s,1-s)}{n}-(ac)^\frac{\min(s,1-s)}{n}\right|^\frac{1}{\min(s,1-s)}.$$
			Moreover,  $G_s^{0,n}(a,b,c,d)=0$ if and only if $ac=bd.$
			\item[{(iii)}] If $p=+\infty$ $($thus  $\tilde p=\frac{1}{n}),$ then 
			$$	\mathcal M_s^{+\infty}(a,b)\mathcal M_s^{-\frac{1}{n}}(c,d)\geq \mathcal M_s^{-\frac{1}{n}}(ac,bd)\left[1+G_s^{+\infty,n}(a,b,c,d) \right],$$
			where $$G_s^{+\infty,n}(a,b,c,d)=n\min(s,1-s)\frac{|a^\frac{1}{n}-b^\frac{1}{n}|}{(ab)^\frac{1}{n}\left[\mathcal M_s^{+\infty}(c,d)\right]^\frac{1}{n}}\left[\mathcal M_s^{-\frac{1}{n}}(ac,bd)\right]^\frac{1}{n}.$$
			Moreover,  $G_s^{+\infty,n}(a,b,c,d)=0$ if and only if $a=b.$
			\item[{(iv)}] If  $p= -\frac{1}{n}$, $($thus  $\tilde p=-\infty),$ then 
			$$	\mathcal M_s^{-\frac{1}{n}}(a,b)\mathcal M_s^	{+\infty}(c,d)\geq \mathcal M_s^{-\frac{1}{n}}(ac,bd)\left[1+G_s^{-\frac{1}{n}, n}(a,b,c,d) \right],$$
			where  $G_s^{-\frac{1}{n},n}(a,b,c,d)=G_s^{+\infty,n}(c,d,a,b).$
			Moreover, $G_s^{-\frac{1}{n},n}(a,b,c,d)=0$ if and only if $c=d.$
		\end{itemize}
	\end{lemma}

	\begin{remark}\rm \label{remark-1}
		 Let us observe the  homogeneity property of $G_s^{p,n}(\cdot,\cdot,\cdot,\cdot)$, i.e., for every $\lambda,\mu>0$ and $a,b,c,d>0$, one has 
		\begin{equation}\label{Gap-homogen}
		G_s^{p,n}(\lambda a,\lambda b,\mu c,\mu d)=G_s^{p,n}(a,b,c,d).
		\end{equation}
		%	Checking the cases in Lemma \ref{lemma-p-mean}, one has that 
		%	 $\mathcal G_s^{p,N}(a,b)=0$ if and only if $a=b.$
	\end{remark}
	
%	\medskip
	
%	\section{The Borell-Brascamp-Lieb deficit: proof of main results
	%	Theorems \ref{Theorem-Riemannian}  \& \ref{Theorem-Riemannian-egyenloseg} 
%}
%	\label{section-3}
	
%	Let $(M,d,m)$ be a metric measure space with topological dimension $N\in \mathbb N$, $p\geq -\frac{1}{N}$, $s\in (0,1),$ and three nonzero integrable functions $f,g,h:M\to [0,\infty)$. Let us consider the {\it Borell-Brascamp-Lieb deficit}, given by
%	$$\delta_{M,s}^p(f,g,h)=\frac{\|h\|_1}{\mathcal M_s^\frac{p}{pN+1}(\|f\|_1,\|g\|_1)}-1. $$
%	In particular, (\ref{elso-BBL-feltetel}) implies that $\delta_{\mathbb R^n,s}^p(f,g,h)\geq 0.$
	
	\medskip
Before  the proof of  Theorems \ref{Theorem-Riemannian} and \ref{Theorem-Riemannian-egyenloseg}, we recall some indispensable notions/results from the theory of optimal mass transportation on Riemannian manifolds. To do this, let $(M,w)$ be a complete $n$-dimensional Riemannian manifold and $d:M\times M\to \mathbb R$ be its distance function. 
Fixing $f$ and $g$ as in Theorem \ref{Theorem-Riemannian}, there exists a unique optimal transport map $\psi:M\to M$  from the measure $\mu=\tilde f {\rm d}\sf{m}$ to $\nu=\tilde g{\rm d}
\sf{m}$ with densities   $\tilde f=f/\|f\|_1$, $\tilde g=g/\|g\|_1$, see McCann \cite{McCann}. The optimal Brenier-type map has the form $\psi(x)=\exp_x(-\nabla \varphi(x))$, where $\varphi:M\to \mathbb R$  is a $d^2/2$-concave function, i.e., there exists a function $\eta:Y\to \mathbb R\cup \{-\infty\}$ with $\emptyset \neq Y\subset M$ such that
$$\varphi(x)=\inf_{y\in Y}\left(\frac{d^2(x,y)}{2}-\eta(y)\right),\ x\in M.$$
The $s$-interpolant optimal transport map 	$\psi_s:M\to M$ is  given by  $\psi_s(x)=\exp_x(-s\nabla \varphi(x))$ for a.e. $x\in \operatorname{supp} f$. For further use, let $A=\operatorname{supp}f$. We also have the injectivity of the interpolant  $\psi_s$ on $A$, see \cite[Lemma 5.3]{CMS}, i.e., if $\psi_s(x)=\psi_s(x')$ at two points $x,x'\in A$ of differentiability for the $d^2/2$-concave function $\varphi$, then $x=x'.$

	Let $B(x,r)=\{y\in
	M:d(x,y)<r\}$ be the geodesic ball with center $x\in M$ and radius  $r>0$. Fix $s\in (0,1).$ 
	According to Cordero-Erausquin, McCann and Schmuckenschl\"ager
	\cite{CMS},  the volume distortion coefficient in
	$(M, w)$ is defined by
	\begin{equation}\label{vol-distortion}
	v_s(x,y) = \lim\limits_{r \to 0}\frac{\textsf{m} \left( Z_s(x, B(y, r))\right)}{\textsf{m} \left( B(y, sr)\right)},
	\end{equation}
	where $\textsf{m}$ is the Riemannian measure.  	The Jacobian determinant inequality on $(M,w)$, cf. Cordero-Erausquin, McCann and Schmuckenschl\"ager
	\cite[Lemma 6.1]{CMS}, reads as
	\begin{equation}\label{Jacobian-inequality}
	{\rm Jac}(\psi_s)(x)\geq \mathcal M_s^\frac{1}{n}(v_{1-s}(\psi(x),x),v_{s}(x,\psi(x)){\rm Jac}(\psi)(x))\ \ {\rm for\ a.e.}\ x\in A.
	\end{equation}
	We notice that the Monge-Amp\`ere equation holds, i.e.  
	\begin{equation}\label{Monge-Ampere}
	\tilde f(x)=\tilde g(\psi(x)){\rm Jac}(\psi)(x)\ {\rm for \ a.e.}\ x\in A.
	\end{equation}
\\

	{\it Proof of Theorems \ref{Theorem-Riemannian}   and \ref{Theorem-Riemannian-egyenloseg}.} The proof  of these results will be presented simultaneously.
	   	Let $$\tilde h(z)=\frac{h(z)}{\mathcal M_s^\frac{p}{pn+1}(\|f\|_1,\|g\|_1)},\ z\in M.$$ 
	   We first notice that $$\psi_s(A)\subseteq \operatorname{supp}h,$$ up to a null measure set. Indeed, if $x\in A$, then $\psi(x)\in \operatorname{supp}g$ and by the hypothesis (\ref{ConditionRescaledBBLWithWeights-vege}) and  convention on $\mathcal M_s^p$, it follows that $h(\psi_s(x))>0.$
%	 and its Jacobian is 
%	\begin{equation}\label{jacobi-s}
%	{\rm Jac}(\psi_s)(x)={\rm det}\left[Y(s){\rm Hess}\left.\left[\frac{d^2(\psi_s(x),\cdot)}{2}-s\varphi(\cdot) \right]\right|_{x}\right],
%	\end{equation}
%where	$Y(s)=d(\exp_x)_{-s\nabla \varphi(x)}$ is the Jacobian of the exponential map at $-s\nabla \varphi(x)\in T_xM.$ 
%	Similarly, we have that  $$ {\rm Jac}(\psi)(x)={\rm det}\left[Y(1){\rm Hess}\left.\left[\frac{d^2(\psi(x),\cdot)}{2}-\varphi(\cdot) \right]\right|_{x}\right]. $$	
%	
%	\textcolor{red}{Ugy tÃŒnik hogy az elÃ¶bi ket kepletet a Jacobiannak egyaltalan nem hasznaljuk a bizonyitasban. Talan el is hagyhatjuk akkor.... Mit szolsz hozza?}
	We consider several cases according to the values of $p$.
	
	\textbf{\underline{{Case 1}}:}   $p\in (-\frac{1}{n},\infty)\setminus \{0\}.$ 
{Integrating with respect to $\textsf{m}$, by the change} of variable $z=\psi_s(x)$ (since $\psi_s$ is injective), we obtain 
	\begin{eqnarray*}
		\|\tilde h\|_1&=&\int_M\tilde h =\int_{\operatorname{supp}h}\tilde h\\&\geq& \int_{\psi_s(A)}\tilde h=\int_{A}\tilde h(\psi_s(x)){\rm Jac}(\psi_s)(x)\\
		&\geq & \int_{A}\mathcal M^{p}_s
		\left(\frac{f(x)}{v_{1-s}(\psi(x),x)},\frac{g(\psi(x))}{v_s(x,\psi(x))} \right)\mathcal M_s^{-\frac{p}{pn+1}}\left(\frac{1}{\|f\|_1},\frac{1}{\|g\|_1}\right){\rm Jac}(\psi_s)(x)\ \ \ \ \ \ \  {\rm (see \ (\ref{ConditionRescaledBBLWithWeights-vege}))}  \\&\geq& 
		\int_A \mathcal  M^{-\frac{1}{n}}_s
		\left(\frac{\tilde f(x)}{v_{1-s}(\psi(x),x)},\frac{\tilde g(\psi(x))}{v_s(x,\psi(x))} \right) {\rm Jac}(\psi_s)(x)\times \\&&  \ \ \ \times \left(1+G_s^{p,n}\left(\frac{f(x)}{v_{1-s}(\psi(x),x)},\frac{g(\psi(x))}{v_s(x,\psi(x))},\frac{1}{\|f\|_1},\frac{1}{\|g\|_1}\right) \right) \ \ \ \ \ \ \ \ \ \ \ \ \ {\rm (see \ Lemma\  \ref{lemma-p-mean}\ (i))}
		\\&\geq& 
		\int_A \mathcal  M^{-\frac{1}{n}}_s
		\left(\frac{\tilde f(x)}{v_{1-s}(\psi(x),x)},\frac{\tilde g(\psi(x))}{v_s(x,\psi(x))} \right) \mathcal M_s^\frac{1}{n}\left(v_{1-s}(\psi(x),x),v_{s}(x,\psi(x))\frac{\tilde f(x)}{\tilde g(\psi(x))}\right)\times \\&&  \ \ \ \times \left(1+G_s^{p,n}\left(\frac{f(x)}{v_{1-s}(\psi(x),x)},\frac{g(\psi(x))}{v_s(x,\psi(x))},\frac{1}{\|f\|_1},\frac{1}{\|g\|_1}\right) \right)\ \ \ \ \ \ \ \ \ \ \ \ \ \ \ \ \ \ \ \ \ \ \  \ \ {\rm (see \ (\ref{Jacobian-inequality}))}
		\\& = &1+\int_M \tilde f(x)G_s^{p,n}\left(\frac{f(x)}{v_{1-s}(\psi(x),x)},\frac{g(\psi(x))}{v_s(x,\psi(x))},\frac{1}{\|f\|_1},\frac{1}{\|g\|_1}\right),\ \ \ \ \ \ \ \ \ \ \ \ \ \ \ \ \ \ \ \ \ \  \ \ {\rm (see \ (\ref{Monge-Ampere}))}
	\end{eqnarray*}
	where the last equality   follows by the relation 
	$$\mathcal M_s^{-\frac{1}{n}}(\lambda a,\lambda b)\mathcal M_s^{\frac{1}{n}}(1/ a,1/ b)=\lambda,\ \ \ \ \ \ a,b,\lambda>0,$$
	which proves {Theorem \ref{Theorem-Riemannian}}.

	 Now, assume that (a) holds, i.e., $\delta_{M,s}^p(f,g,h)=	\|\tilde h\|_1-1=0$.  It follows directly that 
	$$G_s^{p,n}\left(\frac{f(x)}{v_{1-s}(\psi(x),x)},\frac{g(\psi(x))}{v_s(x,\psi(x))},\frac{1}{\|f\|_1},\frac{1}{\|g\|_1}\right)=0\ \ \ {\rm for \ a.e.}\ x\in A,$$
	and there are equalities in the above estimates. In particular, $$\operatorname{supp}\tilde h=\operatorname{supp}h=\psi_s(A),$$ up to a null measure set of $M$, which gives property (i) of Theorem
	\ref{Theorem-Riemannian-egyenloseg}. By the characterization of $G_s^{p,n}(a,b,c,d)=0$ (see Lemma \ref{lemma-p-mean} (i)), the latter relation is equivalent to $$\frac{f(x)}{
		v_{1-s}(\psi(x),x)\|f\|_1^\frac{1}{pn+1}}=\frac{g(\psi(x))}{ v_s(x,\psi(x))\|g\|_1^\frac{1}{pn+1}}\ \ {\rm for\ a.e.}\ x\in A.$$  By (\ref{ConditionRescaledBBLWithWeights-vege}) and the above estimate we necessarily have for a.e. $x\in A$ that 
	\begin{eqnarray*}
	h(\psi_s(x))&=&\mathcal M^{p}_s
	\left(\frac{f(x)}{v_{1-s}(\psi(x),x)},\frac{g(\psi(x))}{v_s(x,\psi(x))} \right)=\frac{f(x)}{
		v_{1-s}(\psi(x),x)\|f\|_1^\frac{1}{pn+1}}\left[\mathcal M_s^\frac{p}{pn+1}(\|f\|_1,\|g\|_1)\right]^\frac{1}{pn+1},
	\end{eqnarray*}
	which is (iii) of Theorem
	\ref{Theorem-Riemannian-egyenloseg}. Since we also have equality in the Jacobi determinant inequality  (\ref{Jacobian-inequality}), 
	property (ii) {of Theorem
	\ref{Theorem-Riemannian-egyenloseg} } directly follows by (iii); thus every item of (b) holds true.  
	The reverse implication is trivial. 
	
	\textbf{\underline{{Case 2}}:}   $p=+\infty.$ A similar reasoning as in Case 1 and Lemma \ref{lemma-p-mean} (iii) give that 
	$$\|\tilde h\|_1\geq 1+\int_M \tilde f(x)G_s^{+\infty,n}\left(\frac{f(x)}{v_{1-s}(\psi(x),x)},\frac{g(\psi(x))}{v_s(x,\psi(x))},\frac{1}{\|f\|_1},\frac{1}{\|g\|_1}\right).$$
	
	If $\delta_{M,s}^{+\infty}(f,g,h)=0$, the latter integrand is necessarily zero.  Since $G_s^{+\infty,n}(a,b,c,d)=0$ if and only if $a=b$  (see Lemma \ref{lemma-p-mean} (iii)), we obtain $$\frac{f(x)}{
		v_{1-s}(\psi(x),x)}=\frac{g(\psi(x))}{ v_s(x,\psi(x))}\ \ {\rm for\ a.e.}\ x\in A.$$ Furthermore, in order to have the equality case, by  (\ref{ConditionRescaledBBLWithWeights-vege}) and the latter relation we necessarily have for a.e. $x\in A$ that 
	$$h(\psi_s(x))=\mathcal M^{+\infty}_s
	\left(\frac{f(x)}{v_{1-s}(\psi(x),x)},\frac{g(\psi(x))}{v_s(x,\psi(x))} \right)=\frac{f(x)}{
		v_{1-s}(\psi(x),x)}=\frac{g(\psi(x))}{ v_s(x,\psi(x))},$$ 
	which corresponds  to (iii) of Theorem
	\ref{Theorem-Riemannian-egyenloseg}. Clearly, one also has (i) and by the equality in (\ref{Jacobian-inequality}) we necessarily have for  a.e. $ x\in A$ that
	$${\rm Jac}(\psi_s)(x)= \mathcal M_s^\frac{1}{n}\left(v_{1-s}(\psi(x),x),v_{s}(x,\psi(x))\frac{\tilde f(x)}{\tilde g(\psi(x))}\right)=\frac{v_{1-s}(\psi(x),x)}{\|f\|_1}\mathcal M_s^\frac{1}{n}(\|f\|_1,\|g\|_1),$$ 
	which is precisely (ii) of Theorem
	\ref{Theorem-Riemannian-egyenloseg}.  The converse is trivial again. 
	
	\textbf{\underline{{Case 3}}:}   $p=0.$ Similarly as above, by Lemma \ref{lemma-p-mean} {(ii)} we have
	\begin{eqnarray*}
		\|\tilde h\|_1&\geq& \int_{\psi_s(A)}\tilde h=\int_{A}\tilde h(\psi_s(x)){\rm Jac}(\psi_s)(x)\\
		&\geq & \int_{A}\mathcal M^{0}_s
		\left(\frac{f(x)}{v_{1-s}(\psi(x),x)},\frac{g(\psi(x))}{v_s(x,\psi(x))} \right)\mathcal M_s^{0}\left(\frac{1}{\|f\|_1},\frac{1}{\|g\|_1}\right){\rm Jac}(\psi_s)(x)\\
		&\geq &1+\int_M \tilde f(x)  G_s^{0,n}\left(\frac{ f(x)}{v_{1-s}(\psi(x),x)},\frac{ g(\psi(x))}{v_s(x,\psi(x))},\frac{1}{\|f\|_1},\frac{1}{\|g\|_1}\right).
	\end{eqnarray*}
	Let us assume that $\delta_{M,s}^{0}(f,g,h)=0$; thus, the latter integrand is zero.  Note that $G_s^{0,n}(a,b,c,d)=0$ if and only if $ac=bd$  (see Lemma \ref{lemma-p-mean} (ii));  therefore, we obtain $$\frac{\tilde f(x)}{
		v_{1-s}(\psi(x),x)}=\frac{\tilde g(\psi(x))}{ v_s(x,\psi(x))}\ \ {\rm for\ a.e.}\ x\in A.$$
	Having equality in (\ref{ConditionRescaledBBLWithWeights-vege}), from the latter relation we obtain 
	for a.e. $x\in A$ that 
	$$h(\psi_s(x))=\mathcal M^{0}_s
	\left(\frac{f(x)}{v_{1-s}(\psi(x),x)},\frac{g(\psi(x))}{v_s(x,\psi(x))} \right)=\frac{\tilde f(x)}{
		v_{1-s}(\psi(x),x)}\mathcal M_s^0(\|f\|_1,\|g\|_1),$$ 
	which is  (iii) of Theorem
	\ref{Theorem-Riemannian-egyenloseg}. Property (i) follows trivially, while (ii) comes from (iii) and the equality in (\ref{Jacobian-inequality}), i.e., 
	$${\rm Jac}(\psi_s)(x)= \mathcal M_s^\frac{1}{n}\left(v_{1-s}(\psi(x),x),v_{s}(x,\psi(x))\frac{\tilde f(x)}{\tilde g(\psi(x))}\right)={v_{1-s}(\psi(x),x)}\ \ {\rm for\ a.e.}\ x\in A.$$ 
	
	\textbf{\underline{{Case 4}}:}   $p=-\frac{1}{n}.$ The proof is similar to the case $p=+\infty$; indeed, 
	one has 
	$$\|\tilde h\|_1\geq 1+\int_M \tilde f(x)G_s^{-\frac{1}{n},n}\left(\frac{f(x)}{v_{1-s}(\psi(x),x)},\frac{g(\psi(x))}{v_s(x,\psi(x))},\frac{1}{\|f\|_1},\frac{1}{\|g\|_1}\right).$$
	By Lemma \ref{lemma-p-mean} (iv), the claim follows. 	 The equality case is treated in the following result. 
	\hfill $\square$
	
	\begin{theorem}\label{Theorem-Riemannian-2} {\bf (Equality in Borell-Brascamp-Lieb inequality; $p=-\frac{1}{n}$)}  Let us assume that the assumptions  in Theorem \ref{Theorem-Riemannian} are fulfilled.  
		%	if
		%	\begin{eqnarray}\label{ConditionRescaledBBLWithWeights-vege-2} h(z)
		%	\geq M^{-\frac{1}{n}}_s
		%	\left(\frac{f(x)}{v_{1-s}(y,x)},\frac{g(y)}{v_s(x,y)} \right) \
		%	\  { for\ all}\ (x,y)\in M\times M, z\in Z_s(x,y),
		%	\end{eqnarray}
		%	then 
		%	$$\displaystyle\delta_{M,s}^{-\frac{1}{n}}(f,g,h)\geq \displaystyle\int_M \tilde f(x)G_s^{-\frac{1}{n},n}\left(\frac{f(x)}{v_{1-s}(\psi(x),x)},\frac{g(\psi(x))}{v_s(x,\psi(x))},\frac{1}{\|f\|_1},\frac{1}{\|g\|_1}\right)$$
		%	 %and  $\mathcal G_s^{p,n}(\cdot,\cdot)$ comes from {\rm (\ref{R-definicioja})}. 
		Then the following assertions are equivalent$:$ 
		\begin{itemize}
			\item[(a)] $\delta_{M,s}^{-\frac{1}{n}}(f,g,h)=0;$ 
			\item[(b)]	the following statements simultaneously hold$:$ 
			\begin{itemize}
				\item[(i)] $\operatorname{supp}h=\psi_s(\operatorname{supp}f)$ up to a null measure set$;$
				\item[(ii)] $h(\psi_s(x))= \mathcal M_s^{-\frac{1}{n}}\left(\frac{f(x)}{v_{1-s}(\psi(x),x)},\frac{g(\psi(x))}{v_{s}(x,\psi(x))}\right)=\frac{f(x)}{ {\rm Jac}(\psi_s)(x)}$
				for a.e. $x\in \operatorname{supp}f;$
				\item[(iii)] $\|f\|_1=\|g\|_1.$
			\end{itemize}
		\end{itemize}
	\end{theorem}
	
	{\it Proof.} Assume first that $\delta_{M,s}^{-\frac{1}{n}}(f,g,h)=0,$ i.e.,  $\|\tilde h\|_1=1.$ It follows by Case 4  of the previous proof and the characterization of the equality $G_s^{-\frac{1}{n},n}(a,b,c,d)=0$ (see  Lemma \ref{lemma-p-mean} (iv)) that $\|f\|_1=\|g\|_1,$ which is precisely property (iii). As in the previous cases, we also have that $\operatorname{supp}h=\psi_s(\operatorname{supp}f)$ up to a null measure set$.$
	In order to prove (ii), we have by (\ref{ConditionRescaledBBLWithWeights-vege}) that
%	we observe that  
%	$$1=\|\tilde h\|_1=\frac{\| h\|_1}{\mathcal  M_s^{-\infty}\left(\|f\|_1,\| g\|_1\right)}=\frac{\| h\|_1}{\|f\|_1},$$ i.e., $\|h\|_1=\|f\|_1=\|g\|_1.$ We also have that 
	$h(\psi_s(x))= \mathcal M_s^{-\frac{1}{n}}\left(\frac{f(x)}{v_{1-s}(\psi(x),x)},\frac{g(\psi(x))}{v_{s}(x,\psi(x))}\right)$ for a.e. $ x\in A$ and 
	\begin{eqnarray*}
	{\rm Jac}(\psi_s)(x)&=& \mathcal M_s^\frac{1}{n}\left(v_{1-s}(\psi(x),x),v_{s}(x,\psi(x))\frac{\tilde f(x)}{\tilde g(\psi(x))}\right)\\&=&\mathcal M_s^\frac{1}{n}\left(v_{1-s}(\psi(x),x),v_{s}(x,\psi(x))\frac{ f(x)}{ g(\psi(x))}\right)\ \ \ \ \ \ \ \ \ \ \ \ \ \ \ \ \ \    ({\rm since}\ \|f\|_1=\|g\|_1) \\
	&=&f(x)\mathcal M_s^\frac{1}{n}\left(\frac{v_{1-s}(\psi(x),x)}{f(x)},\frac{ v_{s}(x,\psi(x))}{ g(\psi(x))}\right)\\&=&f(x)\left[\mathcal M_s^{-\frac{1}{n}}\left(\frac{f(x)}{v_{1-s}(\psi(x),x)},\frac{ g(\psi(x))}{ v_{s}(x,\psi(x))}\right)\right]^{-1}\\&=&\frac{f(x)}{h(\psi_s(x))}.
	\end{eqnarray*}
	
	 Conversely, we assume that (i)-(iii) hold. Then
	\begin{eqnarray*}
	\|\tilde h\|_1&=&\int_{\operatorname{supp}h}\tilde h=\int_{\psi_s(\operatorname{supp}f)}\tilde h(z){\rm d}{\sf m}(z)\ \ \ \ \ \ \ \ \ \ \ \ \ \ \ \ \ \ \ \ \ \ \ \ \ \ \ \ \ \ \ \ \ \  \ \ \ \ \ \ \ \ \ \ \ \ \ \ \ \ \ \ \ \ \ \ \ ({\rm see}\ ({\rm i}))\\&=& \frac{\displaystyle\int_{\operatorname{supp}f} h(\psi_s(x)){\rm Jac}(\psi_s)(x){\rm d}{\sf m}(x)}{\mathcal  M_s^{-\infty}\left(\|f\|_1,\| g\|_1\right)}\ \ \ \ \ \ \ \ \ \ \ \    \ \ \ \ \ \ \ \ \ \ \ \ \ \  ({\rm change\ of\ variables}\ z=\psi_s(x))
	\\&=& \frac{\displaystyle\int_{\operatorname{supp}f} f(x){\rm d}{\sf m}(x)}{\|f\|_1}\ \ \ \ \ \ \ \ \ \ \ \ \ \ \ \ \ \ \ \   \ \ \ \ \ \ \ \ \ \ \ \ \ \ \ \ \ \ \ \ \ \ \ \ \ \ \ \  \ \ \ \ \ \ \ \ \ \ \ \ \ \ \ \ \ \ \ \ \   ({\rm see\ (ii)\&(iii)})\\&=&1,
	\end{eqnarray*}
which concludes the proof. \hfill $\square$
	
\begin{remark}\rm \label{remark-4-14} 
	Note the difference in equality cases between $p > -\frac{1}{n}$ and $p = -\frac{1}{n}$, respectively. However, in all cases, it holds (see Cases 1-4) that
$$	{h(\psi_s(x))}=\mathcal M^p_s
	\left(\frac{f(x)}{v_{1-s}(\psi(x),x)},\frac{g(\psi(x))}{v_s(x,\psi(x))} \right)\ \ {\rm for\ a.e.}\ x\in A.$$
\end{remark}

\medskip
	Our main results (Theorems \ref{Theorem-Riemannian} and \ref{Theorem-Riemannian-egyenloseg}) can be efficiently applied to establish various rigidity results. In \S\ref{section-Euclidean-0} we consider the Euclidean case, in \S\ref{section-Riemannian-0}  the case of Riemannian manifolds, while in \S\ref{section-5-0} we discuss the case of Finsler manifolds. The notations are kept from the previous sections.

	\section{Dubuc's result recovered via optimal mass transportation}\label{section-Euclidean-0}

	Let $t\in (0,1)$ and $p\in \mathbb R\cup \{+\infty\}$. We say that a nonnegative integrable function $\Phi:K\to \mathbb R$ is $(t,p)$-\textit{concave} on the convex set $K\subset \mathbb R^n$ if
	$$\Phi((1-t)x+ty)\geq \mathcal M_t^p(\Phi(x),\Phi(y))\ \ {\rm for\ all}\ x,y\in K.$$
	If $\Phi$ is continuous on $K$  then the $(t_0,p)$-concavity of $\Phi$ for some $t_0\in (0,1)$ implies the $(t,p)$-concavity of $\Phi$ for every $t\in (0,1).$ In such a case, the latter notation is simply called $p$-concavity, see Gardner \cite[Section 9]{Gardner}.  In particular, in the latter case, the $p$-concavity of $\Phi$ in $K$ means that 
	$\Phi^p$ is concave  in $K$ if $p>0$, $\Phi^p$ is convex  in $K$ if $p<0$, $\Phi$ is log-concave  in $K$ if $p=0$, and $\Phi$ is constant  in $K$ if $p=+\infty. $ 
	
	The main result of this section provides a novel, qualitative characterization of the equality case in the Borell-Brascamp-Lieb inequality, complementing the result of Dubuc \cite{Dubuc} (see also Rossi \cite{Rossi-PhD} and Rossi and Salani \cite{Rossi-Salani, Rossi-Salani-AA}):

	\begin{theorem}\label{Theorem-Euklidean-egyenloseg} {\bf (Equality in Borell-Brascamp-Lieb inequality; Euclidean case)} 
		Let $s\in (0,1),$ $p\geq-\frac{1}{n}$ and  $f,g,h:\mathbb R^n\to [0,\infty)$ be three nonzero, compactly supported integrable functions  satisfying $(\ref{elso-BBL-feltetel}).$ Then 
		the following two assertions are equivalent: 
		\begin{itemize}
			\item[(a)] $\delta_{\mathbb R^n,s}^p(f,g,h)=0,$ i.e., equality holds in the Borell-Brascamp-Lieb inequality {\rm (\ref{BBL-euklidesz})}$;$
			%			\item[(b)]  $\rho_s=\mathcal \mathcal M_s^\frac{p}{pn+1}(\|f\|_1,\|g\|_1)h;$
			\item[(b)] there exist an element $x_0\in \mathbb R^n$, a convex set $K\subset \mathbb R^n$ with $K=\operatorname{supp}f$ up to a null measure set and a $(t,p)$-concave function $\Phi:K\to \mathbb R$ with $t=\frac{sc_0}{1-s+sc_0}$ and  $c_0=\left(\frac{\mathcal L^n(\operatorname{supp}g)}{\mathcal L^n(\operatorname{supp}f)}\right)^\frac{1}{n}$ such that up to null measure sets 
			\begin{equation}\label{transzlacio}
			\operatorname{supp}g=c_0\operatorname{supp}f +x_0\ \  \ {\rm and}\ \  \operatorname{supp}h=(1-s+sc_0)\operatorname{supp}f +sx_0,
			\end{equation}
			and for a.e. $x\in K,$
			\begin{equation}\label{Dubuc-fuggvenyek}
			\left\{\begin{array}{lll}
			f(x)=\Phi(x); && \\
			% u\geq 0 &\mbox{in} &   \Omega;\\
			g(c_0x+x_0)=c_0^\frac{1}{p}\Phi(x); &&\\
			h((1-s+sc_0)x +sx_0)=\left[\mathcal  M_s^\frac{p}{pn+1}\left(1,c_0^\frac{pn+1}{p}\right)\right]^\frac{1}{pn+1}\Phi(x).
			\end{array}\right.
			\end{equation}

		\end{itemize}
	\end{theorem}
\noindent 	Hereafter, the following two conventions are used: 

$\bullet$ if $p=0$ then it will turn out by the proof that  $c_0=1$, thus we may consider $c_0^\frac{1}{p}=1;$ 

$\bullet$ if $p=-\frac{1}{n}$, we consider $$\lim_{p\to -\frac{1}{n}}\left[\mathcal  M_s^\frac{p}{pn+1}\left(1,c_0^\frac{pn+1}{p}\right)\right]^\frac{1}{pn+1}=\mathcal M_s^{-\frac{1}{n}}(1,c_0^{-n}).$$
%	which is again in a perfect concordance with (\ref{Dubuc-egyenloseg-2}). 
\vspace{0.9cm}

	{\it Proof of Theorem \ref{Theorem-Euklidean-egyenloseg}.} (a)$\implies$(b) 	We distinguish two cases. 
	
	\textbf{\underline{Case 1}:} $p>-\frac{1}{n}.$ {Taking into consideration that in the Euclidean case the distortion coefficients $v_s(x,y)$ are identically equal to $1$}, according to Theorem \ref{Theorem-Riemannian-egyenloseg}, the  equality in the Borell-Brascamp-Lieb inequality, i.e., $\delta_{\mathbb R^n,s}^p(f,g,h)=0,$ is characterized by:

	\begin{itemize}
		\item[(i)] $\operatorname{supp}h=\psi_s(\operatorname{supp}f)$ up to a null measure set$;$
		\item[(ii)] ${\rm Jac}(\psi_s)(x)=\left[\mathcal  M_s^\frac{p}{pn+1}\left(1,\frac{\|g\|_1}{\|f\|_1}\right)\right]^\frac{pn}{pn+1}$ for a.e. $x\in \operatorname{supp}f;$
		\item[(iii)]
		for
		a.e. $x\in \operatorname{supp} f$, one has $$\frac{h(\psi_s(x))}{\left[\mathcal  M_s^\frac{p}{pn+1}(\|f\|_1,\|g\|_1)\right]^\frac{1}{pn+1}}=\frac{ f(x)}{
			\|f\|_1^\frac{1}{pn+1}}=\frac{ g(\psi(x))}{ \|g\|_1^\frac{1}{pn+1}}.$$
	\end{itemize}
	For simplicity, let $A=\operatorname{supp}f$ and $B=\psi(A).$ We also recall that  $\psi:A\to B$ is the optimal transport map from the measure $\mu=\tilde f{\rm d}\mathcal L^n$ to $\nu=\tilde g{\rm d}\mathcal L^n$, where $\tilde f=f/\|f\|_1$ and $\tilde g=g/\|g\|_1.$
	% where $\mathcal L^n$ is the usual $n$-dimensional outer Lebesgue measure; 
	In fact, $$\psi(x)=\exp_x(-\nabla \varphi(x))=x-\nabla \varphi(x)$$ for some  $|\cdot|^2/2$-concave function $\varphi:\mathbb R^n\to \mathbb R$.  Equivalently,  there exists a convex function   $\eta:\mathbb R^n\to \mathbb R,$ $\eta(x)=\frac{|x|^2}{2}-\varphi(x)$ such that $\psi=\nabla \eta$ and  $\psi_\# \mu=\nu$, see Villani \cite[p.187]{Villani-2}. Accordingly, $$\psi_s(x)=x-s\nabla \varphi(x)=(1-s)x+s\nabla \eta(x).$$ 
	
	It is clear by (\ref{elso-BBL-feltetel}) (or (\ref{ConditionRescaledBBLWithWeights-vege}))
	and the definition of $\mathcal M_s^p$ that $$\psi_s(\operatorname{supp}f)= \psi_s(A)\subseteq Z_s(A,B)\subseteq \operatorname{supp}h.$$ Now, in particular, (i) implies that $\mathcal L^n(Z_s(A,B))=\mathcal L^n(\psi_s(A)).$ By a change of variables and (ii), it follows that
	\begin{eqnarray*}
		\mathcal L^n(Z_s(A,B))&=&\mathcal L^n(\psi_s(A))=\int_{\psi_s(A)}{\rm d}\mathcal L^n=\int_{A}{\rm Jac}(\psi_s)(x){\rm d}\mathcal L^n(x)\\&=&\left[\mathcal  M_s^\frac{p}{pn+1}\left(1,\frac{\|g\|_1}{\|f\|_1}\right)\right]^\frac{pn}{pn+1}\mathcal L^n(A)\\&=&\left(1-s+s\left(\frac{\|g\|_1}{\|f\|_1}\right)^\frac{p}{pn+1}\right)^n\mathcal L^n(A).
	\end{eqnarray*}
	On the other hand, by the Monge-Amp\`ere equation (\ref{Monge-Ampere}) for $\tilde f$ and $\tilde g$, one has $\tilde f(x)=\tilde g(\psi(x)){\rm Jac}(\psi)(x)$  for a.e. $ x\in A;$
	in particular, by the last relation of (iii) we have that 
	\begin{equation}\label{regularitashoz-Caffarelli}
	{\rm Jac}(\psi)(x)=\left(\frac{\|g\|_1}{\|f\|_1}\right)^\frac{pn}{pn+1}\ \ {\rm for\ a.e.}\ x\in A.
	\end{equation}
	Therefore, by (\ref{regularitashoz-Caffarelli}) one has
	\begin{eqnarray*}
		\mathcal L^n(B)&=&\mathcal L^n(\psi(A))=\int_{\psi(A)}{\rm d}\mathcal L^n=\int_{A}{\rm Jac}(\psi)(x){\rm d}\mathcal L^n(x)\\&=&\left(\frac{\|g\|_1}{\|f\|_1}\right)^\frac{pn}{pn+1}\mathcal L^n(A).
	\end{eqnarray*}
	Combining the above two relations, we obtain that 
	$$\mathcal L^n(Z_s(A,B))^\frac{1}{n}=(1-s)\mathcal L^n(A)^\frac{1}{n}+s\mathcal L^n(B)^\frac{1}{n},$$
	i.e., we have equality in the Brunn-Minkowski inequality.  It is known (see, e.g.,  Gardner \cite[p. 363]{Gardner}) that then $A$ and $B$ are homothetic convex bodies (i.e., compact convex sets with non-empty interior) from which sets of measure zero are removed. Let $K$ and $S$ be these convex bodies   which differ from $A$ and $B$ by null sets, respectively, and let $c_0>0$ and $x_0\in \mathbb R^n$ such that
	$S=c_0K+x_0.$  Without loss of generality, we may consider in the sequel the set of interior points $K^{\circ}$ and $S^{\circ}$ instead of the sets themselves, having the same measures as $K$ and $S$, respectively. It is clear that $c_0=\left(\frac{\mathcal L^n(B)}{\mathcal L^n(A)}\right)^{1/n}$. Since $S$ is convex, relation (\ref{regularitashoz-Caffarelli}) and the interior regularity result of Caffarelli \cite{Caffarelli}  imply that $\eta$ is of class $C^2$ in the interior of this set. Thus the Aleksandrov second derivative $D_o^2 \eta$ becomes the usual Hessian of $\eta$, see Villani \cite[Theorem 4.14]{Villani-2}.   Moreover, (ii) and (\ref{regularitashoz-Caffarelli}) imply that
	$${\rm det}^\frac{1}{n}[(1-s)I_n+s{\rm Hess} \eta(x)]={\rm Jac}(\psi_s)(x)^\frac{1}{n}=(1-s){\rm det}^\frac{1}{n}[I_n]+s{\rm det}^\frac{1}{n}[{\rm Hess} \eta(x)],\ x\in K^{\circ};$$
	we emphasize that the latter relation is valid for all $x\in K^{\circ}$ (since $\eta \in C^2$ on $K^{\circ}$) and not only for a.e. $x\in K^{\circ}$. 
	This relation and the strict concavity of det$^\frac{1}{n}(\cdot)$ over the cone of nonnegative definite symmetric matrices  give that {\rm Hess}$\eta(x)=c_0 I_n$ for every $x\in K^{\circ}$, where 
	\begin{equation}\label{c0-definicio}
c_0=\left(\frac{\mathcal L^n(B)}{\mathcal L^n(A)}\right)^{1/n}=\left(\frac{\|g\|_1}{\|f\|_1}\right)^\frac{p}{pn+1}.
	\end{equation}
	Therefore, 
	$$\psi(x)=\nabla \eta(x)=c_0x+x_0\ \ \ {\rm and}\ \ \ \psi_s(x)=(1-s+sc_0)x+sx_0,\ \ x\in K^{\circ}. $$
	%}
	By continuity, these relations hold true for all $x\in K$. 
	Accordingly, by (iii) we have that
	\begin{equation}\label{Dubuc-egyenloseg}
	\frac{h((1-s+sc_0)x+sx_0)}{\left[\mathcal  M_s^\frac{p}{pn+1}(\|f\|_1,\|g\|_1)\right]^\frac{1}{pn+1}}=\frac{ g(c_0x+x_0)}{ \|g\|_1^\frac{1}{pn+1}}=\frac{ f(x)}{
		\|f\|_1^\frac{1}{pn+1}},\ \ \ x\in K.
	\end{equation}
	
	Now, let $x_1,x_2\in K$ be arbitrarily fixed elements. Let $y_2:=\psi(x_2)=c_0x_2+x_0\in S$. By (\ref{Dubuc-egyenloseg}) we have that
	$$g(y_2)=g(c_0x_2+x_0)=\left(\frac{\|g\|_1}{\|f\|_1}\right)^\frac{1}{pn+1}f(x_2).$$
	Let $z:=(1-s)x_1+sy_2=(1-s)x_1+sc_0x_2+sx_0\in Z_s(x_1,y_2)$; if we denote  $\tilde x=\frac{1-s}{1-s+sc_0}x_1+\frac{sc_0}{1-s+sc_0}x_2$, then $\tilde x\in K$ and $z=(1-s+sc_0)\tilde x+sx_0.$ Applying again (\ref{Dubuc-egyenloseg}), it turns out that 
	$$h(z)=\left[\mathcal  M_s^\frac{p}{pn+1}\left(1,\frac{\|g\|_1}{\|f\|_1}\right)\right]^\frac{1}{pn+1}f(\tilde x).$$
	Replacing now the above relations into (\ref{elso-BBL-feltetel}) for $x_1$ and $y_2$, it follows that 
	{\small \begin{equation}\label{fp-concave}
	\left[\mathcal  M_s^\frac{p}{pn+1}\left(1,\frac{\|g\|_1}{\|f\|_1}\right)\right]^\frac{1}{pn+1}f\left(\frac{1-s}{1-s+sc_0}x_1+\frac{sc_0}{1-s+sc_0}x_2\right )\geq \mathcal M_s^p\left(f(x_1),\left(\frac{\|g\|_1}{\|f\|_1}\right)^\frac{1}{pn+1}f(x_2)\right).
	\end{equation}}
	We distinguish two cases: 
	
\underline{Case 1a:} $p=0$. Note that by (\ref{c0-definicio}) one has $c_0=1$ and relation (\ref{fp-concave}) reduces to 
		$$f\left((1-s)x_1+sx_2\right )\geq \mathcal M_s^0\left(f(x_1),f(x_2)\right),$$
		i.e., $f$ is a $(s,0)$-concave function in $K.$
		
		\underline{Case 1b:} $p\neq 0$. Again by (\ref{c0-definicio}), a simple computation and relation (\ref{fp-concave}) give that  $$f\left(\frac{1-s}{1-s+sc_0}x_1+\frac{sc_0}{1-s+sc_0}x_2\right )\geq \mathcal M_{\frac{sc_0}{1-s+sc_0}}^p(f(x_1),f(x_2)),$$
		i.e., $f$ is a $(t,p)$-concave function in $K$ with $t=\frac{sc_0}{1-s+sc_0}$. 
	
	The rest of the proof of (\ref{Dubuc-fuggvenyek}) follows by (\ref{Dubuc-egyenloseg}). 
	\medskip
	%which is a kind of concavity of $f^p$ on $\tilde A.$
	
%	\underline{Case 2:} $p\in (-\frac{1}{n},0)$. By (\ref{fp-concave}) we have  that  $$f^p\left(\frac{1-s}{1-s+sc_0}x_1+\frac{sc_0}{1-s+sc_0}x_2\right )\leq \frac{1-s}{1-s+sc_0}f^p(x_1)+\frac{sc_0}{1-s+sc_0}f^p(x_2),$$ where
%	$c_0=\left(\frac{\|g\|_1}{\|f\|_1}\right)^\frac{p}{pn+1}.$
%	
%
%	
%	
%	\underline{Case 4:} $p=\infty$. Relation (\ref{fp-concave}) implies that
%	$$f\left(\frac{1-s}{1-s+sc_0}x_1+\frac{sc_0}{1-s+sc_0}x_2\right )\geq \mathcal M_s^{+\infty}\left(f(x_1),f(x_2)\right),$$
%	where $c_0=\left(\frac{\|g\|_1}{\|f\|_1}\right)^\frac{1}{n}.$
	
	\textbf{\underline{Case 2}:} $p=-\frac{1}{n}.$ 
	To treat this case, we need the following H\"older-type inequality
	\begin{equation}\label{Holder-integral}
	\int_{A}\mathcal M_s^{\frac{1}{n}}\left(f_1(x),f_2(x)\right){\rm d}\mathcal L^n(x)\leq  \mathcal M_s^{\frac{1}{n}}\left(\int_{A}f_1(x){\rm d}\mathcal L^n(x),\int_{A}f_2(x){\rm d}\mathcal L^n(x)\right),
	\end{equation}
	where $f_1,f_2:A\to \mathbb R$ are nonnegative, integrable functions on a measurable set $A\subset \mathbb R.$ The proof of (\ref{Holder-integral}) follows by the Newton binomial expansion and the classical H\"older inequality for integrals; indeed, 
	\begin{eqnarray*}
	L&:=&\int_{A}\mathcal M_s^{\frac{1}{n}}\left(f_1(x),f_2(x)\right){\rm d}\mathcal L^n(x)=\int_{A}\left((1-s)f_1^\frac{1}{n}(x)+sf_2^\frac{1}{n}(x)\right)^n{\rm d}\mathcal L^n(x)\\&=&\sum_{k=0}^n\left(\begin{matrix}
		n \\
		k
	\end{matrix}\right)(1-s)^ks^{n-k}\int_{A}f_1^\frac{k}{n}(x)f_2^\frac{n-k}{n}(x){\rm d}\mathcal L^n(x)\\&\leq & \sum_{k=0}^n\left(\begin{matrix}
	n \\
	k
\end{matrix}\right)(1-s)^ks^{n-k}\left(\int_{A}f_1(x){\rm d}\mathcal L^n(x)\right)^\frac{k}{n}\left(\int_{A}f_2(x){\rm d}\mathcal L^n(x)\right)^\frac{n-k}{n}\\&=&\mathcal M_s^{\frac{1}{n}}\left(\int_{A}f_1(x){\rm d}\mathcal L^n(x),\int_{A}f_2(x){\rm d}\mathcal L^n(x)\right).
	\end{eqnarray*}
	Moreover, equality holds in (\ref{Holder-integral}) if and only if for some $c>0$ we have $f_2(x)=cf_1(x)$ for a.e. $x\in A.$
	
	Due to Theorem \ref{Theorem-Riemannian-2}, the  equality in the Borell-Brascamp-Lieb inequality, i.e., $\delta_{\mathbb R^n,s}^{-\frac{1}{n}}(f,g,h)=0,$ is characterized by: 
	\begin{itemize}
		\item[(i)] $\operatorname{supp}h=\psi_s(\operatorname{supp}f)$ up to a null measure set$;$
		\item[(ii)] $h(\psi_s(x))= \mathcal M_s^{-\frac{1}{n}}\left({f(x)},g(\psi(x))\right)=\frac{f(x)}{ {\rm Jac}(\psi_s)(x)}$
		for a.e. $x\in \operatorname{supp}f;$
		\item[(iii)] $\|f\|_1=\|g\|_1.$
	\end{itemize}
	Let us keep the previous notations, i.e., $A=\operatorname{supp}f$, $B=\psi(A)$ and the
	convex function $\eta:\mathbb R^n\to \mathbb R$ with $\psi=\nabla \eta.$ By (ii) we have that $f(x)=h(\psi_s(x)){ {\rm Jac}(\psi_s)(x)}$ for a.e. $x\in A$, thus $\|f\|_1=\|h\|_1.$ Moreover, by (iii) and the Monge-Amp\`ere equation (\ref{Monge-Ampere})
	it follows that $f(x)=g(\psi(x)){\rm Jac}(\psi)(x)$ for a.e. $x\in A.$ In particular, 
	\begin{equation}\label{B-mertek}
	\mathcal L^n(B)=\mathcal L^n(\psi(A))=\int_{\psi(A)}{\rm d}\mathcal L^n=\int_{A}{\rm Jac}(\psi)(x){\rm d}\mathcal L^n(x)=\int_{A}\frac{f(x)}{g(\psi(x))}{\rm d}\mathcal L^n(x).
	\end{equation}
	Since 
	$ \psi_s(A)\subseteq Z_s(A,B)\subseteq \operatorname{supp}h,$ by (i) it follows  that $\mathcal L^n(Z_s(A,B))=\mathcal L^n(\psi_s(A)).$ Therefore, 
	\begin{eqnarray*}
		\mathcal L^n(Z_s(A,B))&=&\mathcal L^n(\psi_s(A))=\int_{\psi_s(A)}{\rm d}\mathcal L^n=\int_{A}{\rm Jac}(\psi_s)(x){\rm d}\mathcal L^n(x)\\&=&\int_{A}f(x)\mathcal M_s^{\frac{1}{n}}\left(\frac{1}{f(x)},\frac{1}{g(\psi(x))}\right){\rm d}\mathcal L^n(x)=\int_{A}\mathcal M_s^{\frac{1}{n}}\left(1,\frac{f(x)}{g(\psi(x))}\right){\rm d}\mathcal L^n(x)\ \ \ \  \  {\rm (see \ (ii))} \\&\leq & \mathcal M_s^{\frac{1}{n}}\left(\int_{A}{\rm d}\mathcal L^n(x),\int_{A}\frac{f(x)}{g(\psi(x))}{\rm d}\mathcal L^n(x)\right)\ \ \ \ \ \ \ \ \  \ \ \ \  \ \ \ \ \ \ \ \ \ \ \ \ \ \ \ \ \  \ \ \ \  \ \ \ \ \ \ \ \ \ \  {\rm (see \ (\ref{Holder-integral}))}  \\&=&\mathcal M_s^{\frac{1}{n}}\left(\mathcal L^n(A),\mathcal L^n(B)\right)\ \ \ \ \ \ \ \ \  \ \ \ \  \ \ \ \ \ \ \ \ \ \ \ \ \ \ \ \ \  \ \ \ \  \ \ \ \ \ \ \ \ \ \ \  \ \ \ \  \ \ \ \ \ \ \ \ \ \ \ \ \ \ \ \ \ \ {\rm (see \ (\ref{B-mertek}))} \\&=&\left((1-s)\mathcal L^n(A)^\frac{1}{n}+s\mathcal L^n(B)^\frac{1}{n})\right)^n\\&\leq &	\mathcal L^n(Z_s(A,B)). \ \ \ \ \  \ \ \ \ \ \ \ \ \ \ \ \ \ \ \ \ \ \ \ \ \ \ \ \ \ \ \ \ \ \ \ \ \ \ \ \ \  {\rm (cf. \ Brunn-Minkowski\ inequality)}
	\end{eqnarray*}
	Consequently, in the latter estimates we necessarily have equalities. First, being equality in the Brunn-Minkowski inequality, the sets $A$ and $B$ are homothetic convex bodies up to a null measure set;  let $K$ and $S$ be the convex bodies which differ from $A$ and $B$ by null sets,  respectively, and $c_0>0$ and $x_0\in \mathbb R^n$ such that
	$S=c_0K+x_0.$ As before,  we may consider the interior points of $K^{\circ}$ and $S^{\circ}$ instead of the sets $K$ and $S$ themselves.
	Second,  by the equality case in (\ref{Holder-integral}), we have  for some $c>0$ that $\frac{f(x)}{g(\psi(x))}=c$ for a.e. $x\in A.$ In particular, by (ii) we have
	\begin{equation}\label{J-c}
	{\rm Jac}(\psi)(x)=c\ \ {\rm and}\ \ {\rm Jac}(\psi_s)(x)=\mathcal M_s^\frac{1}{n}(1,c)\ \ {\rm for\ a.e.}\ x\in A.
	\end{equation}
	It is clear that $c=\frac{\mathcal L^n(B)}{\mathcal L^n(A)}=c_0^n$.
	
	By the convexity of $S$, relation (\ref{J-c}) and the interior regularity result of Caffarelli \cite{Caffarelli}, it turns out that $\eta$ is of class $C^2$ on $K^{\circ}$. Furthermore,  by (\ref{J-c}) we have 
	$${\rm det}^\frac{1}{n}[(1-s)I_n+s{\rm Hess} \eta(x)]=\left(\mathcal M_s^\frac{1}{n}(1,c)\right)^\frac{1}{n}=1-s+sc^\frac{1}{n}=(1-s){\rm det}^\frac{1}{n}[I_n]+s{\rm det}^\frac{1}{n}[{\rm Hess} \eta(x)],\ x\in K^{\circ},$$
	thus the strict concavity of det$^\frac{1}{n}(\cdot)$ on the cone of nonnegative definite symmetric matrices implies that {\rm Hess}$\eta(x)=c_0 I_n$ for every $x\in K^{\circ}$. 
	Accordingly, 
	$$\psi(x)=\nabla \eta(x)=c_0x+x_0\ \ \ {\rm and}\ \ \ \psi_s(x)=(1-s+sc_0)x+sx_0,\ \ x\in K^{\circ}. $$
	By continuity reason, the latter relations hold true for all $x\in K$ and by (ii) we have
	\begin{equation}\label{Dubuc-egyenloseg-2}
	{\mathcal  M_s^\frac{1}{n}(1,c_0^n)}{h((1-s+sc_0)x+sx_0)}=c_0^n{ g(c_0x+x_0)}{ }={ f(x)}
	,\ \ \ x\in K.
	\end{equation}
	
	Let $x_1,x_2\in K$ be two arbitrarily fixed elements. Let $y_2:=\psi(x_2)=c_0x_2+x_0\in S$. By (\ref{Dubuc-egyenloseg-2}) we have that
	$$g(y_2)=g(c_0x_2+x_0)=c_0^{-n}f(x_2).$$
	Let $z:=(1-s)x_1+sy_2=(1-s)x_1+sc_0x_2+sx_0\in Z_s(x_1,y_2)$; if   $\tilde x=\frac{1-s}{1-s+sc_0}x_1+\frac{sc_0}{1-s+sc_0}x_2$, then $\tilde x\in K$ and $z=(1-s+sc_0)\tilde x+sx_0.$ By (\ref{Dubuc-egyenloseg-2}), we have that 
	$$h(z)=\mathcal  M_s^{-\frac{1}{n}}(1,c_0^{-n})f(\tilde x).$$
	Replacing the above expressions into (\ref{elso-BBL-feltetel}),  it follows that 
	\begin{equation*}%\label{fp-concave-2}
	\mathcal  M_s^{-\frac{1}{n}}(1,c_0^{-n})f\left(\frac{1-s}{1-s+sc_0}x_1+\frac{sc_0}{1-s+sc_0}x_2\right )\geq \mathcal M_s^{-\frac{1}{n}}\left(f(x_1),c_0^{-n}f(x_2)\right),
	\end{equation*}
	which is equivalent to 
	$$f^{-\frac{1}{n}}\left(\frac{1-s}{1-s+sc_0}x_1+\frac{sc_0}{1-s+sc_0}x_2\right )\leq \frac{1-s}{1-s+sc_0}f^{-\frac{1}{n}}(x_1)+\frac{sc_0}{1-s+sc_0}f^{-\frac{1}{n}}(x_2),$$
	which means that $f$ is $(t,-\frac{1}{n})$-concave in $K$ with $t=\frac{sc_0}{1-s+sc_0}$. The relations for $g$ and $h$ from (\ref{Dubuc-fuggvenyek}) easily follow by (\ref{Dubuc-egyenloseg-2}).  
	
	(b)$\implies$(a) This implication trivially holds; indeed, the inequality in (\ref{elso-BBL-feltetel}) and the equality in (\ref{BBL-euklidesz}) easily follow by the $(t,p)$-concavity of $\Phi$ and relation (\ref{Dubuc-fuggvenyek}), respectively. 
	\hfill $\square$\\

%	
%	-------------------------------------------------------------------------------------------------------------------------------
%	\begin{remark}\rm 
%		(a) Amit nem ertettunk Udvarhelyen az az volt, hogy a Dubuc cikk elso oldalan levo Theoreme $B_n$, ami epp a $p=-\frac{1}{n}$ esetnek felel meg, hogyan fugg ossze a Theoreme 12-vel (160 oldal a Dubuc cikkben). Nem ertettuk azt sem, hogy a mi eredmenyunk azt mondja ki, hogy $\|f\|_1=\|g\|_1$, ami az o Theorem 12-ben nem jelent meg. Viszont, a Theoreme $B_n$-ben pontosan az van, hogy az $f, g$ es $h$ fuggvenyek integraljai megegyeznek, ami egy valtozocserebol jon ki! Ez nyilvan, pontosan a $p=-\frac{1}{n}$ esetben ervenyes es lathatjuk a Theorem \ref{Theorem-Euklidean-egyenloseg}-bol is, hogy pontosan akkor van Dubuc-tipusu egyutthatonk a $g$ es $h$ elott, amikor $p=-\frac{1}{n}$.  Ehhez kepest, a fotetelben, Theoreme 12, altalaban megjelennek az integralok, ami nalunk is megjelenik a $c_0$ suly formajaban, ami nem mas, mint $c_0=\left(\frac{\mathcal L^n(B)}{\mathcal L^n(A)}\right)^{1/n}=\left(\frac{\|g\|_1}{\|f\|_1}\right)^\frac{p}{pn+1},$ amikor $p>-\frac{1}{n}.$  
%		
%		(b) Szinten a Dubuc cikk: a Theorem 12-nel NEM mond semmit a $\phi$ fuggvenyrol (mint ahogyan tette a Theorem $B_n$-ben, ahol a $\phi$ konvex volt), hanem ez helyett felir egy altalanos egyenlotlenseget (161 oldal), ami epp a mi $(t,p)-$konkavitasunknak felel meg megfelelo $t$ sullyal tekintve. Nekunk azert is jo a helyzet, mert effektiv megmondjuk a megfelelo sulyokat is pl. a Rossi-Salani cikkhez kepest.  
%	\end{remark}
%	
%	
%	HATRA TEVE
	
		Although our approach is more appropriate for characterizing equality cases, we conclude the present section by stating  weak stability results for Brunn-Minkowski-type inequalities, e.g. for the dimension-free-Brunn-Minkowski inequality (\ref{log-BM});  an exhaustive study of the latter inequality can be found in  B\"or\"oczky,  Lutwak,  Yang and Zhang \cite{BLYZ}.   
%		 For simplicity,  we shall state a quantitative result for the log-Brunn-Minkowski inequality  in $\mathbb R^n$ as an application of Theorem \ref{Theorem-Riemannian};  an exhaustive study of this inequality can be found in  B\"or\"oczky,  Lutwak,  Yang and Zhang \cite{BLYZ}. 
	
%	\begin{corollary}\label{log-BM-stability} {\bf (Quantitative log-Brunn-Minkowski inequality)}
%		Let $n\geq 2$ and $s\in (0,1)$. For every nonempty compact sets $A,B\subset \mathbb R^n$ with $V(A)\neq 0\neq V(B)$ we have  
%		$$\delta_s^0(A,B):= \frac{	V((1-s)A+sB) }{ V(A)^{1-s}V(B)^s}-1\geq n\tilde s\frac{\left|V(A)^\frac{\tilde s}{n}-V(B)^\frac{\tilde s}{n}\right|^\frac{1}{\tilde s}}{(1-s)V(A)^\frac{1}{n}+sV(B)^\frac{1}{n}},$$
%		where $\tilde s=\min(s,1-s).$
%	\end{corollary}

	\begin{proposition}\label{proposition-stability}  {\rm \textbf{(Quantitative $p$-Brunn-Minkowski inequality in $\mathbb{R}^n$)}}
		Let $n\geq 2$, $s\in (0,1)$ and $p\geq -\frac{1}{n}$. For every nonempty compact sets $A,B\subset \mathbb R^n$ with $\mathcal L^n(A)\neq 0\neq \mathcal L^n(B)$ we have  
		\begin{equation}\label{quantit-egy}
		\delta_s^p(A,B):=\frac{\mathcal L^n((1-s)A+sB)}{\mathcal M_s^\frac{p}{1+pn}\left(\mathcal L^n(A),\mathcal L^n(B)\right)}-1\geq G_s^{p,n}\left(1,1,{\mathcal L^n(B)},{\mathcal L^n(A)}\right).
		\end{equation}
		In particular, the quantitative dimension-free-Brunn-Minkowski inequality reads as $$\delta_s^0(A,B):= \frac{	\mathcal L^n((1-s)A+sB) }{ \mathcal L^n(A)^{1-s}\mathcal L^n(B)^s}-1\geq n\tilde s\frac{\left|\mathcal L^n(A)^\frac{\tilde s}{n}-\mathcal L^n(B)^\frac{\tilde s}{n}\right|^\frac{1}{\tilde s}}{(1-s)\mathcal L^n(A)^\frac{1}{n}+s\mathcal L^n(B)^\frac{1}{n}},$$
		where $\tilde s=\min(s,1-s).$ 
		 Moreover, %$\delta_s^0(A,B)=0$ if and only if $A$ and $B$ are translates. 
		 \begin{itemize}
		 	\item[(i)] $\delta_s^{+\infty}(A,B)=0$   if and only if $A$ and $B$
		 	are homothetic  convex  bodies  up to a null measure set$;$ 
		 	\item[(ii)] if $p<+\infty$, then $\delta_s^{p}(A,B)=0$  if and only if $A$ and $B$
		 	are translated  convex  bodies up to a null measure set. 
		 \end{itemize}
	\end{proposition}

\begin{remark}\rm 
We point out that Proposition \ref{proposition-stability} is a simple consequence of Theorem \ref{Theorem-Euklidean-egyenloseg} rather than an independent
verification of the equality
case in the classical Brunn-Minkowski inequality. Indeed, in the proof of Theorem \ref{Theorem-Euklidean-egyenloseg}  we used the characterization of equality cases in the Brunn-Minkowski inequality. 
\end{remark}

	{\it Proof of Proposition \ref{proposition-stability}.} Let $f=	\mathbbm{1}_A$, $g=\mathbbm{1}_B$  and $h=\mathbbm{1}_{(1-s)A+sB}$; then  
	$$\delta_{\mathbb R^n,s}^p(f,g,h)=\frac{	\mathcal L^n((1-s)A+sB) }{ \mathcal M_s^\frac{p}{1+pn}\left(\mathcal L^n(A),\mathcal L^n(B)\right)}-1.$$
	On the other hand, for a.e. $x\in A$, we have 
	\begin{eqnarray*}
	G_s^{p,n}\left(\frac{ f(x)}{v_{1-s}(\psi(x),x)},\frac{ g(\psi(x))}{v_s(x,\psi(x))},\frac{1}{\|f\|_1},\frac{1}{\|g\|_1}\right)&=&G_s^{p,n}\left(1,1,\frac{1}{\mathcal L^n(A)},\frac{1}{\mathcal L^n(B)}\right)\\&=&G_s^{p,n}\left(1,1,{\mathcal L^n(B)},{\mathcal L^n(A)}\right).
	\end{eqnarray*}
	It remains to apply Theorem \ref{Theorem-Riemannian} and relation (\ref{Gap-homogen}) to conclude the proof of (\ref{quantit-egy}).  
	
	Moreover, if $\delta_s^{p}(A,B)=0$ for some $p\geq -\frac{1}{n}$,  by Theorem \ref{Theorem-Euklidean-egyenloseg} (more precisely, by (\ref{transzlacio})) we have that $A$  and $B$ are convex (up to a null measure set) and there exists $x_0\in \mathbb R$ such that  $B=c_0A+x_0,$  where $c_0=\left(\frac{\mathcal L^n(B)}{\mathcal L^n(A)}\right)^\frac{1}{n}.$  In particular, if $p<+\infty$, by (\ref{Dubuc-fuggvenyek}) it turns out that $\mathbbm{1}_B(c_0x+x_0)=c_0^\frac{1}{p}\mathbbm{1}_A(x)$ for a.e. $x\in A$, which implies that $c_0=1.$ The converse is trivial. 
	\hfill $\square$
	
%	\medskip
%	
%	\textcolor{red}{Erdekesnek tartom hogy egyenlÃ¶seg kell fennaljon a mertekekben ha $p<\infty$ esetben 
%	egyenlÃ¶seget akarunk itt es igy csak transzlaciok engedelyeztek; ugyanakkor a $p= \infty$ esetben az egyenlÃ¶seg megenged homotetiakat is. Ismert mar ez az eredmeny vagy uj? Mit gondolsz? } 
%	\medskip
	
	\begin{remark}\rm 
%		The first part of the quantitative log-Brunn-Minkowski inequality (see Corollary \ref{log-BM-stability}) directly follows by (\ref{quantit-egy}) for $p=0;$ 	the equality case yields by Corollary \ref{Cor-Euclidean-2}(ii), see below. 	
	Note that the right hand side of (\ref{quantit-egy}) measures the  difference between the volumes $\mathcal L^n(A)$ and $\mathcal L^n(B)$ whenever $p<+\infty$. For $p=+\infty$, inequality (\ref{quantit-egy}) reduces precisely to the usual Brunn-Minkowski inequality (\ref{BM-1}) since $G_s^{+\infty,n}\left(1,1,{\mathcal L^n(B)},{\mathcal L^n(A)}\right)=0,$  thus no meaningful stability can be obtained in this case.  
	\end{remark}

		\section{Equality in Borell-Brascamp-Lieb inequality: Riemannian case}\label{section-Riemannian-0} 
		%\subsubsection{Equality in}
		
		In subsection \ref{subsection-4.1}  we shall discuss the consequences of the equality case in Borell-Brascamp-Lieb inequality in the Riemannian setting, while in subsection \ref{subsection-4.2} we characterize the equality case in the distorted Brunn-Minkowski inequality. 
		
		\subsection{Curvature rigidity}\label{subsection-4.1}

		  We begin this section with an important notation to be used in the sequel. For every $k\in \mathbb R$, let ${\bf s}_{k}:[0,\infty)\to \mathbb
		R$ be the function defined by
		\begin{equation*}
		{\bf s}_{k}(r)=\left\{
		\begin{array}{lll}
		\displaystyle\frac{\sinh(\sqrt{-k}r)}{\sqrt{-k}r} & \hbox{if} & {k}<0,\\
		1
		& \hbox{if} &  {k}=0, \\
		\displaystyle\frac{\sin(\sqrt{k}r)}{\sqrt{k}r} & \hbox{if} & {k}>0,
		\end{array}\right.\ \ \ \ r>0.
		\end{equation*}
		By taking the limit $r\to 0$, one may choose ${\bf s}_k(0)=1$. 
		
	Let $(M,w)$ be a complete $n$-dimensional Riemannian manifold with  Ricci curvature {\rm Ric}$(M)\geq (n-1)k$  for some $k\in \mathbb R.$ Let  $s\in (0,1)$, $p\geq -\frac{1}{n}$ and $f,g,h:M\to [0,\infty)$ be three nonzero, compactly supported integrable functions with $\operatorname{supp}f=A$ and $\operatorname{supp}g=B$, verifying 
		\begin{equation}\label{gorbulet}
		h(z)\geq \mathcal M^{p}_s
		\left(\left(\frac{{\bf s}_{k}(d(x,y))}{{\bf s}_{k}((1-s)d(x,y))}\right)^{n-1}{f(x)},\left(\frac{{\bf s}_{k}(d(x,y))}{{\bf s}_{k}(sd(x,y))}\right)^{n-1}{g(y)} \right)
		\end{equation}
		for all $(x,y)\in A\times B, z\in Z_s(x,y).$ Since {\rm Ric}$(M)\geq (n-1)k$, Bishop's comparison principle implies that for every $x\in M$, $y\in M\setminus {\rm cut}(x)$ and $s\in (0,1)$, 
		\begin{equation}\label{vess}
		v_s(x,y)\geq \left(\frac{{\bf s}_{k}(sd(x,y))}{{\bf s}_{k}(d(x,y))}\right)^{n-1},
		\end{equation}
		see e.g. Bishop and Crittenden \cite{Bishop-Crittenden}, and  Cordero-Erausquin, McCann and Schmuckenschl\"ager
		\cite[Corollary 2.2]{CMS}.   Here, cut$(x)\subset M$ denotes the cut locus of $x\in M,$ which is a null set of $M$, see Sakai \cite[Lemma III. 4.4 (c)]{Sakai}. 
	The estimate (\ref{vess}) and assumption  (\ref{gorbulet}) imply through the monotonicity of $\mathcal M_s^p(\cdot,\cdot)$ the validity of (\ref{ConditionRescaledBBLWithWeights-vege}). Consequently, Theorem \ref{Theorem-Riemannian} implies that $$\delta_{M,s}^p(f,g,h)\geq 0,$$
	which is precisely Corollary 1.1 in Cordero-Erausquin, McCann and Schmuckenschl\"ager
	\cite{CMS}.
	
	Within the aforementioned geometric setting  we establish the following rigidity result  appearing whenever the Borell-Brascamp-Lieb deficit vanishes.

	\begin{theorem}\label{Theorem-rigiditas} {\bf (Curvature rigidity; Riemannian case)} 
	Under the above assumptions,	if $$\delta_{M,s}^p(f,g,h)=0$$ then for a.e. 
		$x\in \operatorname{supp} f=A$ one has$:$
		\begin{itemize}
			\item[(i)] the sectional curvature is equal to the constant $k$ 	
			 along the geodesic  $ t\mapsto \psi_t(x),$ $t\in [0,1];$ 
			\item[(ii)] if $p>-\frac{1}{n}$ and $d_x=d(x,\psi(x))$, then $$\frac{h(\psi_s(x))}{\left[\mathcal M_s^\frac{p}{pn+1}(\|f\|_1,\|g\|_1)\right]^\frac{1}{pn+1}}=\left(\frac{{\bf s}_{k}(d_x)}{{\bf s}_{k}((1-s)d_x)}\right)^{n-1}\frac{ f(x)}{
				\|f\|_1^\frac{1}{pn+1}}=\left(\frac{{\bf s}_{k}(d_x)}{{\bf s}_{k}(sd_x)}\right)^{n-1}\frac{ g(\psi(x))}{ \|g\|_1^\frac{1}{pn+1}}.$$
			\item[(iii)] if $p=-\frac{1}{n}$ and $d_x=d(x,\psi(x))$, then $\|f\|_1=\|g\|_1$ and $$	h(\psi_s(x))= \mathcal M^{-\frac{1}{n}}_s
			\left(\left(\frac{{\bf s}_{k}(d_x)}{{\bf s}_{k}((1-s)d_x)}\right)^{n-1}{f(x)},\left(\frac{{\bf s}_{k}(d_x)}{{\bf s}_{k}(sd_x)}\right)^{n-1}{g(\psi(x))} \right).$$
		\end{itemize}
	\end{theorem}
	
%	\textcolor{red}{ Can we say something about the geodesic convexity of $A,B$ and $(t,p)$-convexity of the functions $f,g,h$ similar to Euclidean case? I think that would be the analogue of Dubuc's result in the 
%	constant sectional curvature case. }
	
	{\it Proof.} 
	Assume that the Borell-Brascamp-Lieb deficit vanishes, i.e. $\delta_{M,s}^p(f,g,h)= 0.$  By Remark \ref{remark-4-14}, one has that 
	\begin{equation}\label{kell-vegere}
	{h(\psi_s(x))}=\mathcal M^p_s
	\left(\frac{f(x)}{v_{1-s}(\psi(x),x)},\frac{g(\psi(x))}{v_s(x,\psi(x))} \right)\ \ {\rm for\ a.e.}\ x\in A.
	\end{equation}
	On the other hand, by relations (\ref{gorbulet})-(\ref{kell-vegere}) and the monotonicity of $\mathcal M_s^p(\cdot,\cdot)$ we have for a.e. $x\in A$ that
	\begin{eqnarray*}%\label{eziskell}
		h(\psi_s(x))&\geq& \mathcal M^{p}_s
		\left(\left(\frac{{\bf s}_{k}(d_x)}{{\bf s}_{k}((1-s)d_x)}\right)^{n-1}{f(x)},\left(\frac{{\bf s}_{k}(d_x)}{{\bf s}_{k}(sd_x)}\right)^{n-1}{g(\psi(x))} \right)\\&\geq&
		\mathcal M^p_s
		\left(\frac{f(x)}{v_{1-s}(\psi(x),x)},\frac{g(\psi(x))}{v_s(x,\psi(x))} \right)\\&=&h(\psi_s(x)).
	\end{eqnarray*}
	Consequently, we have equalities in the above estimates. Again, by the 
	monotonicity of $\mathcal M_s^p(\cdot,\cdot)$ we necessarily have for a.e. $x\in A$ that
	\begin{equation}\label{eziskellegyenloseg}
	\left(\frac{{\bf s}_{k}((1-s)d_x)}{{\bf s}_{k}(d_x)}\right)^{n-1}=v_{1-s}(\psi(x),x)\ \  {\rm and}\ \  \left(\frac{{\bf s}_{k}(sd_x)}{{\bf s}_{k}(d_x)}\right)^{n-1}=v_s(x,\psi(x)),
	\end{equation}
	which proves (ii)\&(iii) through Theorem \ref{Theorem-Riemannian-egyenloseg} (b)(iii) and Theorem \ref{Theorem-Riemannian-2}, respectively. 
	
	If $Y(s)=d(\exp_x)_{-s\nabla \varphi(x)}$ denotes the   differential of the exponential map at $-s\nabla \varphi(x)\in T_xM,$  relation (\ref{eziskellegyenloseg}) implies in particular that  for a.e. $x\in A$, 
	$$v_s(x,\psi(x))=\frac{{\rm det}[Y(s)]}{{\rm det}[Y(1)]}=\left(\frac{{\bf s}_{k}(sd_x)}{{\bf s}_{k}(d_x)}\right)^{n-1}.$$
	 By using the equality case in the comparison principle of Bishop and Crittenden \cite[\S 11.10]{Bishop-Crittenden}, the latter relation implies that for a.e. $x\in A$  the sectional curvature  along the geodesics $ t\mapsto \psi_t(x),$ $t\in [0,1]$ is constant, having its value $k$;   the detailed proof is given in Cordero-Erausquin, McCann and Schmuckenschl\"ager
	 \cite[Corollary 2.2]{CMS}. 
	\hfill $\square$
	
	\begin{remark}\rm Theorem \ref{Theorem-rigiditas} complements both Th\`eor\'em 1 from Cordero-Erausquin \cite{Cordero-CRAS} and Corollary 2.2 from Cordero-Erausquin, McCann and Schmuckenschl\"ager
		\cite{CMS} where the Pr\'ekopa-Leindler inequalities are considered  (i.e., $p=0$). 
	\end{remark}
	
	\subsection{Equality in distorted Brunn-Minkowski inequality}\label{subsection-4.2}

Let $(M,w)$ be a complete $n$-dimensional Riemannian manifold with  {\rm Ric}$(M)\geq (n-1)k$  for some $k\in \mathbb R.$ It is well known that the Borell-Brascamp-Lieb inequality (\ref{BBL-eredeti-Riemann}) implies the distorted Brunn-Minkowski inequality $(\ref{BM-eredeti})$ (see e.g. Bacher \cite{Bacher}), i.e., for every compact sets $A,B\subset M$ and $s\in (0,1)$ one has 
$$\textsf{m}(Z_s(A,B))^\frac{1}{n}\geq \tau_{1-s}^{k,n}(\Theta_{A,B})\textsf{m}(A)^\frac{1}{n}+\tau_{s}^{k,n}(\Theta_{A,B})\textsf{m}(B)^\frac{1}{n},$$ 
where $\Theta_{A,B}$ is given by (\ref{theta-nak}). The latter inequality follows by choosing 
	   \begin{equation}\label{f-g-h}
	   h:=\mathbbm{1}_{Z_s(A,B)}, \ \ f:=\left(\frac{{\bf s}_k((1-s)\Theta_{A,B})}{{\bf s}_k(\Theta_{A,B})}\right)^{n-1}\mathbbm{1}_A\ \ {\rm and}\ \  g:=\left(\frac{{\bf s}_k(s\Theta_{A,B})}{{\bf s}_k(\Theta_{A,B})}\right)^{n-1}\mathbbm{1}_B,
	 \end{equation}
	 which verify inequality (\ref{gorbulet}) for $p=+\infty$.
	  Due to the definition of $\Theta_{A,B}$ and monotonicity properties of $\mathcal M_s^{+\infty}(\cdot,\cdot)$ and the function $r\mapsto \frac{{\bf s}_{k}(r)}{{\bf s}_{k}(sr)}$, respectively, we obtain $\delta_{M,s}^{+\infty}(f,g,h)\geq 0$. Since  $\tau_s^{k,n}(\theta)=s\left(\frac{{\bf s}_{k}(s\theta)}{{\bf s}_{k}(\theta)}\right)^{1-\frac{1}{n}}$, the non-negativity of the Borell-Brascamp-Lieb deficit is equivalent to  the distorted Brunn-Minkowski inequality (\ref{BM-eredeti}). We notice that the same choice for $f,g$ and $h$ also provide for every $p\geq -\frac{1}{n}$ that
	  %\begin{equation}\label{BBL-BM}
	  $$
	  \textsf{m}(Z_s(A,B))\geq \mathcal M_s^{\frac{p}{1 + np}}\left(\left(\frac{{\bf s}_k((1-s)\Theta_{A,B})}{{\bf s}_k(\Theta_{A,B})}\right)^{n-1}\textsf{m}(A),\left(\frac{{\bf s}_k(s\Theta_{A,B})}{{\bf s}_k(\Theta_{A,B})}\right)^{n-1}\textsf{m}( B)\right).
	  $$
	  %\end{equation}
	  For $p=+\infty$ the latter inequality reduces to the above distorted Brunn-Minkowski inequality.

	In the sequel, we shall provide a  complete characterization of the equality in the distorted Brunn-Minkowski inequality $(\ref{BM-eredeti})$.  To do this, we recall that a set  $A\subset M$ contains a cut locus pair if there exist $x,y\in A$ such that $x$ belongs to the cut locus of $y$. As usual, $A\subset M$ is geodesic convex if every two points of $A$ can be joined by a unique minimizing geodesic whose image belongs entirely to $A.$

		\begin{theorem}\label{Theorem-Riemannian-rigiditas-gorbulet-CD} {\bf (Equality in distorted Brunn-Minkowski  inequality)} 	Let $(M,w)$ be a complete $n$-dimensional Riemannian manifold,  
			%		with Ricci curvature {\rm Ric}$(M)\geq (n-1)k$ for some $k\in \mathbb R$.   
			$A,B\subset M$ be compact sets with %\textcolor{red}{with nonempty interior} 
			$ {\rm \sf{m}}(A)\neq 0\neq {\rm \sf{m}}(B)$  
			and $s\in (0,1)$.	
			Then the following statements hold$:$ 
			\begin{itemize}
				\item[(i)] {\rm (Positively curved case)} If ${\rm Ric}(M)\geq (n-1)k$ for some $k>0$, equality holds in $(\ref{BM-eredeti})$ if and only if  
				$Z_s(A,B)=A=B$
				up to a null measure set$;$ moreover, if the sets $A$ and $B$ do not contain cut locus pairs,   equality holds in $(\ref{BM-eredeti})$ if and only if  there exists an open, geodesic convex set in $M$ which differs from $A$ and $B$ by a null set$;$  
				%up to a null measure set $A=B$ and these sets are geodesics convex
				\item[(ii)] {\rm (Negatively curved case)} If $(M,g)$ has nonpositive, nonzero sectional curvature and ${\rm Ric}(M)\geq (n-1)k$ for some $k<0$, equality cannot hold in $(\ref{BM-eredeti});$
				%			\item[(b)]  $\rho_s=\mathcal \mathcal M_s^\frac{p}{pn+1}(\|f\|_1,\|g\|_1)h;$
				\item[(iii)]  {\rm (Null curved case)} Let ${\rm Ric}(M)\geq 0$ and $\pi: \tilde{M} \to M$  be the universal covering of $M$. Assume that the sets $A$ and $B$ are small enough and sufficiently close to each other in the sense that there exist open sets $\tilde{U} \subseteq \tilde{M}$ and $U\subseteq M$ such that every geodesic segment with ends in the sets $A$ and  $B$ is unique and belongs to $U$, and 		
				 $\pi:\tilde{U} \to U$ is a homeomorphism. Then  equality holds in $(\ref{BM-eredeti})$ if and only if $\tilde{U}$ is isometrically identified with an open subset of $\R^n$ and  $\pi^{-1}(A)$ and $\pi^{-1}(B)$ are convex sets up to null measure sets which are homothetic to each other in $\mathbb R^n$.
%					$($The latter conclusion is relevant whenever 
%				%${\rm \sf{m}}(A\triangle B)\neq 0$; 
%				$A\neq B$ up to a null measure set$;$
%				otherwise, the optimal transport map is the identity, $\Theta_{A,B}=0$, and equality in  $(\ref{BM-eredeti})$ trivially holds if and only if $Z_s(A,B)=A=B$ up to a null measure set$).$
			\end{itemize}
		\end{theorem}
	
	\begin{remark}\rm  (a) Let us note that the different nature of statements (i) and (ii) in the above theorem is due to the fact that the definition of $\Theta_{A,B}$ changes according to the sign of the lower bound of the 
	Ricci curvature. In this sense it is not expected that (i) is a particular case of (ii). 
	
%	(b) The lack of cut locus pairs in the sets $A$ and $B$ in the second part of (i) is crucial in our argument.  In particular, when $M=\mathbb S^n$ is the standard unit round sphere and 
%	there exists $x\in \mathbb S^n$ such that $x$ and $-x$ are interior points of $A$ or $B$, one can easily prove that equality holds in $(\ref {BM-eredeti})$ if and only if $A=B=\mathbb S^n.$ 
	
 	(b) 
 %By using covering maps and the fact that a certain region of the covering space $\tilde M$  has null sectional curvature 
 %(coming from the equality in (\ref{BM-eredeti})), the sets $\pi^{-1}(A)$ and $\pi^{-1}(B)$ in (iii) can be considered subsets of %$\mathbb R^n$; see the proof for details.
  The assumptions that  $A$ and $B$ are small enough and sufficiently close to each other are crucial for the third statement. Indeed, let us consider the cylinder $M=\mathbb S^1\times \mathbb R \subseteq \R^3$ with the induced Euclidean metric and two (small) congruent curvilinear rectangles $A$ and $B$ in the opposite sides of the cylinder. Then $\pi^{-1}(A)$ and $\pi^{-1}(B)$ are congruent rectangles in $\tilde M=\mathbb R^2$, and we have a strict inequality  in $(\ref{BM-eredeti})$ since ${\rm \sf{m}}(A)= {\rm \sf{m}}(B)=\frac{{\rm \sf{m}}(Z_{1/2}(A,B))}{2}.$ 
	\end{remark}

	{\it Proof of Theorem \ref{Theorem-Riemannian-rigiditas-gorbulet-CD}.} 
Let us suppose that $A$ and $B$ are two compact subsets of $M$ and $s\in (0,1)$  such that equality holds in $(\ref{BM-eredeti})$. As we shall see, the most difficult part will be to prove the statement about the geodesic convexity of $A$ and $B$ in part (i).

In the sequel, let us assume that we have equality in  
	(\ref{BM-eredeti}), i.e., 
	\begin{equation}\label{egyeloseg-BM}
	\textsf{m}(Z_s(A,B))^\frac{1}{n}= \tau_{1-s}^{k,n}(\Theta_{A,B})\textsf{m}(A)^\frac{1}{n}+\tau_{s}^{k,n}(\Theta_{A,B})\textsf{m}(B)^\frac{1}{n}.
\end{equation}
Moreover, by  Theorem \ref{Theorem-rigiditas} (ii), we also have for a.e. $x\in A$ that
\begin{equation}\label{A-B-s}
1=\frac{\left(\frac{{\bf s}_{k}(d_x)}{{\bf s}_{k}((1-s)d_x)}\right)^{n-1}}{\left(\frac{{\bf s}_k(\Theta_{A,B})}{{\bf s}_k((1-s)\Theta_{A,B})}\right)^{n-1}}=\frac{\left(\frac{{\bf s}_{k}(d_x)}{{\bf s}_{k}(sd_x)}\right)^{n-1}}{\left(\frac{{\bf s}_k(\Theta_{A,B})}{{\bf s}_k(s\Theta_{A,B})}\right)^{n-1}}.
\end{equation}

%\textcolor{red}
{	
 	(i) (Positively curved case)  Two cases are distinguished.}
	 
%\textcolor{red}
\textbf{\underline{{Case 1}}}: {	 $A\cap B\neq \emptyset.$ Clearly, by (\ref{theta-nak}) we have  $\Theta_{A,B}=0.$ Therefore, due to the monotonicity of $r\mapsto \frac{{\bf s}_{k}(r)}{{\bf s}_{k}(sr)}$, relation  (\ref{A-B-s}) and $\Theta_{A,B}=0$ give that $d_x=d(x,\psi(x))=0$ for a.e. $x\in A$. Thus, $\psi(x)=x$ for a.e. $x\in A$ which implies that $B=A$ up to a null measure set. Thus, (\ref{egyeloseg-BM}) reduces to $\textsf{m}(Z_s(A,B))=\textsf{m}(A)=\textsf{m}(B)$. Let $S=A\cap B.$ It is clear that $\textsf{m}(S)=\textsf{m}(A).$ By the definition of the $s$-intermediate set $Z_s$, we have that  $S\subseteq Z_s(S,S)\subseteq Z_s(A,B)$. Moreover,  $\textsf{m}(Z_s(A,B)\setminus S)=\textsf{m}(Z_s(A,B))-\textsf{m}( S)=0,$ i.e. $Z_s(A,B)$ is equal to $A\cap B$ up to a null measure set. }
	
%\textcolor{red}
\textbf{\underline{{Case 2}}}: {  $A\cap B= \emptyset.$  By the monotonicity of  $r\mapsto \frac{{\bf s}_{k}(r)}{{\bf s}_{k}(sr)}$ and (\ref{A-B-s}) we have 
	\begin{equation}\label{dx-kell}
	d_x=d(x,\psi(x))=\Theta_{A,B}=\min\{d(x,y):x\in A,y\in B\}>0\ {\rm for\ a.e.}\ x\in A.
	\end{equation}
	 For simplicity of notation, let $t_0:=\Theta_{A,B}$ and $$B_{t_0}=\{x\in M:\textrm{there exists}\ y\in B\ \textrm{such that}\ d(x,y)<t_0\}=\bigcup_{y\in B}B(y,t_0)$$ be the $t_0$-neighborhood of $B$. }	
%\textcolor{red}
It is clear that $A\cap B_{t_0}=\emptyset$. Indeed, if we assume that $x\in A\cap B_{t_0}$, then there exists $y\in B$ such that $d(x,y)<t_0$, which contradicts the fact that $t_0=\Theta_{A,B}$. 

Now, let us fix $x\in A$ such that $d(x,\psi(x))=t_0$; due to (\ref{dx-kell}), the latter happens for a.e. $x\in A$.  
	By construction, we have that $B(\psi(x),t_0)\subset {\rm int}B_{t_0}=B_{t_0}$, thus $B(\psi(x),t_0)\cap A=\emptyset.$ Fix $r_0\in (0,t_0)$. Then, for every $0<r<r_0$ let us fix $z_r\in Z_{\frac{r}{2t_0}}(x,\psi(x))$; then  $B(z_r,\frac{r}{2})\subset B(x,r)\cap B(\psi(x),t_0).$ Therefore, $B(z_r,\frac{r}{2})\subset B(x,r)\setminus A$, i.e., $A$ is $\frac{1}{2}$-porous at $x$.   Since $x\in A$ is arbitrarily fixed and porous sets have zero measure (see e.g. Rajala \cite{Rajala}) it follows that $A$ has null measure,  $\textsf{m}(A)=0$, which contradicts our assumption, proving the first part of the assertion.

Now,  we assume  the sets $A$ and $B$ do not contain cut locus pairs and   (\ref{egyeloseg-BM}) holds.   By Cases 1\&2 we know that  $A$, $B$ and $Z_s(A,B)$ coincide up to a null measure set. Accordingly, without loss of generality we may consider the case that $Z_s(A,A)=A\cup C$ where ${\sf{m}}(C)=0$. The proof of the geodesic convexity of $A$ (up to a null measure set) is divided into several steps.  Before to do this, let $A_*$ be the density one set of $A$. Clearly, $A_*\subset \overline A= A$ by the closedness of $A$ and ${\sf{m}}(A\setminus A_*)=0$
by Lebesgue's theorem.

{\it \underline{Claim 1}:} $Z_s(A_*,A_*)\subseteq A.$ %where $A_*$ denotes the density set of $A$. 

%\textcolor{blue}{	
	Let $x,y\in A_*$ be arbitrarily fixed; we shall prove  that $\{z\}= Z_s(x,y)\subset A.$ Note that $z\in Z_s(x,y)$ is unique since $x\notin {\rm cut}(y)$. Moreover, the latter fact also implies  that there are neighborhoods $U$ and $V$ of $x$ and $y$, respectively, such that $x'\notin {\rm cut}(y')$ for every $(x',y')\in U\times V$. Clearly, we may choose $U:=B(x,\frac{1}{m})$ and $V:=B(y,\frac{1}{m})$ for $m\in \mathbb N$ sufficiently large. Let $A_x^m=A\cap B(x,\frac{1}{m})$ and $A_y^m=A\cap B(y,\frac{1}{m})$. Since $x,y\in A_*$, it follows that $${\sf{m}}(A_x^m)\geq \frac{1}{2}{\sf{m}}( B(x,\frac{1}{m}))\ \ {\rm and}\ \  {\sf{m}}(A_y^m)\geq \frac{1}{2}{\sf{m}}( B(y,\frac{1}{m}))$$ for $m\in \mathbb N$ sufficiently large. Thus, by the Brunn-Minkowski inequality (\ref{BM-eredeti}) applied to the sets $A_x^m$ and $A_y^m$ in $M$, and by using the fact that $\tau_s^{k,n}\geq s$ for every $s\in (0,1),$ we have for sufficiently large $m$ that 
	\begin{equation}\label{bm-becsles}
	{\sf{m}}(Z_s(A_x^m,A_y^m))^\frac{1}{n}\geq (1-s){\sf{m}}(A_x^m)^\frac{1}{n}+s{\sf{m}}(A_y^m)^\frac{1}{n}>0.
	\end{equation}
	Since $Z_s(A,A)=A\cup C$ with ${\sf{m}}(C)=0$, the estimate (\ref{bm-becsles}) shows that $Z_s(A_x^m,A_y^m)$
	contains a positively measured subset of $A$. Therefore, for every $m\in \mathbb N$ large enough,  let us choose such a triplet $(x_m,y_m,z_m)$ with $x_m\in A_x^m$, $y_m\in A_y^m$ and $\{z_m\}=Z_s(x_m,y_m)\subset A;$   the element $z_m$ is also uniquely determined since $x_m\notin {\rm cut}(y_m)$. 
%}

%\textcolor{blue}{
	We shall prove that the sequence $(z_m)_m$ converges to $z$  (up to a subsequence) and $z\in A$.  Since $M$ is compact (following by the Bonnet-Myers theorem) and   $z_m\in A$ for every $m\in \mathbb N$, there exists $\tilde z\in M$ such that $\lim_{m\to \infty}z_m=\tilde z\in \overline A=A$. It remains to prove that $\tilde z=z$.  By $Z_s(x_m,y_m)=\{z_m\}$, we have that   $d(x_m,z_m) = s d(x_m,y_m)$ and $
	d(z_m,y_m) = (1-s) d(x_m,y_m)$. Taking the limit as $m\to \infty$, it follows that $	d(x,\tilde z) = s d(x,y) $ and $
	d(\tilde z,y) = (1-s) d(x,y)$, i.e. $\tilde z\in Z_s(x,y)$. By uniqueness, we have $\tilde z=z$, which concludes the proof of Claim 1. 
%}

{\it \underline{Claim 2}:} $A_*$ is open. 

This statement can be seen as a curved version of the Steinhaus theorem, see \cite{Steinhaus}. First, let us observe that $(A_*)_*=A_*.$ Indeed, since $A_*\subseteq A$,  the inclusion $(A_*)_*\subseteq A_*$ is trivial. Conversely, if we assume that there exists $x\in A_*\setminus (A_*)_*$, it follows that for every $\varepsilon>0$ sufficiently small there exists $r_\varepsilon>0$ such that for every $0<r<r_\varepsilon$ we have ${\sf{m}}(A\cap B(x,r))\geq (1-\varepsilon){\sf{m}}( B(x,r))$ and ${\sf{m}}(A_*\cap B(x,r))\leq (\eta+\varepsilon){\sf{m}}( B(x,r))$ for some $\eta\in [0,1)$. Therefore, one has 
$$0={\sf{m}}((A\setminus A_*)\cap B(x,r))={\sf{m}}(A\cap B(x,r))-{\sf{m}}(A_*\cap B(x,r))\geq (1-\eta-2\varepsilon){\sf{m}}( B(x,r)),$$
a contradiction. 

Let $p\in A_*=(A_*)_*$ and fix $r>0$ such that $B(p,2r)$ is a totally normal neighborhood of $p$.  First, let us assume that $\frac{1}{2}\leq s<1.$ We introduce the function $R_p: B(p,r)\to  B(p,r)$ which associates to each $x\in B(p,r)$ the point $R_p(x)$ by reflecting $x$ through $p$ such that $p\in Z_s(x,R_p(x))$. We notice that   $R_p(x)\in B(p,r)$ (since $\frac{1}{2}\leq s<1$) and the point $R_p(x)$ is uniquely determined, i.e., $R_p$ is well defined. 
Fix $0<\delta<r$ sufficiently small that will specified later; performing the same construction for every $q\in B(p,\delta)$ instead of $p$, we defined the function $R_q:B(p,r)\to B(p,r+\delta)$ such that 
\begin{equation}\label{q-kesobb-majd}
q\in Z_s(x,R_q(x))\  {\rm for\ all}\ x\in B(p,r).
\end{equation} 
Since $p\in (A_*)_*$, for every $\varepsilon>0$ sufficiently small there exists $r_\varepsilon>0$ such that for every $0<r<r_\varepsilon$, we have  ${\sf{m}}(A_*\cap B(p,r))\geq (1-\varepsilon){\sf{m}}( B(p,r)).$ By the Borel regularity of the measure ${\sf{m}}$ one can find a compact set $K\subset A_*\cap B(p,r)$ such that ${\sf{m}}(K)\geq (1-2\varepsilon){\sf{m}}( B(p,r)).$ Now, choose $\delta<r$ so small that  
\begin{equation}\label{q2-kesobb-majd}
R_q(K)\subset B(p,r)\ \ 
 {\rm and} \ \ {\sf{m}}(R_q(K))\geq \frac{1-s}{2s}{\sf{m}}(K) \ {\rm for\ all}\ q\in B(p,\delta).
 \end{equation}
  The inclusion $R_q(K)\subset B(p,r)$ follows by a continuity reason. In order to verify the inequality in (\ref{q2-kesobb-majd}), let us observe first that $R_q(x)=\exp_q\circ R\circ \exp_q^{-1}(x)$, $x\in B(p,r)$, where $R:T_qM\to T_qM$ is the $s$-reflection given by $R(y)=-\frac{1-s}{s}y,$ $y\in T_qM.$ Since $\exp_q$ is a diffeomorphism on $B(p,r)$ and $d(\exp_q)_0={\rm id}$, the map $\exp_q$ is a local bi-Lipschitz map with bi-Lipschitz constant arbitrarily close to 1, which concludes the proof of  (\ref{q2-kesobb-majd}). 
  
  With this choice of $\delta>0$, we shall prove that $B(p,\delta)\subset A.$ By contradiction, let us assume that there exists $q\in B(p,\delta)$ such that $q\notin A$. We notice that there is no $x\in K$ such that $R_q(x)\in K$. Indeed, by contrary, we would have that $x\in A_*$ and  $R_q(x)\in A_*$, thus by (\ref{q-kesobb-majd}) and Claim 1 we get $q\in Z_s(x,R_q(x))\subset Z_s(A_*,A_*)\subseteq A$,  which contradicts  $q\notin A$. Therefore, for every $x\in K$ one has that $R_q(x)\notin K$, i.e., $K\cap R_q(K)=\emptyset$. On the other hand, since $K\cup R_q(K)\subseteq B(p,r)$, by (\ref{q2-kesobb-majd}) we have 
  $${\sf{m}}(B(p,r))\geq {\sf{m}}(K)+{\sf{m}}(R_q(K))\geq \left(1+\frac{1-s}{2s}\right){\sf{m}}(K)\geq \left(1+\frac{1-s}{2s}\right) (1-2\varepsilon){\sf{m}}( B(p,r)),$$
  a contradiction. Accordingly, $B(p,\delta)\subset A.$ Since $B(p,\delta)$ is open, one has that $B(p,\delta)=B(p,\delta)_*\subseteq A_*$.  
  
 The case  $0<s<\frac{1}{2}$ works similarly by interchanging $(s,1-s)$ with $(1-s,s)$. Accordingly, the function  $R_p: B(p,r)\to  B(p,r)$  will be defined by reflecting  $x$ through $p$ with the property that $p\in Z_{1-s}(x,R_p(x))=Z_s(R_p(x),x)$ (instead of $p\in Z_s(x,R_p(x))$); the same should be performed in (\ref{q-kesobb-majd}) for $R_q$, $q\in B(p,\delta)$, i.e.,  $q\in Z_s(R_q(x),x)$.

{\it \underline{Claim 3}:} $Z_s(A_*,A_*)\subseteq A_*.$

%\textcolor{red}{	
	Since $A_*$ is open (cf. Claim 2), one can prove that  $ Z_s(A_*,A_*)$ is also open. Indeed, let $z\in Z_s(A_*,A_*)$ be fixed arbitrarily. Then there exists $x,y\in A_*$ such that $\{z\}= Z_s(x,y)$. Let $V\subset A_*$ be an open neighborhood of $y$. Due to the lack of cut locus pairs in $A$, the map $\exp_x:\exp_x^{-1}(V)\to V$ is a diffeomorphism. Therefore, the set $U=\exp_x(s\exp_x^{-1}(V))$ is open,  $U=Z_s(x,V)\subset Z_s(A_*,A_*)$ and $z=\exp_x(s\exp_x^{-1}( y))\in U$. Accordingly, by Claim 1 one has $Z_s(A_*,A_*)=Z_s(A_*,A_*)_*\subseteq A_*.$  
%	In fact, since the reverse inclusion trivially holds, we have equality. 
%}

{\it \underline{Claim 4}:}  $A_*$ is geodesic convex. 

%\textcolor{red}{
 Let $x,y\in A_*$ $(x\neq y)$, and  $d_0:=d(x,y)$. Since $x\notin {\rm cut}(y)$, let $\gamma:[0,1]\to M$ be the unique minimal geodesic  joining $x$ and $y$, parametrized 
	proportionally to arc-length. Since $A_*$ is open, there exists $\delta>0$ such that $B(x,\delta)\cup B(y,\delta)\subset A_*$. If $\delta\geq d_0$, we have nothing to prove, since ${\rm Im}(\gamma)\subset B(x,\delta)\subset A_*$. If $\delta< d_0$, let $s_0<\delta$ and  let $I=\overline B(x,s_0)\cap {\rm Im}(\gamma)$ and $J=\overline B(y,s_0)\cap {\rm Im}(\gamma)$ be two geodesic segments in $\gamma$ with lengths $s_0,$ i.e., $I=\gamma([0,\frac{s_0}{d_0}])$ and $J=\gamma([1-\frac{s_0}{d_0},1])$. Hereafter,  $\overline B(x,r)=\{y\in M:d(x,y)\leq r\},$ $r>0.$ By the minimality of $\gamma$, we clearly have that $Z_s(I,J)\subset {\rm Im}(\gamma)$; more precisely, by the parametrization we have that $Z_s(I,J)=\gamma([s(1-\frac{s_0}{d_0}),(1-s)\frac{s_0}{d_0}+s])$ and its length is $s_0$.  Moreover, since $I\cup J\subset A_*$, by Claim 3 we also have that $Z_s(I,J)\subseteq A_*$. Repeating this argument, we cover the whole geodesic segment $\gamma$ after finitely many steps with such pieces of geodesic segments of length $s_0$, all of them belonging to $A_*$. 
	
	%Thus $A_*$ is geodesic convex. 
%}

%The proof of (i) is concluded since the measure of $A_*\subset A$ is equal to the measure of $A$. 
	
	\medskip 
	
%\textcolor{red}
%\textcolor{red}
{	(ii) (Negatively curved case)
Due to (\ref{theta-nak}),  one has $\Theta_{A,B}=\max\{d(x,y):x\in A,y\in B\}>0.$ Similarly as above, relation   (\ref{A-B-s}) implies that 
\begin{equation}\label{dx-kell-2}
d_x=d(x,\psi(x))=\Theta_{A,B}=:t^0\ {\rm for\ a.e.}\ x\in A.
\end{equation}
The proof  is 'dual' to (i); for completeness, we provide it. Let  $$B^{t^0}=\bigcup_{y\in B}(M\setminus \overline B(y,t^0)).$$ 
 Since $t^0>\inf_{x\notin B}\max_{y\in B}d(x,y)$, it turns out that $\bigcap_{y\in B}\overline B(y,t^0)\neq \emptyset$; thus $B^{t^0}$ is a proper open subset of $M.$ }
 
%\textcolor{red}

	We claim that $A\cap  B^{t^0}=\emptyset$; indeed, if $x\in A\cap  B^{t^0}$, it follows that there exists $y\in B$ such that $x\in M\setminus \overline B(y,t^0)$, i.e., $d(x,y)>t^0,$ which contradicts the definition of $t^0=\Theta_{A,B}$.

	According  to (\ref{dx-kell-2}), for a.e. $x\in A$, one has  $d(x,\psi(x))=t^0$ and $x\notin {\rm cut}(\psi(x))$; let us choose such an $x\in A$.   
	It is clear that $M\setminus \overline B(\psi(x),t^0)\subset {\rm int}B^{t^0}=B^{t^0}$, thus $(M\setminus \overline B(\psi(x),t^0))\cap A=\emptyset.$ Since $x\notin {\rm cut}(\psi(x))$,  we may extend the minimal geodesic joining the point $\psi(x)$ to $x$ beyond $x$ such that the extended geodesic is still minimizing between $\psi(x)$ and points in a small neighborhood of $x$. Let $z_r\in M$ be such a point belonging to the extended geodesic with $d(z_r,x)=\frac{r}{2}$ for sufficiently small $r>0$; thus, $d(z_r,\psi(x))=d(z_r,x)+d(x,\psi(x))=\frac{r}{2}+t^0.$ This construction shows that
	$B(z_r,\frac{r}{2})\subset B(x,r)$ and $B(z_r,\frac{r}{2})\subset M\setminus \overline B(\psi(x),t^0)$, i.e., $B(z_r,\frac{r}{2})\subset B(x,r)\setminus A$, which means that $A$ is  $\frac{1}{2}$-porous at $x$. Consequently,  one has $\textsf{m}(A)=0$, which contradicts our assumption.

	 \medskip 
	
 	(iii) (Null curved case) Let $\pi: \tilde{M}\to M$ be the universal covering of $M$, see Boothby \cite[Corollary 9.8]{Boothby}. We consider the pull-back metric on $\tilde{M}$ such that $\pi$ becomes a local isometry.

Let  $A$ and $B$ be two sets in $M$  which are small enough, sufficiently close to each other as in the assumption and $\textsf{m}(A)\neq 0\neq \textsf{m}(B)$; let $U$ and $\tilde U$ the sets from the statement of the theorem. Since ${\rm Ric}(M)\geq 0$ (thus $k=0$)  the  equality in (\ref{BM-eredeti}) reads as 
\begin{equation}\label{BM-0}
\textsf{m}(Z_s(A,B))^\frac{1}{n}= (1-s)\textsf{m}(A)^\frac{1}{n}+s\textsf{m}(B)^\frac{1}{n}.
\end{equation}
 By Theorem \ref{Theorem-rigiditas} (ii), applied to $f$, $g$ and $h$ from  (\ref{f-g-h}), it turns out that $Z_s(A,B)=\psi_s(A)$ up to a null measure set and by  Theorem \ref{Theorem-rigiditas} (i) we have that for a.e. $x\in A$ the sectional curvature is zero  along the geodesic  $ t\mapsto \psi_t(x),$ $t\in [0,1].$   Let $A_0\subset A$ be such that at any point of $A_0$ the above property holds, i.e., for every $x\in A_0$ the sectional curvature is zero  along the geodesic  $ t\mapsto \psi_t(x),$ $t\in [0,1]$; clearly, ${\sf{m}}(A)={\sf{m}}(A_0).$  Thus $C:=\pi^{-1}(\{\psi_t(x):x\in A_0,t\in [0,1]\})\subset \pi^{-1}(U)=\tilde U$ which is isometric to a proper subset of $\mathbb R^n$ endowed with the usual Euclidean metric.  In fact, the map  $\pi: \tilde{U} \to U= \pi(\tilde{U})$ is an isometry and the sets $\pi^{-1}(A),$ $\pi^{-1}(B)$ and $\pi^{-1}(Z_s(A,B))=\pi^{-1}(\psi_s(A))$ are subsets of $C$ up to null measure sets. By the isometric property of the covering map $\pi: \tilde{U} \to U$ and relation (\ref{BM-0}) it turns out that 
\begin{equation}\label{BM-13}
\mathcal L^n(\pi^{-1}(Z_s(A,B)))^\frac{1}{n}= (1-s)\mathcal L^n(\pi^{-1}(A))^\frac{1}{n}+s\mathcal L^n(\pi^{-1}(B))^\frac{1}{n}.
\end{equation}
Note that 
\begin{equation}\label{inkluzio}
(1-s)\pi^{-1}(A)+s\pi^{-1}(B)\subseteq \pi^{-1}(Z_s(A,B)).
\end{equation}
 Indeed, if $\tilde a\in \pi^{-1}(A)$ and $\tilde b\in \pi^{-1}(B)$ are arbitrarily fixed, then the geodesic segment $t\mapsto (1-t)\tilde a+t\tilde b\subset \mathbb R^n$, $t\in [0,1]$,  is mapped by $\pi$ to the (minimal) geodesic  segment $t\mapsto \pi((1-t)\tilde a+t\tilde b)\subset M$, $t\in [0,1]$, joining the points $\pi(\tilde a)\in A$ and $\pi(\tilde b)\in B$; thus, $\pi((1-s)\tilde a+s\tilde b)\in Z_s(A,B)$, which concludes the proof of (\ref{BM-13}). By  (\ref{BM-13}), (\ref{inkluzio}) and the usual Brunn-Minkowski inequality (\ref{BM-1}), we necessarily obtain that
$$ \mathcal L^n((1-s)\pi^{-1}(A)+s\pi^{-1}(B))^\frac{1}{n}= (1-s)\mathcal L^n(\pi^{-1}(A))^\frac{1}{n}+s\mathcal L^n(\pi^{-1}(B))^\frac{1}{n}.$$
Therefore, by Proposition \ref{proposition-stability}(i) it follows that  the sets $\pi^{-1}(A)$ and $\pi^{-1}(B)$ are homothetic  convex  bodies  from  which  sets  of  measure  zero  have  been  removed. 
%	and equality holds in  $(\ref{BM-eredeti})$ then $(M,w)$ is a.e. flat between $A$ and $B$ in the sense that for a.e. $x\in A$ the sectional curvature vanishes along the geodesics  $ t\mapsto \psi_t(x),$ $t\in [0,1],$ where $t\mapsto\psi_t$ denotes the optimal transport map from the measure  $\mu=\mathbbm{1}_A/{\rm \sf{m}}(A){\rm d}V_w$ to $\nu =\mathbbm{1}_B/{\rm \sf{m}}(B){\rm d}V_w$.
	\hfill $\square$

	\begin{remark}\rm 
		 Let 
	 $\mu=\mathbbm{1}_A/{\rm \sf{m}}(A){\rm d}\sf{m}$ and $\nu =\mathbbm{1}_B/{\rm \sf{m}}(B){\rm d}\sf{m}$ be the measures from the proof of  Theorem \ref{Theorem-Riemannian-rigiditas-gorbulet-CD} and $\psi:M\to M$ be the optimal transport map between them. Then we generically have the two-sided estimate for the Wasserstein distance between $\mu$ and $\nu$; namely, 
	 \begin{equation}\label{Wasserstein}
	 (\Theta_{A,B}^{\rm min})^2\leq \mathcal W(\mu,\nu):=\int_A d^2(x,\psi(x)){\rm d}\mu(x)\leq (\Theta_{A,B}^{\rm max})^2,
	 \end{equation}
	where $$\Theta_{A,B}^{\rm min}=\min\{d(x,y):x\in A,y\in B\}\ \ {\rm and}\ \ \Theta_{A,B}^{\rm max}=\max\{d(x,y):x\in A,y\in B\}.$$
 The proofs of (i)/(ii) in Theorem \ref{Theorem-Riemannian-rigiditas-gorbulet-CD} correspond to the equality cases in the two inequalities of (\ref{Wasserstein}), appearing in the   positively/negatively curved settings.  
For instance, when $A$ and $B$ are two disjoint positive measure sets, the equality at the left hand side cannot hold. Indeed, if we push-forward the initial measure $\mu$ only with $\Theta_{A,B}^{\rm min}$, we cannot reach the target measure $\nu$; the reason is that  the transport cost  $(\Theta_{A,B}^{\rm min})^2$ is not enough to realize this transportation. {A similar explanation works also in the 'dual' case (ii); in this setting, such an equality cannot be realized since by pushing forward the measure $\mu$ to $\nu$ the  transport cost $(\Theta_{A,B}^{\rm max})^2$ is too large.   
}

%\textcolor{red}{	(b) {\bf Cylinder.  Ezt meg vegig kell  gondolni a (iii) szempontjabol. Indukalt gorbulet=0, es ha veszunk egy nagyon kicsi illetve egy nagyon nagy gombot a hengeren (utobbi 'felcsavarodik' jo sokszor a hengerre), akkor nem igazan varhato el az un. homotetikussag a ket halmaz kozott. Nyilvan, ebben az esetben a sokasag nem simply connected, mivel $\pi_1(Cylinder)=\mathbb Z$ es vannak cut pontok is.}}

%\textcolor{red}{Lehet hogy az alabbi  b) reszt mar akar ki is lehetne hagyni, valoszinu hogy azelÃ¶tt let ideirva mielÃ¶tt meg a Dubuc uj bizonyitasa bekerÃŒlt a cikkbe. Mit szolsz?}

%(b) 	The second part of (i) can be seen as the curved counterpart of Dubuc's result \cite{Dubuc}, see also  Figalli and Jerison \cite{Figalli-Jerison-3}. Indeed, in \cite[Th\'eor\`eme 6]{Dubuc} it is proved that for a measurable set $A\subset \mathbb R^n$ one has 
%$\mathcal L^n((1-s)A+sA)= \mathcal L^n(A)$ if and only if $A$ is an open convex set up to a null measure set. A minor modification of our proof provides an alternative 
%way to demonstrate Dubuc's result in $\mathbb R^n$ due to the lack of cut locus points; in fact, the Bonnet-Myers theorem should be avoided in Claim 1 which can be done in a trivial manner.  

%It remains to consider the same problem on generic $n$-dimensional compact Riemannian manifolds with ${\rm Ric}(M)\geq (n-1)k$ for some $k>0$. 

%In particular, if the Ricci curvature {\rm Ric}$(M)\geq 0$, we have that 
%	$$
%	\textsf{m}(Z_s(A,B))\geq\mathcal M_s^{\frac{p}{1 + np}}(\textsf{m}(A),\textsf{m}( B)).
%$$
		\end{remark}

	\medskip 
	
%	\begin{corollary}\label{elozo}
%		Let $(M,w)$ be a complete $n$-dimensional Riemannian manifold with Ricci curvature {\rm Ric}$(M)\geq 0$ and $s\in (0,1)$, $p\geq-\frac{1}{n}$. Then for every nonempty open bounded sets $A,B\subset M$ one has 
%		\begin{equation}\label{volume-BM-2}
%		{\rm Vol}(Z_s(A,B))\geq\mathcal M_s^{\frac{p}{1 + np}}({\rm Vol}(A),{\rm Vol}( B)).
%		\end{equation}
%		Moreover, if equality holds in $(\ref{volume-BM-2})$, the sectional curvature vanishes along the geodesic  $[0,1]\in t\mapsto \psi_t(x)$  for a.e. $x\in A$; in addition, if $p<+\infty$, then ${\rm Vol}(A)={\rm Vol}( B)$. 
%	\end{corollary}
	
%	{\it Proof.} 
%	%	Since Ric$(M)\geq 0$, we have that $v_s(x,y)\geq 1$ for every $x,y\in M$.  
%	The proof follows by Theorem \ref{Theorem-rigiditas} by choosing the indicator functions	
%	$f=\mathbbm{1}_A$, $g=\mathbbm{1}_B$ and $h=\mathbbm{1}_{Z_s(A,B)}$  of the sets $A$, $B$ and $Z_s(A,B),$ respectively.  
%	\hfill $\square$\\

Let $(M,w)$ be a complete $n$-dimensional Riemannian manifold with nonnegative Ricci curvature and $s\in (0,1)$. Then for every nonempty open bounded sets $A,B\subset M$ one has 
	\begin{equation}\label{volume-BM}
	{\sf{m}}(Z_s(A,B))^\frac{1}{n}\geq (1-s){\sf{m}}(A)^\frac{1}{n}+s{\sf{m}}(B)^\frac{1}{n},
	\end{equation}
	which is a particular form of (\ref{BM-eredeti}) for $k=0$. 
We conclude this section by characterizing the equality in (\ref{volume-BM}) via the flatness of the manifold; namely, we have
	
	\begin{corollary}\label{corollary-elso} {\rm \textbf{(Equality in Brunn-Minkowski inequality vs flatness)}} Under the above assumptions, we have: 
	\begin{itemize}
		\item[(i)] if for any points $x,y\in M$ there exist two open sets $A,B\subset M$ with   $x\in A$ and $y\in B$ and verifying the equality  in $(\ref{volume-BM})$,  then $(M,w)$ is flat$;$
		\item[(ii)] if $(M,g)$ is simply connected, equality holds in $(\ref{volume-BM})$ for arbitrary  geodesic balls $A=B(x,r)$ and $B=B(y,R)$ if and only if $(M,w)$ is isometric to $\mathbb R^n.$
	\end{itemize}
	\end{corollary}
	
	{\it Proof}. (i) Fix $x,y\in M$ arbitrarily and assume that we have equality in $(\ref{volume-BM})$ for some open sets $A,B\subset M$ with $x\in A$ and $y\in B$.  Let $\psi:A\to B$ be the optimal transport map from the measure $\mu=\mathbbm{1}_A/{\sf{m}}(A){\rm d}\sf{m}$ to $\nu=\mathbbm{1}_B/{\sf{m}}(B){\rm d}\sf{m}$.  
	By Theorem \ref{Theorem-Riemannian-rigiditas-gorbulet-CD} (iii), the sectional curvature is zero along the geodesics  $ t\mapsto \psi_t(x)$, $t\in[0,1]$,   joining a.e. $x\in A$ to $\psi(x)\in B$.  The arbitrariness of the points $x,y$  and a density argument shows that the sectional curvature on $(M,w)$ is zero. 
	
	(ii) If $(M,w)$ is isometric to $\mathbb R^n$, we have equality in  $(\ref{volume-BM})$ for every balls. Conversely, if $(M,w)$ is simply connected, the equality case in $(\ref{volume-BM})$ for geodesic balls implies that $(M,w)$ has zero sectional curvature (from (i)). By the  Killing-Hopf theorem (see, e.g., do Carmo \cite[Theorem 4.1]{doCarmo}) it follows that $(M,w)$  is isometric to  $\mathbb R^n$.     \hfill $\square$\\

	\medskip
	
	\section{Equality in Borell-Brascamp-Lieb inequality: Finsler case}\label{section-5-0}
	
{
	 Let $M$ be a connected
	$n$-dimensional smooth manifold and $TM=\bigcup_{x \in M}T_{x}
	M $ be its tangent bundle. The pair $(M,F)$ is a \textit{Finsler
		manifold} if the continuous function $F:TM\to [0,\infty)$ satisfies
	the conditions
}
	
	(a) $F\in C^{\infty}(TM\setminus\{ 0 \});$
	
	(b) $F(x,tv)=tF(x,v)$ for all $t\geq 0$ and $(x,v)\in TM;$
	%i.e., $F$ is absolutely homogeneous of degree one;
	
	(c) $g_v:=g_{ij}(x,v)=[\frac12F^{2}%
	]_{v^{i}v^{j}}(x,v)$ is positive definite for all $(x,v)\in
	TM\setminus\{ 0 \}.$
	
	\noindent If $F(x,tv)=|t|F(x,v)$ for all $t\in \mathbb R$ and $(x,v)\in TM,$ then
	$(M,F)$ is a  reversible Finsler manifold. A Finsler manifold $(M,F)$ is a:  
	\begin{itemize}
		\item {\it Riemannian manifold}, whenever  $g_{ij}(x)=g_{ij}(x,v)$ is independent of $v.$
	
		\item {\it locally Minkowski space}, if  there
		exists a local coordinate system $(x^i)$ on $M$ with induced tangent
		space coordinates $(y^i)$ such that $F$ depends only on
		$v=v^i{\partial}/{\partial x^i}$ and not on $x.$ 
		
		\item {\it Minkowski space}, whenever $M$ is a finite dimensional vector space (identified by $\mathbb R^n$) which is endowed by a Minkowski norm, inducing a Finsler metric on $\mathbb R^n$ by	translations. 
		\item {\it Berwald space}, whenever the coefficients  $\Gamma_{ij}^{k}(x,v)$ of the Chern connection $R^v$ are independent of $v$. It is clear that
		Riemannian manifolds and $($locally$)$ Minkowski spaces are Berwald 		spaces.
	\end{itemize}

Let $\sigma: [0,r]\to M$ be a piecewise smooth curve. The value $%
L_F(\sigma)= \displaystyle\int_{0}^{r} F(\sigma(t), \dot\sigma(t))\,{\text d}%
t $ denotes the \textit{integral length} of $\sigma.$ For
$x_1,x_2\in M$,
denote by $\Lambda(x_1,x_2)$ the set of all piecewise $C^{\infty}$ curves $%
\sigma:[0,r]\to M$ such that $\sigma(0)=x_1$ and $\sigma(r)=x_2$.
The \textit{metric function} $d_{F}: M\times M
\to[0,\infty)$ is defined by
\begin{equation}  \label{quasi-metric}
d_{F}(x_1,x_2) = \inf_{\sigma\in\Lambda(x_1,x_2)} L_F(\sigma).
\end{equation}
 A
$C^{\infty}$-curve $\sigma:[0,l] \to M$ is called a
\emph{geodesic} if it is locally $d_F$-minimizing and has a constant
speed (i.e., $F(\sigma,\dot{\sigma})$ is constant). $(M,F)$ is forward (resp. backward) complete if any geodesic $\sigma:[0,l] \to  M$ can be extended to $[0,\infty)$ (resp. $(-\infty,l]$). 
The forward and backward metric balls with center $x\in M$ and radius $r>0$ are defined by $B^+(x,r)=\{y\in M:d_F(x,y)<r\}$ and $B^-(x,r)=\{y\in M:d_F(y,x)<r\}$, respectively.  For two linearly independent vectors
$v,w \in T_xM$ and $\mathcal S=\mathrm{span}\{v,w\}$, the {\it flag
	curvature} of the \emph{flag} $(\mathcal S;v)$ is defined by
\[
K(\mathcal S;v) :=\frac{g_v(R^v(w,v)v,w)}{F(v)^2 g_v(w,w) -
	g_v(v,w)^{2}}.
\]
If $(M,F)$ is Riemannian, then the flag curvature reduces to the
sectional curvature which depends only on $\mathcal S$ (not on the
choice of $v \in \mathcal S$). For further concepts and results from Finsler geometry (as Ricci curvature and mean covariation) we refer to Bao, Chern and Shen \cite{BCS}, Krist\'aly \cite{Kristaly-JGA}, Ohta \cite{Ohta} and Shen \cite{Shen}.

	Given $\mu$ and $\nu$ two absolutely continuous measures on $(M,F)$ w.r.t. the Finsler measure ${\sf{m}}_F$ with compact support, there exists a unique optimal transport map from $\mu$ to $\nu$ of the form $\psi(x)=\exp_x(\boldsymbol{\nabla}(-\varphi(x)))$, where $\varphi:M\to \mathbb R$ is a $d_F^2/2$-concave function and $\boldsymbol{\nabla}$ is the Finslerian gradient on $M,$ see Ohta \cite[Theorem 4.10]{Ohta}.  For $s\in (0,1)$ fixed, let  $\psi_s(x)=\exp_x(s\boldsymbol{\nabla}(-\varphi(x)))$ be the $s$-intermediate optimal transport map. The key tool to prove Borell-Brescamp-Lieb inequalities on Finsler manifolds is the Jacobian inequality
	\begin{equation}\label{Jacobian-inequality-Finsler}
	{\rm Jac}(\psi_s)(x)\geq \mathcal M_s^\frac{1}{n}(v_{s}^>(x,\psi(x)),v_{s}^<(x,\psi(x)){\rm Jac}(\psi)(x))\ \ {\rm for\ a.e.}\ x\in \operatorname{supp}(\mu),
	\end{equation}
	where ${\rm Jac}(\psi_s)(x)$ and ${\rm Jac}(\psi)(x)$ are the Jacobian determinant of $\psi_s$ and $\psi$ at $x$ and 
	$$v_s^>(x,y)=\lim_{r\to 0}\frac{{\sf{m}}_F(Z_s(B^-(x,r),y))}{{\sf{m}}_F(B^-(x,(1-s)r))}\ \ {\rm and}\ \ v_s^<(x,y)=\lim_{r\to 0}\frac{{\sf{m}}_F(Z_s(x,B^+(y,r)))}{{\sf{m}}_F(B^+(x,sr))}, $$
	 see Ohta \cite[Proposition 5.3]{Ohta}. 
	 
 Let $(M,F)$ be a forward geodesically complete, $n$-dimensional Finsler manifold with vanishing mean covariation  and Ricci curvature ${\rm Ric}_F(v)\geq\allowbreak \left(n-1\right)k$ for every unit vector $v\in TM$ and some $k\in \mathbb R.$ Let  $s\in (0,1)$, $p\geq -\frac{1}{n}$ and  
	 $f,g,h:M\to [0,\infty)$ be three nonzero, compactly supported integrable functions with $\operatorname{supp}f=A$ and $\operatorname{supp}g=B$, verifying 
	 \begin{equation}\label{gorbulet-2}
	 h(z)\geq \mathcal M^{p}_s
	 \left(\left(\frac{{\bf s}_{k}(d_F(x,y))}{{\bf s}_{k}((1-s)d_F(x,y))}\right)^{n-1}{f(x)},\left(\frac{{\bf s}_{k}(d_F(x,y))}{{\bf s}_{k}(sd_F(x,y))}\right)^{n-1}{g(y)} \right)
	 \end{equation}
	 for all $(x,y)\in A\times B, z\in Z_s(x,y).$
	 Then it is know (see Ohta \cite[Corollary 9.4]{Ohta}) that 
	 $$\delta_{M,s}^p(f,g,h)\geq 0.$$

	 Following  the arguments from \S \ref{section-2} one can formulate in a natural way the Finslerian counterparts of Theorems \ref{Theorem-Riemannian} \& \ref{Theorem-Riemannian-egyenloseg}. We shall state without proof the Finslerian counterpart of Theorem \ref{Theorem-rigiditas}; we leave the details to the interested reader.

		\begin{theorem}\label{Theorem-rigiditas-Finsler} {\bf (Curvature rigidity; Finsler case)} 
	Under the above assumptions,	if  $$\delta_{M,s}^p(f,g,h)=0$$ then for a.e. 
		$x\in \operatorname{supp} f=A$, one has
		\begin{itemize}
			\item[(i)] the flag curvature is equal to the constant $k$ along the geodesic  $ t\mapsto \psi_t(x),$ $t\in [0,1]$, for flags having the form $\{\mathcal S,v\}$ with $\mathcal S={\rm
				span}\{u,v\}\subset T_{\psi_t(x)}M$ and $v=\frac{\rm d }{{\rm d}t}\psi_t(x);$

			\item[(ii)] if $p>-\frac{1}{n}$ and  $d_x=d_F(x,\psi(x))$, then $$\frac{h(\psi_s(x))}{\left[\mathcal M_s^\frac{p}{pn+1}(\|f\|_1,\|g\|_1)\right]^\frac{1}{pn+1}}=\left(\frac{{\bf s}_{k}(d_x)}{{\bf s}_{k}((1-s)d_x)}\right)^{n-1}\frac{ f(x)}{
				\|f\|_1^\frac{1}{pn+1}}=\left(\frac{{\bf s}_{k}(d_x)}{{\bf s}_{k}(sd_x)}\right)^{n-1}\frac{ g(\psi(x))}{ \|g\|_1^\frac{1}{pn+1}}.$$
			\item[(iii)] if $p=-\frac{1}{n}$ and $d_x=d_F(x,\psi(x))$, then $\|f\|_1=\|g\|_1$ and $$	h(\psi_s(x))= \mathcal M^{-\frac{1}{n}}_s
			\left(\left(\frac{{\bf s}_{k}(d_x)}{{\bf s}_{k}((1-s)d_x)}\right)^{n-1}{f(x)},\left(\frac{{\bf s}_{k}(d_x)}{{\bf s}_{k}(sd_x)}\right)^{n-1}{g(\psi(x))} \right).$$
		\end{itemize}
	\end{theorem}

  \begin{remark}\rm 
  	In Theorem \ref{Theorem-rigiditas-Finsler} (i) we have information only on the flag curvature in specific directions of the flag, and not necessarily for any flag direction.  When $(M,F)$ is Riemannian, Theorem \ref{Theorem-rigiditas-Finsler} reduces to Theorem \ref{Theorem-rigiditas}. 
  \end{remark}

	Let $(M,F)$ be a  forward geodesically complete $n$-dimensional Berwald space with nonnegative Ricci curvature and $s\in (0,1)$. Then for every nonempty open bounded sets $A,B\subset M$ one has 
	\begin{equation}\label{volume-BM-3}
	{\sf{m}}_F(Z_s(A,B))^\frac{1}{n}\geq (1-s){\sf{m}}_F(A)^\frac{1}{n}+s{\sf{m}}_F(B)^\frac{1}{n},
	\end{equation}
	see Ohta \cite{Ohta}. 
	
	\begin{corollary}\label{corollary-Finsler-1}{\bf (Brunn-Minkowski inequality on Berwald spaces)}
		Under the above assumptions, if for any points $x,y\in M$ there exist two open sets $A,B\subset M$ with   $x\in A$ and $y\in B$ that verify the equality in $(\ref{volume-BM-3})$,  then $(M,F)$ is a locally Minkowski space. 
	\end{corollary}

	{\it Proof.} 
%	Since the Ricci curvature is nonnegative, we have that $v_s^<\geq 1$ and $v_s^>\geq 1$. Moreover,  since every Berwald space has vanishing mean covariation,
%	see Shen \cite[Propositions 2.6 \& 2.7]{Shen}, we may apply Theorem
%	\ref{Theorem-rigiditas-Finsler}. Thus, 
%	 (\ref{volume-BM-3}) follows by the first part of Theorem \ref{Theorem-rigiditas-Finsler} by choosing $p=+\infty$ and the indicator functions	
%	$f=\mathbbm{1}_A$, $g=\mathbbm{1}_B$ and $h=\mathbbm{1}_{Z_s(A,B)}$  of the sets $A$, $B$ and $Z_s(A,B),$ respectively.
As in Corollary \ref{corollary-elso}, one can prove by means of  Theorem \ref{Theorem-rigiditas-Finsler} (i) that the flag curvature is identically zero (being zero for every choice of the flag). Since $(M,F)$ is a Berwald space, the vanishing of the flag curvature implies that $(M,F)$ is locally Minkowski, see Bao, Chern and Shen \cite[Section 10.5]{BCS}. 
%	Since $v_s^<=v_s^>\equiv 1$ (see (\ref{volume-dist-egyenlo})), we have in addition that ${\sf{m}}_F(B^+(x,r)=\omega_n r^n$ for every $x\in M$ and $r>0$, see Shen \cite{Shen}; in particular, $(M,F)$ is isometric to an $n$-dimensional Minkowski space. 
%	
%	The converse easily follows, taking into account the definition of Minkowski spaces, the particular form of (\ref{volume-form}) and Corollary \ref{Cor-Euclidean-2}, respectively. 
\hfill $\square$

\begin{example} \rm  On $\mathbb R^{n-1}$ ($n\geq 2$) we introduce a
complete Riemannian metric $w$ such that $(\mathbb R^{n-1},w)$ has
nonnegative Ricci curvature, and for every $\varepsilon\geq 0,$ we
define on $\mathbb R^{n}=\mathbb R^{n-1}\times \mathbb R$ the metric
$F_\varepsilon:T\mathbb R^{n}=\mathbb R^{2n}\to [0,\infty)$ for
every $(x,t)\in \mathbb R^{n}$ and $ (y,v)\in T_{x}\mathbb
R^{n-1}\times  T_t\mathbb R=\mathbb R^{n}$ by
$$F_\varepsilon((x,t),(y,v))=\sqrt{w_x(y,y)+v^2 + \varepsilon \sqrt{w_x(y,y)^2+v^4}}.$$
 $(\mathbb
R^{n},F_\varepsilon)$ is a Riemannian manifold if and only if
$\varepsilon=0$; however, if $\varepsilon>0$, then $(\mathbb
R^{n},F_\varepsilon)$ is a non-compact, complete, reversible {non-Riemannian
	Berwald} space with nonnegative Ricci curvature.

Fix $\varepsilon>0$. According to Corollary \ref{corollary-Finsler-1}, if equality holds in $(\ref{volume-BM-3})$ for some open sets $A$ and $B$ in $(\mathbb	R^{n},F_\varepsilon),$ then $(\mathbb R^{n},F_\varepsilon)$ is a (locally) Minkowski space, i.e., $w_x$ is independent of $x$.
%\begin{itemize}
%	\item[(i)]    equality holds in $(\ref{volume-BM-3})$ for arbitrary	geodesic balls $A$ and $B$ in $(\mathbb
%	R^{n},F_\varepsilon);$
%	\item[(ii)] $(\mathbb R^{n},F_\varepsilon)$ is a Minkowski space (i.e., $g_x$ is independent of $x$).
%\end{itemize}
\end{example}

% The reversibility of a Minkowski space can be characterized via the equality in the Brunn-Minkowski inequality: 
%The equality in (\ref{volume-BM-3}) produces surprising scenario even on Minkowski spaces. 

Minkowski spaces are the simplest non-Euclidean Finsler structures, e.g., 	 geodesics  are straight lines, the flag curvature is zero,  ${\sf{m}}_F(S)=\mathcal L^n(S)$ for every measurable set $S\subset M,$ and for every $x,y\in \mathbb R^n$ the Finslerian distance function is given by $d_F(x,y)=F(y-x)$, see Bao, Chern and Shen \cite{BCS}.   However, it turns out that the equality in the Brunn-Minkowski inequality on a generic Minkowski space $(\mathbb R^n,F)$ is not automatically verified even for forward and backward geodesic balls. In addition,  in Example \ref{example-Matsumoto} we provide two classes of Minkowski spaces where equalities  in the Brunn-Minkowski inequality generate two genuinely different scenarios. 

	\begin{corollary}\label{corollary-Finsler-24}{\bf (Brunn-Minkowski inequality  on Minkowski spaces)}
	Let $(\mathbb R^n,F)$ be a Minkowski space, $s\in (0,1)$ and $A,B\subset M$ nonempty open bounded sets. Then the inequality $(\ref{volume-BM-3})$ holds; moreover, if $A$ and $B$ are convex sets $($in the usual sense$)$,  equality holds in $(\ref{volume-BM-3})$ if and only if $A$ and $B$ are homothetic. 
	If $x,y\in \mathbb R^n$ and $r,R>0$ are fixed, the following statements are equivalent: 
	\begin{itemize}
		\item[(i)] equality holds in $(\ref{volume-BM-3})$ for $A=B^+(x,r)$ and $B=B^-(y,R)$$;$
		\item[(ii)]  $B^-(y-x_0,R)=B^+(\frac{R}{r}x,R)$ for some $x_0\in \mathbb R$.  
	\end{itemize}
\end{corollary}

	{\it Proof.} Inequality  (\ref{volume-BM-3}) trivially holds. 
 Assume that for the convex sets $A$ and $B$ we have equality in  $(\ref{volume-BM-3})$.   The positive homogeneity of $F$ implies that
 $Z_s(A,B)=(1-s)A+sB$.  
	 Accordingly, the  equality in (\ref{volume-BM-3}) can be transposed to an equality in the Euclidean Brunn-Minkowski inequality for $A$ and $B$, obtaining  that $A$ and $B$ are homothetic.
	
In the sequel,  let $A=B^+(x,r)$ and $B=B^-(y,R)$ for some $x,y\in \mathbb R^n$ and $r,R>0$.

 (i)$\Rightarrow$(ii). Assume we have equality in (\ref{volume-BM-3}) for $A$ and $B$.  Note that these sets are strictly convex domains of $\mathbb R^n$ in the usual sense, both of them inheriting the convexity of the Minkowski norm $F$, see e.g. Bao, Chern and Shen \cite[p. 12]{BCS}. Accordingly, from the first part of the proof, the sets $A$ and $B$ are homothetic, i.e., $B^-(y,R)=c_0B^+(x,r)+x_0,$ for some $c_0>0$ and $x_0\in \mathbb R^n.$ Moreover, it follows that   $c_0=\frac{R}{r}$, thus    $B^-(y-x_0,R)=B^+(\frac{R}{r}x,R)$. 
 
 (ii)$\Rightarrow$(i). Since $B^-(y-x_0,R)=B^+(\frac{R}{r}x,R),$ we have that $Z_s(A,B)=(1-s)B^+(x,r)+sB^-(y,R)=sx_0+((1-s)\frac{r}{R}+s)B^+(\frac{R}{r}x,R)$. Thus, if $\omega_n$ deotes the volume of the unit ball
 in $\mathbb R^n,$ then $${\sf{m}}_F(Z_s(A,B))^\frac{1}{n}=\left((1-s)\frac{r}{R}+s\right)R\omega_n^\frac{1}{n}=
 (1-s)r\omega_n^\frac{1}{n}+sR\omega_n^\frac{1}{n}=(1-s){\sf{m}}_F(A)^\frac{1}{n}+s{\sf{m}}_F(B)^\frac{1}{n},$$
 which concludes the proof.
\hfill $\square$\\

For simplicity, in the following example we restrict our argument to two-dimensional objects. 

\begin{example}\label{example-Matsumoto} \rm  (a) (\textit{Randers-type Minkowski plane}) 
	Let  $F_b:T\mathbb R^2\to [0,\infty)$ be defined
	by
	\begin{equation}\label{Randers-metrika}
	F_b(x,y):=F_b(y)=\sqrt{\langle Qy,y\rangle}+\langle b,y\rangle,\ \ (x,y)\in T\mathbb R^2,
	\end{equation}
	where $Q$ is a $2\times 2$ positive definite symmetric matrix, $\langle \cdot,\cdot\rangle$ is the usual scalar product in $\mathbb R^2$ and $b\in \mathbb R^2$ is fixed such that $\langle Q^{-1}b,b\rangle<1;$ here $Q^{-1}$ denotes the inverse of $Q$. The pair $(\mathbb R^2,F_b)$ is a  Randers-type
		Minkowski plane which describes  the anisotropic Luneburg-type refraction in optical crystals or the
	electromagnetic field of the physical space-time in general
	relativity (in higher dimension), see Randers \cite{Randers}. Note that $(\mathbb R^2,F_b)$  is reversible if and only if $b=(0,0).$ 
	
	Let $R,r>0$ and $x,y\in \mathbb R^2$ be arbitrarily fixed. Since the forward and backward indicatrices $I^+(x,r)=\partial B^+(x,r)$ and $I^-(y,R)=\partial B^-(y,R)$ are ellipses which can be obtained from each other by translation and dilation, equality in the Brunn-Minkowski inequality (\ref{volume-BM-3}) holds for any choice of $A=B^+(x,r)$ and $B=B^-(y,R)$ in $(\mathbb R^2,F_b)$, due to Corollary \ref{corollary-Finsler-24}; see also Figure \ref{abra-elso}(a). 
	
	\medskip
	
(b)	(\textit{Matsumoto mountain slope metric}) Let 
	$F_\alpha:T\mathbb R^2\to [0,\infty)$ be defined by
 \begin{equation}\label{Matsumoto_metrika} F_\alpha(x,y):=F_\alpha(y)=\left\{
		\begin{array}{ll}
		\frac{y_{1}^{2}+y_{2}^{2}}{v\sqrt{y_{1}^{2}+y_{2}^{2}%
			}+\frac{g}{2}y_{2}\sin\alpha}, & \hbox{}  y=(y_{1},y_{2})\in\mathbb{R}^{2}\setminus\{(0,0)\}; \\
		0, & \hbox{} y=(y_{1},y_{2})=(0,0),
		\end{array}
		\right.
		\end{equation}
	 where  $\alpha\in [0,\pi/2)$,  $v>0$ and 
	$g\approx9.81$. If we assume that $g\sin\alpha < v$, it turns out that $(\mathbb R^2,F_\alpha)$ is a Minkowski plane, describing the law of walking with a constant speed $v$$[m/s]$ under the effect of gravity on a mountain slope having the angle $\alpha$ w.r.t. the horizontal plane,  see	Matsumoto \cite{Matsumoto}.  It is clear that $(\mathbb R^2,F_\alpha)$  is reversible if and only if $\alpha=0,$ which corresponds to the Euclidean setting and $F_b$ reduces to the standard (reversible) metric $F_0(y_1,y_2)=\sqrt{y_1^2+y_2^2}/v.$

	Let   $x,y\in \mathbb R^n$ and $r,R>0$ be arbitrarily fixed. We notice that the indicatrices $$I^+(x,r)=\partial A=\{z\in \mathbb R^2: F_\alpha(z-x)=r\}\ \ {\rm and}\ \ I^-(y,R)=\partial B=\{z\in \mathbb R^2: F_\alpha(y-z)=R\}$$  are convex lima\c cons which cannot be obtained from each other by dilation and translation, unless $\alpha=0$ (i.e., the mountain slope vanishes), see also Figure \ref{abra-elso}(b).  Thus,  due to Corollary \ref{corollary-Finsler-24}, for $\alpha\neq 0$ (i.e.,  we are in the non-Euclidean setting) any choice of  $A=B^+(x,r)$ and $B=B^-(y,R)$ in $(\mathbb R^2,F_\alpha)$ provides strict inequality in the Brunn-Minkowski inequality (\ref{volume-BM-3}). 
	
\begin{figure}[H]
	\includegraphics
	{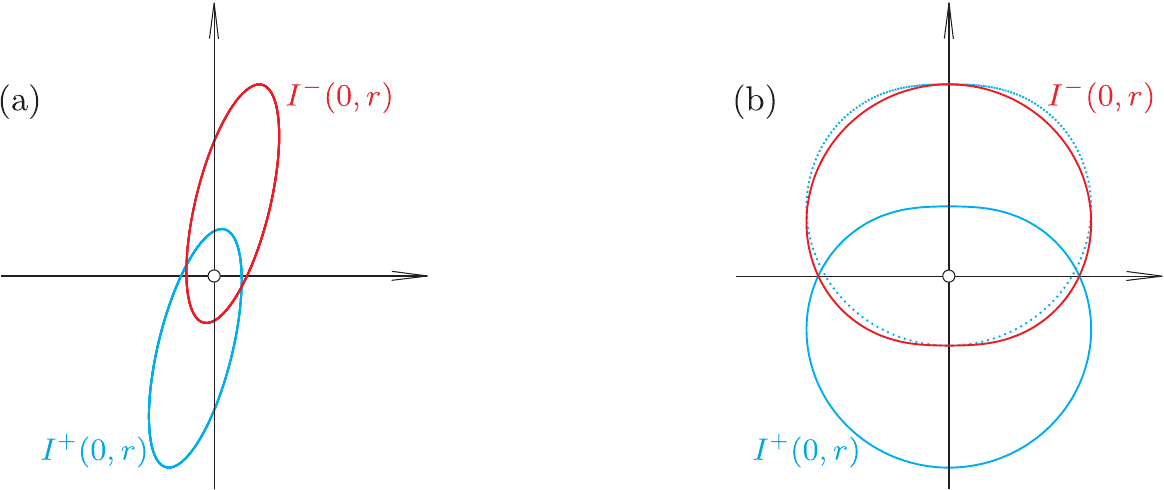}
	\caption{(a) Let $Q=[5,~-1;-1,~1]$ and $b=(1/5,1/2)$ in Example \ref{example-Matsumoto} (a). The forward and backward balls (which are ellipses) with the same radius can be translated to each other. (b) Let $\alpha\approx 35^o$ and $v=6$ in Example \ref{example-Matsumoto} (b). The forward and backward balls with the same radius cannot be translated to each other.}\label{abra-elso}
\end{figure}
\end{example}

\section{Appendix: Proof of Lemma \ref{lemma-p-mean}}\label{appendix}

\noindent 	We first recall  the quantitative Young inequality (see e.g. Cianchi \cite{Cianchi}), i.e.,  if $r\geq 2$ and $\frac{1}{r}+\frac{1}{r'}=1$, one has
\begin{equation}\label{Cianchi}
uv\leq \frac{1}{r}u^r+\frac{1}{r'}v^{r'}-\frac{1}{r}|u-v^\frac{1}{r-1}|^r\ \ {\rm for\ every}\ u,v\geq 0.
\end{equation}

{(i)} Let $p\in (0,\infty)$ and let us assume first that $pn\geq 1.$    Applying inequality (\ref{Cianchi}) for
$r=\frac{p}{\tilde p}= pn+1\geq 2$ and $r'=\frac{1}{\tilde p n}$, we
have that
\begin{eqnarray*}
	% \nonumber to remove numbering (before each equation)
	\frac{\left[\mathcal M_s^{\tilde p}(\frac{1}{c},\frac{1}{d})\right]^{\tilde p}}{\left[\mathcal M_s^p(a,b)\right]^{\tilde p}
		\left[\mathcal M_s^{\frac{1}{n}}(\frac{1}{ac},\frac{1}{bd})\right]^{\tilde p}} &=&
	\frac{(1-s)^\frac{1}{r}a^{\tilde p}}{\left[\mathcal M_s^p(a,b)\right]^{\tilde p}}\frac{(1-s)^\frac{1}{r'}(\frac{1}{ac})^{\tilde p}}{\left[\mathcal M_s^{\frac{1}{n}}(\frac{1}{ac},\frac{1}{bd})\right]^{\tilde p}}+
	\frac{s^\frac{1}{r}b^{\tilde p}}{\left[\mathcal M_s^p(a,b)\right]^{\tilde p}}\frac{s^\frac{1}{r'}(\frac{1}{bd})^{\tilde p}}{\left[\mathcal M_s^{\frac{1}{n}}(\frac{1}{ac},\frac{1}{bd})\right]^{\tilde p}} \\
	&\leq& \frac{1}{r} \frac{(1-s)a^{\tilde p r}}{\left[\mathcal M_s^p(a,b)\right]^{\tilde p r}}+\frac{1}{r'} \frac{(1-s)(\frac{1}{ac})^{\tilde p r}}{\left[\mathcal M_s^{\frac{1}{n}}(\frac{1}{ac},\frac{1}{bd})\right]^{\tilde p r}}  \\&&-\frac{1}{r}\left|\frac{(1-s)^\frac{1}{r}a^{\tilde p}}{\left[\mathcal M_s^p(a,b)\right]^{\tilde p}}-\left(\frac{(1-s)^\frac{1}{r'}(\frac{1}{ac})^{\tilde p}}{\left[\mathcal M_s^{\frac{1}{n}}(\frac{1}{ac},\frac{1}{bd})\right]^{\tilde p}}\right)^\frac{1}{r-1}\right|^r\\
	&&+\frac{1}{r} \frac{sb^{\tilde p r}}{\left[\mathcal M_s^p(a,b)\right]^{\tilde p r}}+\frac{1}{r'} \frac{s(\frac{1}{bd})^{\tilde p r}}{\left[\mathcal M_s^{\frac{1}{n}}(\frac{1}{ac},\frac{1}{bd})\right]^{\tilde p r}}  \\&&-\frac{1}{r}\left|\frac{s^\frac{1}{r}b^{\tilde p}}{\left[\mathcal M_s^p(a,b)\right]^{\tilde p}}-\left(\frac{s^\frac{1}{r'}(\frac{1}{bd})^{\tilde p}}{\left[\mathcal M_s^{\frac{1}{n}}(\frac{1}{ac},\frac{1}{bd})\right]^{\tilde p}}\right)^\frac{1}{r-1}\right|^r\\
	&=&
	1-\frac{1}{r}(1-s)\left|\left[\mathcal M_s^{-p}\left(1,\frac{a}{b}\right)\right]^{\tilde p}-\left[\mathcal M_s^{-\frac{1}{n}}\left(1,\frac{bd}{ac}\right)\right]^\frac{\tilde p}{pn}\right|^r\\
	&&-\frac{1}{r}s\left|\left[\mathcal M_s^{-p}\left(\frac{b}{a},1\right)\right]^{\tilde p}-\left[\mathcal M_s^{-\frac{1}{n}}\left(\frac{ac}{bd},1\right)\right]^\frac{\tilde p}{pn}\right|^r.
\end{eqnarray*}
By rearranging the latter estimate, we obtain
\begin{eqnarray*}
	\left[\mathcal M_s^p(a,b)\right]^{\tilde p}
	\left[\mathcal M_s^{\frac{1}{n}}\left(\frac{1}{ac},\frac{1}{bd}\right)\right]^{\tilde p}&\geq& \left[\mathcal M_s^{\tilde p}\left(\frac{1}{c},\frac{1}{d}\right)\right]^{\tilde p}+ \frac{1}{r}\left[\mathcal M_s^p(a,b)\right]^{\tilde p}
	\left[\mathcal M_s^{\frac{1}{n}}\left(\frac{1}{ac},\frac{1}{bd}\right)\right]^{\tilde p}\times\\
	&&\times \left[(1-s)\left|\left[\mathcal M_s^{-p}\left(1,\frac{a}{b}\right)\right]^{\tilde p}-\left[\mathcal M_s^{-\frac{1}{n}}\left(1,\frac{bd}{ac}\right)\right]^\frac{\tilde p}{pn}\right|^r\right.+ \\
	&&   \ \ \ \ \ \  +s \left.\left|\left[\mathcal M_s^{-p}\left(\frac{b}{a},1\right)\right]^{\tilde p}-\left[\mathcal M_s^{-\frac{1}{n}}\left(\frac{ac}{bd},1\right)\right]^\frac{\tilde p}{pn}\right|^r\right].
\end{eqnarray*}
Since $\frac{1}{p} +n = \frac{1}{\tilde p}$ we can apply (\ref{MspIneq}) to get
$$\mathcal M_s^p(a,b)
\mathcal M_s^{\frac{1}{n}}\left(\frac{1}{ac},\frac{1}{bd}\right)\geq \mathcal M_s^{\tilde p}\left(\frac{1}{c},\frac{1}{d}\right).$$
Using this estimate on the right hand side of the above inequality and the definition of  $G_s^{p,n}(a,b,c,d)$, 
%	combined with the elementary inequality
%	$$ |A^{\alpha} -B^{\alpha}| ^\frac{1}{\alpha} \geq |A-B| \ \text{for} \ A, B \geq 0 \ \text{and} \  \alpha \geq 1,$$
it follows that 
\begin{eqnarray*}
	\left[\mathcal M_s^p(a,b)\right]^{\tilde p}
	\left[\mathcal M_s^{\frac{1}{n}}\left(\frac{1}{ac},\frac{1}{bd}\right)\right]^{\tilde p}&\geq& \left[\mathcal M_s^{\tilde p}\left(\frac{1}{c},\frac{1}{d}\right)\right]^{\tilde p}\left(1+\frac{p}{r} G_s^{p,n}(a,b,c,d)\right).
\end{eqnarray*}
Note that $ \frac{1}{
	\tilde p}=\frac{pn+1}{p}> 1;$ then we may apply to the latter estimate Bernoulli's inequality   (i.e., $(1+z)^r\geq 1+rz$ for $r\geq 1$ and $z\geq -1$), obtaining
\begin{eqnarray*}
	\mathcal M_s^p(a,b)
	\mathcal M_s^{\frac{1}{n}}\left(\frac{1}{ac},\frac{1}{bd}\right)&\geq& \mathcal M_s^{\tilde p}\left(\frac{1}{c},\frac{1}{d}\right)\left(1+\frac{p}{r\tilde p} G_s^{p,n}(a,b,c,d)\right).
\end{eqnarray*}
Since $r\tilde p =p$ and  $$\mathcal M_s^{\tilde p}\left(\frac{1}{c},\frac{1}{d}\right)=\left[\mathcal M_s^{-\tilde p}(c,d)\right]^{-1},$$
the desired relation {follows}. The  case $pn\leq 1$ follows in the same way. 

If $p\in (-\frac{1}{n},0).$ Since $-\tilde p=-\frac{p}{pn+1}\in (0,\infty)$ and $ p=\frac{\tilde p}{\tilde pn+1}$, we can apply the previous estimate by reversing the roles of the means. 

A simple computation shows that $G_s^{p,n}(a,b,c,d)=0$ if and only if $\frac{a}{b}=\left(\frac{d}{c}\right)^\frac{1}{pn+1}.$

{(ii)}   We first assume that $s\geq \frac{1}{2}$, i.e., $\min(s,1-s)=1-s$. We apply (\ref{Cianchi}) with $r=\frac{1}{1-s}\geq 2$ and $r'=\frac{1}{s}$, obtaining 
\begin{eqnarray*}
	% \nonumber to remove numbering (before each equation)
	\frac{[\mathcal M_s^0(\frac{1}{c},\frac{1}{d})]^\frac{1}{n}}{\left[\mathcal M_s^0(a,b)\right]^\frac{1}{n}
		[\mathcal M_s^\frac{1}{n}(\frac{1}{ac},\frac{1}{bd})]^\frac{1}{n}} &=&
	\frac{(\frac{1}{ac})^{\frac{1-s}{n}}(\frac{1}{bd})^{\frac{s}{n}}}{(1-s)(\frac{1}{ac})^\frac{1}{n}+s(\frac{1}{bd})^\frac{1}{n}} \\
	&\leq& \frac{(1-s)(\frac{1}{ac})^\frac{1}{n}+s(\frac{1}{bd})^\frac{1}{n}-(1-s)|(\frac{1}{ac})^{\frac{1-s}{n}}-((\frac{1}{bd})^{\frac{s}{n}})^\frac{1-s}{s}|^\frac{1}{1-s}}{(1-s)(\frac{1}{ac})^\frac{1}{n}+s(\frac{1}{bd})^\frac{1}{n}}\\
	&=&1-(1-s)\frac{|(\frac{1}{ac})^{\frac{1-s}{n}}-(\frac{1}{bd})^{\frac{1-s}{n}}|^\frac{1}{1-s}}{(1-s)(\frac{1}{ac})^\frac{1}{n}+s(\frac{1}{bd})^\frac{1}{n}}\\ &=& 1-(1-s)\left|\left[{\mathcal M_s^{-\frac{1}{n}}\left(1,\frac{bd}{ac}\right)}\right]^{\frac{1-s}{n}}-\left[{\mathcal M_s^{-\frac{1}{n}}\left(\frac{ac}{bd},1\right)}\right]^{\frac{1-s}{n}}\right|^\frac{1}{1-s}.
\end{eqnarray*}
Rearranging the above inequality, and using 
$ \mathcal M_s^0(ac,bd) \geq \mathcal M_s^{-\frac{1}{n}}(ac,bd)$
and the Bernoulli inequality, it follows that 
$$\mathcal M_s^0(ac,bd)
\geq \mathcal M_s^{-\frac{1}{n}}(ac,bd)\left(1+n{(1-s)}\left|\left[{\mathcal M_s^{-\frac{1}{n}}\left(1,\frac{bd}{ac}\right)}\right]^{\frac{1-s}{n}}-\left[{\mathcal M_s^{-\frac{1}{n}}\left(\frac{ac}{bd},1\right )}\right]^{\frac{1-s}{n}}\right|^\frac{1}{1-s}\right).$$ If $s\leq \frac{1}{2}$, we proceed in a similar way as above. Furthermore, $G_s^{0,n}(a,b,c,d)=0$ if and only if $ac=bd.$

{(iii)} We first assume that $a\geq b$. Then we have
\begin{eqnarray*}
	% \nonumber to remove numbering (before each equation)
	\frac{[\mathcal M_s^\frac{1}{n}(\frac{1}{c},\frac{1}{d})]^\frac{1}{n}}{[\mathcal M_s^{+\infty}(a,b)]^\frac{1}{n}
		[\mathcal M_s^{\frac{1}{n}}(\frac{1}{ac},\frac{1}{bd})]^\frac{1}{n}} &=&1+s\frac{(b^\frac{1}{n}-a^\frac{1}{n})}{(abd)^\frac{1}{n}}\left[\mathcal M_s^{-\frac{1}{n}}(ac,bd)\right]^\frac{1}{n}.
\end{eqnarray*}
After a rearrangement,  Bernoulli's inequality and (\ref{MspIneq}) give the required inequality. The same can be done for $a\leq b.$ Moreover, $G_s^{+\infty,n}(a,b,c,d)=0$ if and only if $a=b.$

The proof of (iv) directly follows by (iii); we left it to the interested reader. 
%
%{\bf Case 3:} $p=+\infty$. We assume that $a\geq b.$ Since $\eta=q$ and $$\frac{[\mathcal M_s^q(ac,bd)]^{q}}{[\mathcal M_s^{+\infty}(a,b)]^q
%	[\mathcal M_s^q(c,d)]^q} =1+s\frac{d^q}{	[\mathcal M_s^q(c,d)]^q}\left(\frac{b^q}{a^q}-1\right),$$
%a reorganization of this relation, (\ref{MspIneq}) and the Bernoulli inequality give that
%$$\mathcal M_s^{+\infty}(a,b)
%\mathcal M_s^q(c,d)\geq \mathcal M_s^q(ac,bd)\left(1+\frac{s}{q}\frac{d^q}{[\mathcal M_s^q(c,d)]^q}\frac{a^q-b^q}{a^q}\right).$$
%The case  $a\leq b$ is similar. 
%
%{\bf Case 4:} $p\in (-q,0)$. We have that $\eta<0$ and
%$-p= \frac{(-\eta)q}{(-\eta)+q}\leq 1.$
%We are in the position to apply Case 1 with
%the choices $p:=-\eta>0$, $q:=q$ and $\eta:=-p\leq 1,$ respectively,
%obtaining for every $a,b,c,d>0$ that
%$$ \mathcal M_s^{-\eta}(a,b)\mathcal M_s^q(c,d)\geq \mathcal M_s^{-p}(ac,bd)\left(1+G_s^{-\eta,q}(a,b,c,d) \right).$$
%Since $\mathcal M_s^{-\eta}(a,b)=[\mathcal M_s^{\eta}(1/a,1/b)]^{-1}$, we may rearrange the latter inequality into 
%$$\mathcal M_s^{p}(1/(ac),1/(bd))\mathcal M_s^q(c,d)\geq \mathcal M_s^{\eta}(1/a,1/b)\left(1+G_s^{-\eta,q}(a,b,c,d) \right),$$
%which is the claimed relation. 
\hfill$\square$

\vspace{1cm}

\noindent {\bf Acknowledgement.}  
We express our gratitude to Alessio Figalli and C\'edric Villani  for their useful comments on the earlier version of the manuscript. 
%We express our gratitude to Alessio Figalli  for his suggestion to consider Dubuc's result on curves spaces  
%in the earlier version of the manuscript. 
A. Krist\'aly is grateful to the Mathematisches Institute of Bern for the warm hospitality where this work has been initiated. We also 
thank
the anonymous 
Referee
for
her/his valuable
comments that greatly improved the presentation of the manuscript.\\

%	\begin{itemize}
%	\item[(i)] there exists $x_0\in \mathbb R^n$ such that up to null measure sets  $$\operatorname{supp}g=c_0\operatorname{supp}f +x_0\ \  \ {\rm and}\ \  \operatorname{supp}h=(1-s+sc_0)\operatorname{supp}f +sx_0, $$
%	where $c_0=\left(\frac{\mathcal L^n(\operatorname{supp}g)}{\mathcal L^n(\operatorname{supp}g)}\right)^\frac{1}{n}$. 
%	\item[(ii)] there exist a convex set $S\subset \mathbb R^n$ with $S=\operatorname{supp}f$ up to a null measure set and a $(t,p)-$concave function $\Phi:S\to \mathbb R$ with $t=\frac{sc_0}{1-s+sc_0}$ such that for a.e. $x\in S$ we have 
%	$$f(x)=\Phi(x);$$
%	$$g(c_0x+x_0)=c_0^\frac{1}{p}f(x);$$
%	$$h((1-s+sc_0)x +sx_0)=\left[\mathcal  M_s^\frac{p}{pn+1}\left(1,c_0^\frac{pn+1}{p}\right)\right]^\frac{1}{pn+1}f(x).$$
%\end{itemize}

%	\vspace{0.5cm}

\vspace{0.5cm} 

\noindent {\footnotesize{\sc  Mathematisches Institute,
		Universit\"at Bern,
		Sidlerstrasse 5,
		3012 Bern, Switzerland.}\\
	Email: \textsf{zoltan.balogh@math.unibe.ch}\\

	\noindent {\footnotesize {\sc Department of Economics, Babe\c s-Bolyai University, Str. T. Mihali 58-60, 400591
			Cluj-Napoca, Romania;  \\ Institute of Applied Mathematics, \'Obuda University,
			B\'ecsi \'ut 96, 1034 Budapest, Hungary.}\\ Email:
		{\textsf{alex.kristaly@econ.ubbcluj.ro; kristaly.alexandru@nik.uni-obuda.hu}}\\


\begin{thebibliography}{99}
		
%		\bibitem{ADM} S. Alesker, S. Dar,  V. Milman, 
%		{\it A remarkable measure preserving diffeomorphism between two convex bodies in $\mathbb R^n,$} 
%		Geom. Dedicata 74 (1999), no. 2, 201--212.
		
		
		\bibitem{Bacher} K. Bacher, \textit{On Borell-Brascamp-Lieb inequalities on metric measure spaces},
		Potential Anal.
		 33 (2010), no. 1, 1--15. 
		
		
			\bibitem{BB-1}	K. M.  Ball, K. J.  B\"or\"oczky, {\it  Stability of some versions of the Pr\'ekopa-Leindler inequality}, Monatsh. Math. 163 (2011), no. 1, 1--14.
		
	\bibitem{BB-2}	K. M.  Ball, K. J.  B\"or\"oczky, {\it  Stability of the Pr\'ekopa-Leindler inequality}, Mathematika 56 (2010), no. 2, 339--356.
	
	
		
%		\bibitem{BKS-1} Z. M. Balogh, A. Krist\'aly, K. Sipos, {\it Geodesic interpolation inequalities on Heisenberg groups,} C. R. Math. Acad. Sci. Paris 354 (2016), no. 9, 916--919. 
%		
%			\bibitem{BKS-2} Z. M. Balogh, A. Krist\'aly, K. Sipos,  {\it Geometric inequalities on Heisenberg groups,}  arxiv.org/abs/1605.06839.
%			
%			\bibitem{BKS-3} Z. M. Balogh, A. Krist\'aly, K. Sipos, {\it Jacobian determinant inequality on corank $1$ Carnot groups with applications,} arxiv.org/abs/1701.08831. 
			
			\bibitem{BCS} D.~Bao, S.~S.~Chern, Z.~Shen, {\it Introduction to
				Riemann--Finsler Geometry,} Graduate Texts in Mathematics, 200,
			Springer Verlag, 2000.
		
		\bibitem{Bishop-Crittenden} R. L.  Bishop, R. J.  Crittenden,  {\it Geometry of manifolds.} Reprint of the 1964 original. AMS Chelsea Publishing, Providence, RI, 2001. 
		
		
	\bibitem{Bobkov-Ledoux} S. G. 	Bobkov, M. Ledoux, 
		\textit{From Brunn-Minkowski to Brascamp-Lieb and to logarithmic Sobolev inequalities}, 
		Geom. Funct. Anal. 10 (2000), no. 5, 1028--1052. 
		
		\bibitem{Boothby} W. M. Boothby,  \textit{An introduction to differentiable manifolds and Riemannian geometry.} Second edition. Pure and Applied Mathematics, 120. Academic Press, Inc., Orlando, FL, 1986.
		
		\bibitem{Borell} C. Borell, \textit{Convex set functions in d-space}, Period. Math. Hung. 6 (1975), 111--136. 
		
		\bibitem{BLYZ}  K. J. B\"or\"oczky,  E. Lutwak, D. Yang, G. Zhang,  {\it The log-Brunn-Minkowski inequality}, Adv. Math. 231 (2012), no. 3-4, 1974--1997.
		
		\bibitem{Brascamp-Lieb} H.J. Brascamp, E.H. Lieb, \textit{On extensions of the Brunn-Minkowski and Prékopa-Leindler theorems, including inequalities for log concave functions and with an application to the diffusion equation}, J. Funct. Anal. 22 (1976), no. 4, 366--389. 
		
		\bibitem{BF} D. Bucur, I. Fragal\`a, 
		\textit{Lower bounds for the Pr\'ekopa-Leindler deficit by some distances modulo translations},
		J. Convex Anal. 21 (2014), no. 1, 289--305. 
		
		\bibitem{Caffarelli} L. Caffarelli, {\it The regularity of mappings with a convex potential},
		J. Amer. Math. Soc. 5 (1992), no. 1, 99--104.
		
		\bibitem{Christ} M. Christ, \textit{Near equality in the Brunn-Minkowski inequality,} preprint. Available online at: https://arxiv.org/abs/1207.5062.
		
		\bibitem{Cianchi} A. Cianchi,  {\it Sharp Morrey-Sobolev inequalities and the distance from extremals}, Trans. Amer. Math. Soc. 360 (2008), no. 8, 4335--4347.
		
	\bibitem{CLM} A. Colesanti, G. Livshyts, A. Marsiglietti, 	\textit{On the stability of Brunn-Minkowski type inequalities}, J. Funct. Analysis,  273 (2017), no. 3, 1120--1139.
	
		
		\bibitem{Cordero-CRAS}
		D.~Cordero-Erausquin,
		{\it In\'egalit\'e de Pr\'ekopa-Leindler sur la sph\`ere}, C. R. Acad. Sci. Paris S\'er. I Math. 329 (1999), no. 9, 789--792.
		
		\bibitem{CMS}
		D.~Cordero-Erausquin, R.~J. McCann,  M.~Schmuckenschl{\"a}ger, \emph{A
			{R}iemannian interpolation inequality \`a la {B}orell, {B}rascamp and
			{L}ieb}, Invent. Math. {146} (2001), no.~2, 219--257.
		
		
		\bibitem{DU} S. Dancs, B. Uhrin, {\it On a class of integral inequalities and their measure-theoretic consequences}, J. Math. Anal. Appl. 74 (1980), 388--400.
		
		\bibitem{DU-2} S. Dancs, B. Uhrin, {\it On the conditions of equality in an integral inequality}, Publ. Math. (Debrecen)
		29 		(1982), 117--132.
		
		\bibitem{doCarmo} {M. P. do Carmo}, \textit{{Riemannian Geometry}}, Birkh\"auser,
	Boston, 1992.
		
		\bibitem{Dubuc} S. Dubuc, {\it Crit\`eres de convexit\'e et in\'egalit\'es int\`egrales},   Ann. Inst.
		Fourier (Grenoble) 27 (1977) no. 1,  135--165.
		
	
	
		
	\bibitem{Figalli-Jerison-2} A. Figalli, D. Jerison, 	{\it Quantitative stability of the Brunn-Minkowski inequality for sets of equal volume},
		 Chin. Ann. Math. Ser. B 38 (2017), no. 2, 393--412.
	
		\bibitem{Figalli-Jerison} A. Figalli, D. Jerison,  {\it Quantitative stability for the
		Brunn-Minkowski inequality}, Adv. Math.,  314 (2017),  1--47.	
	
\bibitem{Figalli-Jerison-3}	A. Figalli, D. Jerison, \textit{A sharp Freiman type estimate for semisums in $\mathbb R^n$}, preprint. Available on-line at: 
https://people.math.ethz.ch/~afigalli/submitted-pdf/A-sharp-freiman-type-estimate-for-semisums-in-rn.pdf


\bibitem{FMP-Inventiones}	A. Figalli,	F. Maggi, A. Pratelli, \textit{A mass transportation approach to quantitative isoperimetric inequalities}, 
Invent. Math. 182 (2010), no. 1, 167--211.

\bibitem{FMP-AIHP}	A. Figalli, F. Maggi, A. Pratelli,  \textit{A refined Brunn-Minkowski inequality for convex sets},  Ann. Inst. H. Poincar\'e Anal. Non Lin\'eaire 26 (2009), no. 6, 2511--2519.
		
		\bibitem{Gardner}
		R.~J. Gardner, \emph{The {B}runn-{M}inkowski inequality}, Bull. Amer. Math.
		Soc. (N.S.) {39} (2002), no.~3, 355--405.
		
		\bibitem{GS} D. Ghilli, P. Salani, \textit{Quantitative Borell-Brascamp-Lieb inequalities for power concave functions}, J. Convex Anal. 24 (2017), No. 3. Available online at: http://www.heldermann.de/JCA/JCA24/jca24.htm
		
	\bibitem{Knothe}  H. Knothe,  \textit{Contributions to the theory of convex bodies}, Michigan Math. J. 4 (1957), 39--52.
		
		
		\bibitem{Kristaly-JGA} A. Krist\'aly, {\it A sharp Sobolev interpolation inequality on Finsler manifolds}, J. Geom. Anal. 25 (2015), no. 4, 2226--2240.
		
		 \bibitem{LV} J. Lott, C. Villani, {\it Ricci curvature for metric-measure spaces via optimal transport}, Ann. of
		Math. (2) {169} (2009), no. 3, 903--991. 
		
		\bibitem{Matsumoto} M. Matsumoto, {\it A slope of a mountain is a Finsler surface with respect
			to a time measure}, {J. Math. Kyoto Univ.} 29 (1989), 17--25.
		
		\bibitem{McCann} R. J. McCann, {\it Polar factorization of maps on Riemannian manifolds},
		Geom. Funct. Anal. 11 (2001), no. 3, 589--608.
		
			\bibitem{McCann_Adv_Math} R. J.  McCann, \textit{A convexity principle for interacting gases}, Adv. Math. 128 (1997), no. 1, 153--179.
			
		\bibitem{McCann-PhD}  R. J. McCann, \textit{A convexity theory for interacting gases and equilibrium crystals}, ProQuest LLC,
		Ann Arbor, MI, 1994, Thesis (Ph.D.) Princeton University.
		
			\bibitem{MR} E. Milman, L. Rotem,\textit{ Complemented Brunn-Minkowski inequalities and isoperimetry for homogeneous and non-homogeneous measures}, Adv. Math. 262 (2014), 867--908.
		
		 \bibitem{MS} V. D. Milman, G. Schechtman,  \textit{Asymptotic theory of finite-dimensional normed spaces}, With an appendix by M. Gromov. Springer-Verlag, Berlin, 1986.
		
		
		\bibitem{Ohta} S. Ohta, {\it Finsler interpolation inequalities}, Calc. Var. Partial Differential Equations 36 (2009), no. 2, 211--249.
		
		\bibitem{Randers} G. Randers, \textit{On an asymmetrical metric in the fourspace of general relativity},  Phys. Rev. (2)  59,  (1941), 195--199.
		
		\bibitem{Rajala} T. Rajala, {\it Large porosity and dimension of sets in metric spaces,} Ann. Acad. Sci. Fenn. Math. 34 (2009), no. 2, 565--581.
		
	\bibitem{Rossi-PhD}	A. Rossi, \textit{Borell-Brascamp-Lieb inequalities: rigidity and stability}, PhD Thesis in Mathematics, Informatics, Statistics - Universit\`a di Firenze, 2018.
				
		\bibitem{Rossi-Salani} A. Rossi, P. Salani, {\it Stability for Borell-Brascamp-Lieb inequalities},   Geometric aspects of functional analysis, 339--363,
		Lecture Notes in Math., 2169, Springer, Cham, 2017. 
		
	\bibitem{Rossi-Salani-AA} A. Rossi, P. Salani,	\textit{Stability for a strengthened Borell-Brascamp-Lieb inequality}, Appl. Anal., 2018, to appear. DOI: 10.1080/00036811.2018.1451645. 
		
		\bibitem{Sakai} T. Sakai, \textit{Riemannian geometry}.  Translations of Mathematical Monographs, 149. American Mathematical Society, Providence, RI, 1996.
		
		\bibitem{Shen} Z. Shen, {\it Volume comparison and its applications in Riemann-Finsler geometry}, Adv. Math. 128 (1997), no. 2, 306--328. 
		
		\bibitem{Steinhaus} H. Steinhaus, \textit{Sur les distances des points dans les ensembles de mesure positive}, Fund. Math.  1 (1920), 93--104.
		
		\bibitem{Sturm-2} K.-T. Sturm, {\it On the geometry of metric measure spaces. II}, Acta Math. {196} (2006), no. 1, 133--177.
		
		
		\bibitem{Villani-1} C. Villani, 	{\it  Optimal transport, Old and new.} Grundlehren der Mathematischen Wissenschaften [Fundamental Principles
		of Mathematical Sciences], vol. 338, Springer-Verlag, Berlin, 2009. 
		
		\bibitem{Villani-2} C. Villani, {\it Topics in optimal transportation.} Graduate Studies in Mathematics, 58. American Mathematical Society, Providence, RI, 2003.
		
	\end{thebibliography}
\end{document}